\documentclass[12pt]{article}

\usepackage{graphicx}
\usepackage{amsmath,amsthm,amssymb,amsfonts,color}

\makeatletter
\g@addto@macro\normalsize{%
  \setlength\abovedisplayskip{8pt plus 3pt minus 3pt}
  \setlength\belowdisplayskip{8pt plus 3pt minus 3pt}
  \setlength\abovedisplayshortskip{6pt plus 3pt minus 2pt}
  \setlength\belowdisplayshortskip{6pt plus 3pt minus 2pt}
}
\makeatother

\date{\today}

\def\leq{\leqslant}
\def\geq{\geqslant}
\def\le{\leqslant}
\def\ge{\geqslant}

\normalsize

\numberwithin{equation}{section}

\def\({\bigl(}
\def\){\bigr)}

\newtheorem{theorem}{Theorem}
\newtheorem{lemma}{Lemma}

\newtheorem{remark}{Remark}
\newtheorem{corollary}{Corollary}
\newtheorem{definition}{Definition}
\newtheorem{algorithm}{Algorithm}

\def\dfrac#1#2{\lower0.15ex\hbox{\large$\textstyle\frac{#1}{#2}$}}
\def\({\bigl(}
\def\){\bigr)}

\def\P{\boldsymbol{P}}

\renewcommand{\P}[1]{\Pr\left(#1\right)}

\usepackage{graphicx}
\usepackage{subfigure}
\usepackage[latin1]{inputenc}
\newcommand{\Ii}{ \mathbf{1} }

\def\P{\mathrm P}

\begin{document}
\title{
THRESHOLD SELECTION FOR EXTREMAL INDEX ESTIMATION
} 
\author{Natalia M. Markovich and Igor V. Rodionov\footnote{
              V.A.Trapeznikov Institute of Control Sciences,
                Russian Academy of Sciences,
                Profsoyuznaya 65, 
                117997 Moscow, Russia;
              e-mail: markovic@ipu.rssi.ru, vecsell@gmail.com           
}}


\maketitle

\begin{abstract}
We propose a new threshold selection method for the nonparametric estimation of the extremal index of stochastic processes.
The so-called discrepancy method was proposed 
as a data-driven smoothing tool for estimation of a probability density function. Now it  is modified to select a threshold  parameter of an extremal index estimator. 
To this end, a specific normalization of the discrepancy statistic based on the Cram\'{e}r-von Mises-Smirnov statistic $\omega^2$ is calculated by the $k$ largest order statistics instead of an entire sample. Its asymptotic distribution as $k\to\infty$ is proved to be the same as the $\omega^2$-distribution. The quantiles of the latter distribution are used as discrepancy values. The rate of convergence of an extremal index estimate coupled with the discrepancy method is derived.  The discrepancy method is used as an automatic threshold selection for the intervals and $K-$gaps estimators and it may be applied to other estimators of the extremal index. \\
\\
\textbf{Keywords}:
Threshold selection; Discrepancy method;  Cram\'{e}r-von Mises-Smirnov statistic; Nonparametric estimation;  Extremal index

\end{abstract}
\section{Introduction}\label{intro}
Let $X^n= \{X_i\}_{i=1}^n$ 
be a sample of 
random variables (r.v.s) with cumulative distribution function (cdf) $F(x)$.
By Leadbetter et al. 1983 the stationary sequence  $\{X_n\}_{n\ge 1}$ is said to have extremal index
$\theta\in (0,1]$
 if
for each $0<\tau <\infty$ there is a sequence of real numbers $u_n=u_n(\tau)$ such that it holds
\begin{eqnarray*}\label{1}&&\lim_{n\to\infty}n(1-F(u_n))=\tau, \qquad
\lim_{n\to\infty}P\{M_n\le u_n\}=e^{-\tau\theta},\end{eqnarray*}
where $M_n=\max\{X_1,...,X_n\}$.
The extremal index reflects a cluster structure of an underlying sequence or its local dependence. $\theta=1$ holds if  $X_1,...,X_n$ are independent. For stationary sequences $\theta = 1$ when mixing conditions $D(u_n)$ and $D'(u_n)$ hold (Leadbetter et al. 1983).
\\
Nonparametric estimators of $\theta$ require usually the choice of a threshold and/or a declustering parameter. The well-known  blocks and runs 
estimators of the extremal index  require an appropriate threshold
$u$ and the block size $b$ or the number of consecutive observations $r$ running below $u$ to separate two clusters (Beirlant et al. 2004). A bias-corrected modification of the blocks estimator (Drees 2011)
informs how to avoid the threshold selection by providing  a rather stable plot of the extremal index estimates against $u$ with some remaining uncertainty. In Sun and  Samorodnitsky (2018) the multilevel blocks estimator is proposed where a sequence of increasing levels and a weight function have to be defined. The sliding blocks estimator 
 has asymptotic variance smaller than the disjoint blocks estimator (Robert et al. 2009b) and all of them require the selection of a pair $(u,b)$. The cycles estimator proposed by Ferreira (2018) needs both $u$ and the cycle size $s$ as parameters. The intervals estimator of $\theta$ by Ferro and Segers (2003) and the estimators introduced by Robert (2009b) require the choice of  $u$.  The $K$-gaps estimator is another threshold-based one (S\"{u}veges and Davison 2010).
\\
One of the high quantiles of the sample $X^n$ is  taken usually as $u$ or $u$ is selected visually corresponding to a stability interval of the plot of some estimate $\widehat{\theta}(u)$ against $u$. Following S\"{u}veges and Davison (2010), a list of pairs $(u,K)$ is selected   according to the Information Matrix Test (IMT) in Fukutome et al. (2015). Then $u$ is selected from such a pair that corresponds to the largest number of clusters of exceedances separated by more than $K$ non-exceedances. The semiparametric maxima estimators (Berghaus and B\"{u}cher 2018; Northrop 2015)  depend on the block size only. The choice of the latter remains an open problem.
\\
The objective of this paper is to propose a new nonparametric tool to find the threshold $u$. 
\\
The so-called discrepancy method was proposed in Markovich (1989) and Vapnik et al. (1992) as a data-driven smoothing tool for a probability density function  (pdf) estimation by i.i.d. data.
 We aim to extend this method for an extremal index estimation. 
 The idea 
 was to find an unknown parameter $h$ of the pdf as a solution of the discrepancy equation
\begin{eqnarray*}\label{DiscEq}\rho(\widehat{F}_h,F_n)&=&\delta.\end{eqnarray*}
Here,
 $\widehat{F}_h(x)=\int_{-\infty}^x\widehat{f}_h(t)dt$ holds,  $\hat{f}_h(t)$ is some pdf estimate,
  $\delta$ is a 
  discrepancy value of the estimation of  $F(x)$ by the empirical distribution function $F_n(x)$, i.e.
$\delta=\rho(F,F_n)$.  $\rho(\cdot,\cdot)$ is a metric in the space of cdf's. Since $\delta$ is usually unknown, 
quantiles of the limit distribution of the Cram\'{e}r-von Mises-Smirnov (C-M-S) statistic \begin{eqnarray*}\label{8}\omega^2_n&=&n\int_{-\infty}^{\infty}\left(F_n(x)-F(x)\right)^2dF(x),\end{eqnarray*}
were proposed as $\delta$. The latter limit distribution which is rather complicated can be found in Bolshev and Smirnov (1965) or  Markovich (2007). One can choose other nonparametric statistics like the Kolmogorov-Smirnov or Anderson-Darling ones instead of $\omega^2_n$.
Limit distributions of these statistics are invariant regarding  $F(x)$ (Bolshev and Smirnov 1965). We will focus  on $\omega^2_n$.
Regarding practical applications  the bandwidth $h$ was proposed in Markovich (1989)  as a solution of the equation
\begin{equation}\label{2a}\hat{\omega}^2_n(h)=0.05.\end{equation}
Here,
\[\hat{\omega}^2_n(h) =
\sum_{i=1}^n\left(\widehat{F}_h(X_{i,n})-\frac{i-0.5}{n}\right)^2+\frac{1}{12n}\]
was calculated by the order statistics $X_{1,n}\le...\le X_{n,n}$ corresponding to the sample $X^n$, and the value  $0.05$ 
 corresponding to the  mode of the pdf  of the statistic
$\omega^2_n$ and thus,  the maximum likelihood value of the $\omega^2_n$  
was found by tables of the statistic  $\omega^2_n$ 
(Bolshev and Smirnov 1965) as the discrepancy value $\delta$.
 A similar idea was explored in Markovich (2015) to estimate the extremal index.
  \\
  Following (Ferro and Segers 2003;  Markovich 2014, 2016b, 2017) 
  one can determine a cluster as the number of consecutive observations exceeding the threshold  $u$ between two consecutive non-exceedances.
 Ferro and Segers (2003) state that the times between exceedances of a threshold $u$ by the process $\{X_i\}$, is a random variable $T(u)$ equal in distribution to $T_1(u)$. It holds
  \begin{eqnarray*}T_1(u)&=&\min\{j\ge 1: M_{1,j}\le u, X_{j+1}>u|X_{1}>u\},\end{eqnarray*} where $M_{i,j}=\max\{X_{i+1},...,X_{j}\}$, $M_{1,1}=-\infty$.
 $T_1(u_n)$ normalized by the tail function $\{Y=\overline{F}(u_n)T_1(u_n)\}$, 
 is derived to be asymptotically exponentially distributed with a weight $\theta$ and with an atom at zero with a weight $1-\theta$
   (Ferro and Segers 2003).
\\
Taking the exceedance times $1\le S_1<...<S_{N_u}\le n$, the observed interexceedance times are $T_i=S_{i+1}-S_i$ for $i=1,...,N_u-1$, where $N_u=\sum_{i=1}^n\Ii\{X_i>u\}$ is the number of observations which exceed a predetermined high threshold $u$\footnote{Theoretically,  events $\{T_i=1\}$ 
are allowed. In practice, such cases that mean single inter-arrival times between consecutive exceedances are meaningless.
 } (Ferro and Segers 2003). Denote further $L\equiv L(u)=N_u-1$.
 In case of the statistic  $\omega^2_n$ the discrepancy equation may be calculated by the $k$, $1\le k\le L$, 
 largest order statistics of a sample  $\{Y_i=(N_u/n)T_i\}$  as follows
\begin{equation}\label{20}\hat{\omega}^2_L(u) =
\sum_{i=L-k+1}^L\left(\widehat{G}(Y_{i,L})-\frac{i-0.5}{L}\right)^2+\frac{1}{12L}=\delta.\end{equation}
Here, the distribution model of the normalized inter-exceedance times $\widehat{G}(Y_{i,L})$ is determined by $G(t)=1-\theta\exp(-\theta t)$ 
with a substitution of $\theta$ by some estimate $\widehat{\theta}(u)$ and of $t$ by the order statistic $Y_{i,L}$ (Markovich 2015). 
A value of 
$u$ can be found as a solution of 
(\ref{20}) 
with regard to any consistent nonparametric estimator of $\theta$.
The calculation (\ref{20}) by the entire sample may lead to the lack of a solution of the discrepancy equation regarding $u$ the same way as for the heavy-tailed pdf estimation in (Markovich 2007; Markovich 2016a)  or to too large 
$u$'s which may be not appropriate for the estimation of $\theta$.
\\
The selection of  $k$ and $\delta$ remains a problem.
To overcome this problem
we find a specific normalization 
of the discrepancy statistic $\hat{\omega}^2_L(u)$ in (\ref{20}) such that its limit distribution is the same as for the C-M-S statistic.  
Then its quantiles  may be used as $\delta$. Although the limit distribution of the C-M-S does not depend on $k$ as $k\to\infty$, for moderate samples the selection of $k$ is necessary. From Theorems \ref{T1} and  \ref{TR1}
it is natural to select $k$ such that $k\le \lfloor \theta L(u)\rfloor$. Specifically, for the discrepancy statistic the choice of $k$ has to be modified to satisfy Theorem \ref{TR2}. The discrepancy method can be easily applied not only to the threshold-based but generally to any extremal index estimator. 
In contrast to approaches based on  the minimum of the mean squared error with regard to the threshold or another tuning parameter, e.g. Robert et al. (2009b), the discrepancy method does not attract any knowledge about the asymptotic rates of the variance and the bias of estimates. Similarly to Hall (1990), the minimization of the bootstrap  mean squared error may require two additional parameters to select the size of a resample and to relate  thresholds based on the resamples and the entire sample.
\\
The paper is organized as follows. In Section \ref{Sec2} related work is recalled.  In Section \ref{Sec3} a normalization of
$\hat{\omega}^2_L(\theta)$ denoted as $\widetilde{\omega}^2_L(\theta)$ is found which has the 
limit distribution of $\omega^2_n$ (Theorem \ref{TR1}).
The convergence of $\widetilde{\omega}^2_L(\widehat{\theta})$ to the $\omega^2_n$-distribution is derived when the difference $\sqrt{m_n}(\widehat{\theta}-\theta)$, where $m_n$ is some sequence relating to $k$ and $L$, has a nondegenerate distribution (Theorem \ref{TR2}). In Theorem \ref{TR4} the consistency and the inconsistency conditions for the normalized statistic $\widetilde{\omega}^2_L(\hat{\theta})$ are given.
The  rate of convergence of the extremal index estimates with the threshold selected by the discrepancy method is derived in  Corollary \ref{Cor1}. The choice of the $k$ largest order statistics for  $\widetilde{\omega}^2_L$ by samples of moderate sizes
is discussed. Finally, an algorithm and a simulation study of the discrepancy method based on the normalized $\hat{\omega}^2_L(\widehat{\theta})$ 
statistic
is given in Section \ref{Sec4} and an illustration with real data is stated in Section \ref{Sec5}. 
Proofs can be found in Section \ref{Proofs}.

\section{Important mathematical results}\label{Sec2}
Our results are based on 
Lemmas 2.2.3, 3.4.1  by de Haan and Ferreira (2006) 
concerning the limit distributions of the order statistics and  
Theorem \ref{T1} by Ferro and Segers (2003). 
\begin{lemma}\label{Lem2} (de Haan and Ferreira 2006; Smirnov 1952) Let $U_{1,n}\le U_{2,n}\le ... \le U_{n,n}$ be the $n$th order statistics from a standard uniform distribution. Then, as $n\to\infty$, $k\to\infty$, $n-k\to\infty$,
\[\frac{U_{k,n}-b_{n}}{a_{n}}\]
is asymptotically standard normal with
\begin{eqnarray*}&& b_{n}=\frac{k-1}{n-1},\qquad a_{n}=\sqrt{b_n(1-b_n)\frac{1}{n-1}}.\end{eqnarray*}
\end{lemma}
\begin{lemma}\label{Lem1}
(de Haan and Ferreira 2006) Let $X, X_1, X_2,... ,X_n$ be i.i.d. r.v.s with common cdf $F$, and let $X_{1,n} \le X_{2,n} \le...\le X_{n,n}$ be the $n$th order statistics. The joint distribution of $\{X_{i,n}\}^n_{i=n-k+1}$ given $X_{n-k,n} = t$, for some $k\in\{1,..., n-1\}$,
equals the joint distribution of the set of order statistics $\{X^*_{i,k}\}^k_{i=1}$ of i.i.d. r.v.s $\{X^*_{i}\}^k_{i=1}$ with cdf
\begin{eqnarray*}F_t(x)&=&P\{X \le x|X > t\} = \frac{F(x) - F(t)}{1-F(t)}, ~~~ x>t.\end{eqnarray*}
\end{lemma}
\begin{definition}\label{Def-2} (Ferro and Segers 2003)  For real  $u$ and integers $1\le k\le l$, let
$\mathcal{F}_{k,l}(u)$
be the $\sigma$-field generated by the events $\{X_i>u\}$, $k\le
i\le l$.
Define the mixing coefficients $\alpha_{n,q}(u)$,
\begin{eqnarray*}\label{44}\alpha_{n,q}(u)&=&\max_{1\le k\le n-q}\sup|P(B|A)-P(B)|,
\end{eqnarray*}
where the supremum is taken over all $A\in
\mathcal{F}_{1,k}(u)$ with $P(A)> 0$ and $B\in
\mathcal{F}_{k+q,n}(u)$ and $k$, $q$ are positive integers.
\end{definition}
The next theorem  states that
\[\overline{F}(u_n)T_1(u_n)\to^d T_{\theta}=\left\{
\begin{array}{ll}
\eta, &   \mbox{with probability}~~ \theta,
\\
0, & \mbox{with probability}~~ 1-\theta,
\end{array}
\right.\]where $\eta$ is exponentially distributed with mean $\theta^{-1}$.  The zero asymptotic inter-exceedance times (the intracluster times) imply the times between the consecutive exceedances of the same cluster. The positive asymptotic inter-exceedance times are the inter-cluster times. $\to^d$ denotes convergence in distribution.
\begin{theorem}\label{T1} (Ferro and Segers 2003)
Let  $\{X_n\}_{n\ge 1}$ be a stationary process of r.v.s with tail function $\overline{F}(x)=1-F(x)$.
Let the positive integers  $\{r_n\}$ and the thresholds $\{u_n\}$, $n\ge 1$, be such that
$r_n\to\infty$,
$r_n\overline{F}(u_n)\to\tau$ and
 $P\{M_{r_n}\le u_n\}\to exp(-\theta\tau)$ hold as $n\to\infty$ for some
 $\tau\in(0,\infty)$ and $\theta\in (0,1]$. 
 If there are positive integers $q_n=o(r_n)$ such that
 $\alpha_{cr_n,q_n}(u_n)=o(1)$ for any $c>0$, then we get for $t>0$
\begin{equation}\label{10}P\{\overline{F}(u_n)T_1(u_n)>t\}\to\theta\exp(-\theta
t),\qquad n\to\infty.
\end{equation} 
\end{theorem}
 The intuition for
 declustering of a sample 
 is given in Ferro and Segers (2003). One can assume that the largest $C-1=\lfloor\theta L\rfloor$ inter-exceedance times are approximately independent inter-cluster times.
 The larger $u$ corresponds to the larger inter-exceedance times whose number $L\equiv L(u)$ may be small. It leads to a larger variance of the  estimates based on $\{T_i(u)\}$.
The intervals estimator follows from Theorem \ref{T1}. It is defined as (Beirlant et al. 2004, p.\ 391),
\begin{eqnarray}\label{15}
\hat{\theta}_n(u)&=&\Big\{
\begin{array}{ll}
\min(1,\hat{\theta}_n^1(u)), \mbox{ if } \max\{T_i : 1 \leq i \leq L\} \leq 2,\\
\min(1,\hat{\theta}_n^2(u)), \mbox{ if } \max\{T_i : 1 \leq i \leq L\} > 2,
\end{array}
\end{eqnarray}
where
\begin{eqnarray*}\label{16}&&\hat{\theta}_n^1(u)=\frac{2(\sum_{i=1}^{L}T_i)^2}{L\sum_{i=1}^{L}T_i^2},
~~~\hat{\theta}_n^2(u)=\frac{2(\sum_{i=1}^{L}(T_i-1))^2}{L\sum_{i=1}^{L}(T_i-1)(T_i-2)}.\end{eqnarray*}
 The $K$-gaps  estimator was proposed in S\"{u}veges and Davison (2010) as alternative to the intervals estimator, where the $K$-gaps 
 \[S(u_n)^{(K)}=(\max\left(T_1(u_n)-K,0\right)),~~K=0,1,2,...,\] have 
the same
limiting mixture law (\ref{10}). The $K$-gaps estimator is obtained by the maximum likelihood method using the model (\ref{10}) 
and assuming that the $K$-gaps observations 
are independent. It has the following form
\begin{eqnarray}\label{17}\widehat{\theta}^K&=&0.5\left((a+b)/c+1-\sqrt{((a+b)/c+1)^2-4b/c}\right),
\end{eqnarray}
with $a=L-N_C$, $b=2N_C$, $c=\sum_{i=1}^{L}\overline{F}(u_n)S(u_n)^{(K)}_i$. $N_C$ is the number of non-zero $K$-gaps. 
The $K$-gaps estimator is consistent and asymptotically normal as $L\to\infty$.
Since  the $K$-gaps may have a distribution different from 
(\ref{10}) for moderate samples,  $K$ is selected by a model misspecification test. The iterative weighted least squares estimator of S\"{u}veges (2007) explores the inter-exceedance times with $K=1$.  The automatic selection of an optimal pair $(u,K)$ is proposed in Fukutome et al. (2015) by a choice of   pairs for which 
values of the 
statistic of the information matrix test (the IMT) are less than $0.05$.
The test works satisfactorily when  the number of exceedances is not less than $80$.
\\
The intervals estimator is derived to be consistent  for $m$-dependent processes (Ferro and Segers 2003). Asymptotic normality  $\sqrt{m_n}(\hat{\theta}_n(u)-\theta)\to^d N(0,V)$ as $n\to\infty$ is derived for several extremal index estimators and  different values of variance $V$. 
Here, $m_n=np_n$ holds for the blocks and runs estimators, where $p_n=\P\{X_1>u_n\}$ holds in Weissman and Novak (1978); 
$m_n=L(u_n)$ 
is the number of inter-exceedance times $\{T_i(u_n)\}$ 
 for the intervals estimator in Robert (2009a), $\lfloor x \rfloor$ denotes the integer part of $x$; $m_n=n/p_n$ is taken for the multilevel blocks estimator in Sun and  Samorodnitsky (2018), where $p_n\overline{F}(u_n^s)\to\tau_s$ and $s\in\{1,...,m\}$ is the number of levels $\{u_n^s\}$; $m_n=\lfloor n/r_n\rfloor$ is used for the disjoint  and sliding blocks estimators, where $r_n=o(n)$ are positive integers related to a mixing condition 
 (Robert 2009a); $m_n = k_n$ is taken for the disjoint and
sliding blocks estimators by Berghaus and B\"{u}cher (2018) and Northrop (2015),
where $k_n$ is a number of blocks of length $b_n$ such that $k_n = o(b^2_n)$ holds as
$n \to \infty$.

\section{Main results: $\omega^2$-distribution of the normalized Cram\'{e}r-von Mises-Smirnov statistic}\label{Sec3}
\subsection{Normalized Cram\'{e}r-von Mises-Smirnov statistic for known $\theta$}
Let us rewrite the left-hand side (\ref{20}) in the following form
\begin{equation}\label{21}\hat{\omega}^2_L(u) =
\sum_{i=L-k+1}^L\left(1-\theta\exp(-Y_{i,L}\theta)-\frac{i-0.5}{L}\right)^2+\frac{1}{12L}\end{equation}
and derive its limit distribution.
Note that $L = L(u_{r_n})$ is 
a sequence of r.v.s converging in probability to infinity due to (Theorem 2.1,
Robert 2009a). In sequel,
the limit
distribution of the concerned statistics does not depend on $L$, so
we can neglect its randomness. The threshold sequence $u_{r_n}$ introduced in Robert (2009a) corresponds to the point process of time normalized exceedances  defined on $(0,\infty)$, in contrast to a traditional point process defined on $(0,1)$ (Beirlant et al. 2004).
\\
According to (Martynov 1978; Smirnov 1952) the limit distribution
of the 
C-M-S statistic (\ref{8}) or of $\omega^2_n=n\int_{0}^{1}\left(F_n(t)-t\right)^2dt$
coincides with the distribution of 
\begin{eqnarray*}\Omega &=&\int_0^1 B^2(t) dt,\end{eqnarray*} where $B(t)$ is a Brownian bridge on $[0,1],$
i.e. the Gaussian random process with zero mean and the covariance
function $R(s,t) = \min(s,t) - st,$ $s,t \in [0,1]$. Then
the statistic (\ref{21}) 
built by the $k$ largest order statistics tends 
to $0$ for $k = o(L)$ as $L\to\infty$, 
since
the
interval over which we integrate 
$B^2(t)$ 
tends to an empty set. 
Thus, (\ref{21}) 
has to be 
normalized. 
Let us consider 
the 
normalization of (\ref{21}) 
\begin{eqnarray}\label{omegasquared2}&&\widetilde{\omega}^2_L(\theta)
= \frac{1}{(1-t_k)^2}\cdot
\\
&\cdot &\sum\limits_{i=L-k+1}^{L}\left(1 - \theta\exp(-Y_{i, L}\theta)
-
 t_k - \frac{i - (L-k) -0.5}{k}(1-t_k)\right)^2 +
 \frac{1}{12k},\nonumber\end{eqnarray}
where $t_k = 1 - \theta\exp(-Y_{L-k, L}\theta)$. Let us explain in more detail why
we need 
such normalization.
It follows from (Robert 2009a, Theorem 2.1), 
that there are a probability, $\theta$, of asymptotic positive inter-exceedance times (the inter-cluster times) and a probability, $1-\theta,$ of zero asymptotic inter-exceedance times (the intra-cluster times). Moreover, the inter-cluster times are asymptotically independent exponential with mean 
$1/\theta$. Let us consider the following statistic
 \begin{eqnarray}\label{omegasquared}&&\omega^2_k(\theta) =\frac{1}{Z_{L-k, L}^2}\cdot\nonumber
\\
&\cdot &\sum\limits_{i=L-k+1}^{L}\left(Z_{L-k, L} - Z_{i, L}
  - \frac{i - (L-k) -0.5}{k}Z_{L-k, L}\right)^2 +
 \frac{1}{12k},\end{eqnarray}
 where 
$Z_{i, L}=\theta\exp(-T^*_{i, L}\theta)$, $T^*_{1, L}\leq \ldots \leq T^*_{L,L}$ are order statistics of a sample $\{T^*_i\},$ $\{T^*_i\}$ are independent copies of $T_\theta$. 

It follows
 from (\ref{omegasquared}), Theorem \ref{T1} 
and Lemma \ref{Lem1}, 
that the conditional distribution of the set
of order statistics $\{1 - Z_{i,L}\}_{i=L-k+1}^L$ 
given $1 - Z_{L-k, L} = t_k$ asymptotically agrees for $\limsup_{n\to\infty}k/L<\theta$ 
with the distribution of the set of
order statistics $\{U^\ast_{j,k}\},$ $j=i - (L-k),$ of an i.i.d. sample
$\{U_j^\ast\}$ from the uniform distribution on $[t_k,1]$.
The
asymptotical distribution of 
$\omega^2_k(\theta)$ is given in the next 
theorem.
\begin{theorem}\label{TR1} 
It holds
\begin{eqnarray*}&&\omega^2_k(\theta) \stackrel{d}{\to} \xi
\end{eqnarray*}
as $k\to\infty$, $L\to\infty$, $\limsup_{n\to\infty}k/L<\theta$, 
where $\xi$ has the distribution function $A_1,$ which is the limit
distribution function of the C-M-S 
statistic $\omega^2_n$.
\end{theorem}
\begin{remark} Based on the proof of Theorem \ref{TR1}, one can propose the goodness-of-fit test of von Mises' type
to check the hypothesis $H_0: F(x) = F_0(x)$ for sufficiently large
$x$ using the largest order statistics of a sample $\{X_i\}_{i=1}^n.$ The
test statistic is the following
\begin{eqnarray*}&&\omega^2_k =
\sum_{i=0}^{k-1} \left(\frac{F_0(X_{n-i,n}) - F_0(X_{n-k,n})}{1 - F_0(X_{n-k,n})} - \frac{k-i-0.5}{k}\right)^2+\frac{1}{12k}.\end{eqnarray*}
Theorem \ref{TR1} implies, that the limit
distribution of the statistic $\omega^2_k$ under the hypothesis $H_0$ does not depend on $k$
and $n$. It is equal to the limit distribution of the C-M-S 
statistic, if the hypothesis $H_0$ is true. 
The consistency of the proposed test follows from the equality in
distribution of the test statistic given $F_0(X_{n-k,n}) = t$ and the C-M-S
statistic.
The test is based on the largest order statistics of a
sample. 
It is reasonable
both if only the upper tail of the distribution is of interest
and/or if 
the largest order statistics of a sample are only available.
\end{remark}

\subsection{Normalized Cram\'{e}r-von Mises-Smirnov statistic for unknown $\theta$}
Here, we check whether one can substitute $\theta$ by its estimate $\widehat{\theta}$  
 in  (\ref{omegasquared2})
and find 
conditions imposed on $\widehat{\theta}$ under which
the limit distribution of $\widetilde{\omega}^2_L(\theta)$ will be the same as the limit distribution of $\omega^2_k(\theta).$
Recall again that the number of inter-exceedance times $L = L(u_{r_n})$ is a sequence of r.v.s tending to $+\infty$ (Robert 2009a). It follows from
Theorem \ref{TR1} that the limit distribution of 
$\omega^2_k(\theta)$ does not depend on $L(u_{r_n})$. 
We can assume that $L(u_{r_n})$ is a numerical sequence tending to
 infinity as $n\to\infty.$ We will write $L$ instead of
$L(u_{r_n})$ in the sequel.

In spirit of Theorem 3.2 (Robert 2009a), 
the limit distribution of the following statistic \[\sqrt{L}\left(\sum_{i=1}^{L-1} f(Y_i) - E f(Y_1)\right)\] for some continuous $f$
may not depend on a substitution of the set of r.v.s $\{T^*_i\}_{i=1}^{L-1}$
appearing in (\ref{omegasquared}) 
instead of $\{Y_i\}_{i=1}^{L-1}.$ Moreover, 
$T^*_i \stackrel{d}{=} T_\theta,$ $i\in\{1, \ldots, L-1\}$
and
there is a probability $\theta$ of
the nonzero elements of this set that are  independent exponentially distributed with parameter $\theta$. For convenience, we accept further $L$  instead of $L-1$. For these r.v.s,
Theorem 2.2.1 (de Haan and Ferreira 2006) implies that if $k/L\to 0$ and $k\to\infty$ as $L\to\infty,$ then
\begin{eqnarray*}&&
\sqrt{k}(T^*_{L-k, L}- \ln(L\theta/k)/\theta) = O_P(1).\label{asymp_exp}\end{eqnarray*} 
In light of
these remarks let us assume that there exists a sample of independent exponentially distributed 
r.v.s $\{E_i^{(L)}\}_{i=1}^{l}$ with mean $\theta^{-1}$ for all large enough $L$ such that
\begin{equation}Y_{L-k, L} - E_{l-k, l}^{(L)} = o_P\left(\frac{1}{\sqrt{k}}\right)\label{newcondition}\end{equation}
if $k/L\to 0$ and $k\to\infty$ as $L\to\infty$, where we denote $l = \lfloor\theta L\rfloor$ and assume $k<l$. Here and further, we denote for brevity a sequence of positive integers $\{k_n\}$ as $k$.
Theorem \ref{TR1} remains valid  
when $T^*_{i,L}$, $i\in\{L-k,...,L\}$ 
in $\omega_k^2(\theta)$ are substituted by $E^{(L)}_{i,l}$, $i\in\{l-k,...,l\}$.
\begin{theorem}\label{TR2} Let 
the conditions of Theorem \ref{T1} and the condition (\ref{newcondition}) be fulfilled and the estimator of the
extremal index $\widehat{\theta} = \widehat{\theta}_{n}$ 
 be such
that
\begin{eqnarray}\label{2}&&\sqrt{m_n}(\widehat{\theta}_{n} - \theta)
\stackrel{d}{\to} \zeta,~~ n\to\infty, 
\label{normality}\end{eqnarray}
where the r.v. $\zeta$ has a nondegenerate distribution function $H$. 
Let us assume that the
sequence $m_n$ is such that
\begin{equation}\frac{k}{m_n} = o(1)\ \mbox{ and }\ \frac{(\ln L)^2}{m_n} = o(1)\label{asymp}\end{equation} as $n\to\infty$. 
Then
\begin{eqnarray*}
\widetilde{\omega}^2_L(\widehat{\theta}_{n})& \stackrel{d}{\to}& \xi \sim A_1\end{eqnarray*}
holds, where $A_1$ is the limit distribution function of the C-M-S 
statistic.
\end{theorem}
\begin{remark}Normal distributions 
give examples of $H$ 
regarding the intervals, blocks and sliding blocks estimators of the extremal index (Northrop 2015; Robert 2009a; Robert et al. 2009;   Sun and  Samorodnitsky 2018).
\end{remark}
\begin{remark} The replacement of $o(1)$ by $O(1)$ in  (\ref{asymp}) violates Theorem \ref{TR2}. The assumption $k = O(m_n)$ may lead to a limit distribution of $\widetilde{\omega}^2_L(\widehat{\theta}_{n})$ 
different from $A_1$ that is out of scope of the paper.
\end{remark}
\begin{remark} The limit process of the point process of exceedance times is a compound
Poisson process (Hsing et al. 1988). The condition (\ref{newcondition}) shows in fact the rate of this convergence required for the limit distribution of the normalized C-M-S statistic that is built by the largest $k$ order statistic to preserve the limit distribution of $\omega^2$.
\end{remark}
\begin{theorem}\label{TR4} Let the conditions of Theorem \ref{T1} and (\ref{newcondition}) be fulfilled. Assume that the sequence of estimates $\{\hat{\theta}_{n}\}$ 
is such that for some $\alpha \in [0,1/2]$ 
\begin{eqnarray*}&&
k_{n_s}^\alpha|\hat{\theta}_{n_s} - \theta|
\stackrel{P}{\rightarrow} +\infty,\qquad\mbox{if}\qquad 0<\alpha\le 1/2,
\\
&& 
|\hat{\theta}_{n_s} - \theta|>\varepsilon \qquad\mbox{for some}\qquad \varepsilon>0, \qquad\mbox{if}\qquad \alpha=0
\end{eqnarray*}
hold as $n\to \infty$
for some 
subsequence
$\{k_{n_s}\},$ $s\geq 1,$ of the sequence $\{k_n\}$.
Then for corresponding subsequence $\{L_s\}$ of the sequence $\{L\}$
\begin{eqnarray*}&&\widetilde{\omega}^2_{L_s}(\hat{\theta}_{n_s})/k_{n_s}^{1-2\alpha} \stackrel{P}{\rightarrow} +\infty\end{eqnarray*}
holds as $n\to\infty.$
\end{theorem}
\begin{remark} Theorem \ref{TR4} implies that the non-consistency of the estimator $\widehat{\theta}_{n}$ 
or the consistency with a sufficiently slow rate leads to the non-consistency of $\widetilde{\omega}^2_{L_s}(\widehat{\theta}_{n})$ 
in a sense that its limit distribution does not exist or the latter statistic tends to $+\infty$. In case that $\alpha\neq 0$ holds, the estimator $\widehat{\theta}_{n}$ 
may be consistent but with the rate of convergence slower than 
$k_n^{-\alpha}$.
  Hence, 
 $\widetilde{\omega}^2_L(\widehat{\theta}_{n})$
  may be considered as a quality functional of 
  $\widehat{\theta}_{n}$.
\end{remark}
The consistency of the corresponding extremal index estimates follows from Theorem \ref{TR4}.
The next corollary states, if the solutions of the discrepancy equation exist for each $n$, then the consistency is fulfilled.
\begin{corollary}\label{Cor1} Let $\hat{\theta}_n(u_n)$ be an estimator of $\theta$ and 
$\{\tilde{u}_{k,L}\}$ be some sequence of solutions of the discrepancy equation. Then $\hat{\theta}_n(\tilde{u}_{k,L}) \stackrel{P}{\rightarrow} \theta$
and for arbitrary 
$\varepsilon>0$
\begin{eqnarray*}\label{13}
&& k^{1/2-\varepsilon}|\hat{\theta}_n(\tilde{u}_{k,L}) - \theta| 
\stackrel{P}{\rightarrow} 0\end{eqnarray*} hold as $k\to\infty,$ $L/k\to\infty,$ $L=o(n)$, $n\to\infty.$
\end{corollary}
The proof of the corollary is based on a negation of the assertion of Theorem \ref{TR4}.
\subsection{The choice of $k$ 
}\label{Sec3.3}
According to Theorem \ref{TR2} 
the asymptotic distribution of $\widetilde{\omega}^2_L(\widehat{\theta}_n)$ 
does not depend on $k$. 
The $k$-selection gives another
viewpoint that using only the largest inter-exceedance times screens out the smallest inter-exceedance
times. It is helpful for the reasons discussed in Ferro and Segers (2003) and is
the motivation for the introduction of the tuning parameter $K$ in the $K-$gaps estimator of $\theta$ proposed in S\"{u}veges and Davison (2010).
\\
In practice, for each predetermined $\delta$, $u$ and $L(u)$ one may decrease the $k$-value such that
$k\le \min\{\widehat{\theta}_0 L(u), L(u)^{\beta}\}$, $0<\beta<1$,
($\widehat{\theta}_0$ is some pilot estimate of $\theta$) until the discrepancy equations have solutions and select the largest one among such $k$'s. 
This choice satisfies Theorem \ref{TR2} but it is not unique. For instance, one can select $k=\lfloor(\ln L)^2\rfloor)$. The following simulation shows an ideal case when the accuracy is the best. Namely, the $K-$gaps estimator with $K=0$ and $N_C=k$ coupled with the discrepancy method demonstrates the best choice  when $k=\lfloor\theta L\rfloor$ is chosen. Since $\theta$ is in reality unknown, one has to take $k=\lfloor\widehat{\theta}_0 L\rfloor$. This choice requires an accurate  consistent pilot estimate $\widehat{\theta}_0$. 
\section{Simulation study}\label{Sec4}
In our simulation study we focus on the threshold-based intervals and $K-$gaps estimators.
We propose also a modification of the $K-$gaps estimator with $K=0$ and $N_C=k$ in (\ref{17}) notated as $\widehat{\theta}^{K_0}$. The latter coupled with the discrepancy method demonstrates the best accuracy if an estimate $\widehat{\theta}_0$ is close to $\theta$. The natural drawback of the intervals estimator is that it needs a large sample size $n$ to obtain a moderate size $L(u)$ for a large $u$. The same concerns the $K$-gaps  estimator.
\begin{algorithm}
\begin{enumerate}
\item\label{I1} Using $X^n=\{X_i\}_{i=1}^n$ and taking thresholds $u$ corresponding to quantile levels $q\in\{0.90,0.905,...,0.995\}$, 
generate samples of the inter-exceedance times $\{T_i(u)\}$ and the normalized r.v.s
\begin{equation}\{Y_{i}\}=\{\overline{F}(u)T_i(u)\}=\{(N_u/n)T_i(u)\}, ~ i\in\{1, 2,...,L\}, ~ L=L(u),\end{equation}
where $N_u$  is the number of exceedances over threshold $u$.
\item\label{Item2} For each $u$ select $k=\lfloor \widehat{\theta}_0 L\rfloor$, $k=\min\{\lfloor \widehat{\theta}_0L\rfloor, \sqrt{L}\}$ (in case $\widehat{\theta}_0=1$, accept $k=L-1$) or $k=\lfloor(\ln L)^2\rfloor$, where the intervals estimator (\ref{15}) may be selected as a pilot estimator $\widehat{\theta}_0=\widehat{\theta}_0(u)$ with the same $u$ as in Item \ref{I1}. 
    \item Use a sorted sample $Y_{L-k+1,L}\le ...\le Y_{L,L}$ and find among considered quantiles all  solutions $u_1,...,u_l$ (here, $l$ is a random number) 
     of the following discrepancy equation
\begin{eqnarray}\label{4}
&&\!\!\!\!\!\!\widetilde{\omega}^2_L(\widehat{\theta})=\frac{1}{(1-\widehat{t}_k)^2}\cdot
\\
&&\!\!\!\!\!\!\!\!\!\!\sum\limits_{i=L-k+1}^{L}\left(1 - \widehat{\theta}\exp(-Y_{i, L}\widehat{\theta})
-
 \widehat{t}_k - \frac{i - (L-k) -0.5}{k}(1-\widehat{t}_k)\right)^2
 +  \frac{1}{12k}=\delta_1,\nonumber\end{eqnarray}
where $\widehat{t}_k = 1 - \widehat{\theta}\exp(-Y_{L-k, L}\widehat{\theta})$, $\widehat{\theta}=\widehat{\theta}(u)$ is calculated by (\ref{15}), 
and $\delta_1=0.05$ 
 is the mode
 of the C-M-S statistic. 
If $L<40$ we should replace $\widetilde{\omega}^2_L(\widehat{\theta})$ by
\begin{eqnarray*}
&&(\widetilde{\omega}^2_L(\widehat{\theta}))'
=\left(\widetilde{\omega}^2_L(\widehat{\theta})-\frac{0.4}{L}+\frac{0.6}{L^2}\right)\left(1+\frac{1}{L}\right)
\end{eqnarray*}
and use 
quantiles of the C-M-S statistic as the discrepancy $\delta$ (Kobzar 2006).
    \item For each $u_j$, $j\in\{1,...,l\}$ calculate $\hat{\theta}(u_j)$ and find
    \begin{equation}\label{19}
   \widehat{\theta}_1 =\frac{1}{l}\sum_{i=1}^l\widehat{\theta}(u_i),~ \widehat{\theta}_2 =\widehat{\theta}(u_{min}),~ \widehat{\theta}_3 =\widehat{\theta}(u_{max})
    \end{equation}
    as resulting estimates, where
    $u_{min}=\min\{u_1,...,u_l\}$, 
    $u_{max}=\max\{u_1,...,u_l\}$.
 \end{enumerate}
 \end{algorithm}
 \begin{figure}[tbp]
\centering
 \begin{minipage}{0.98\textwidth}
 \includegraphics[width=0.5\textwidth]{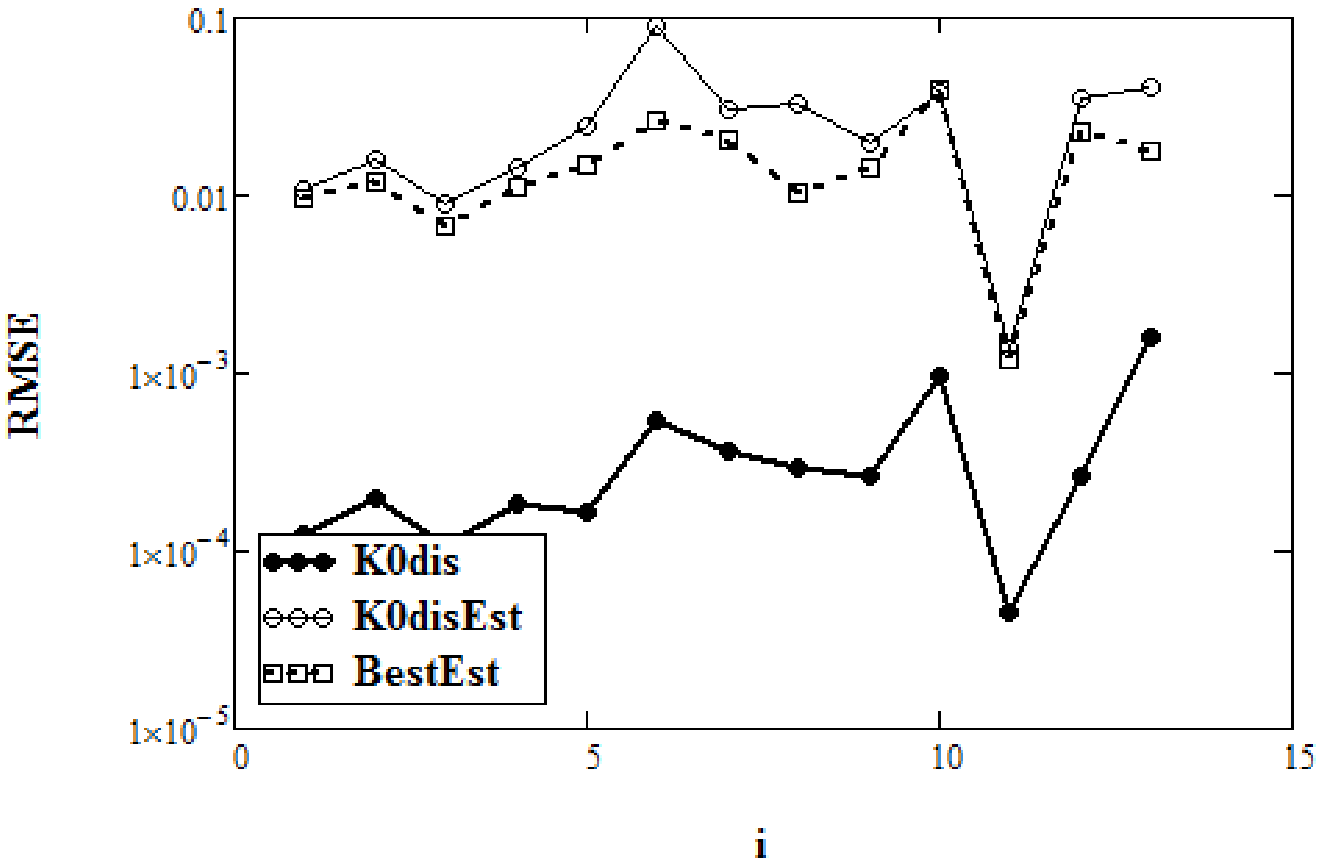}
\includegraphics[width=0.5\textwidth]{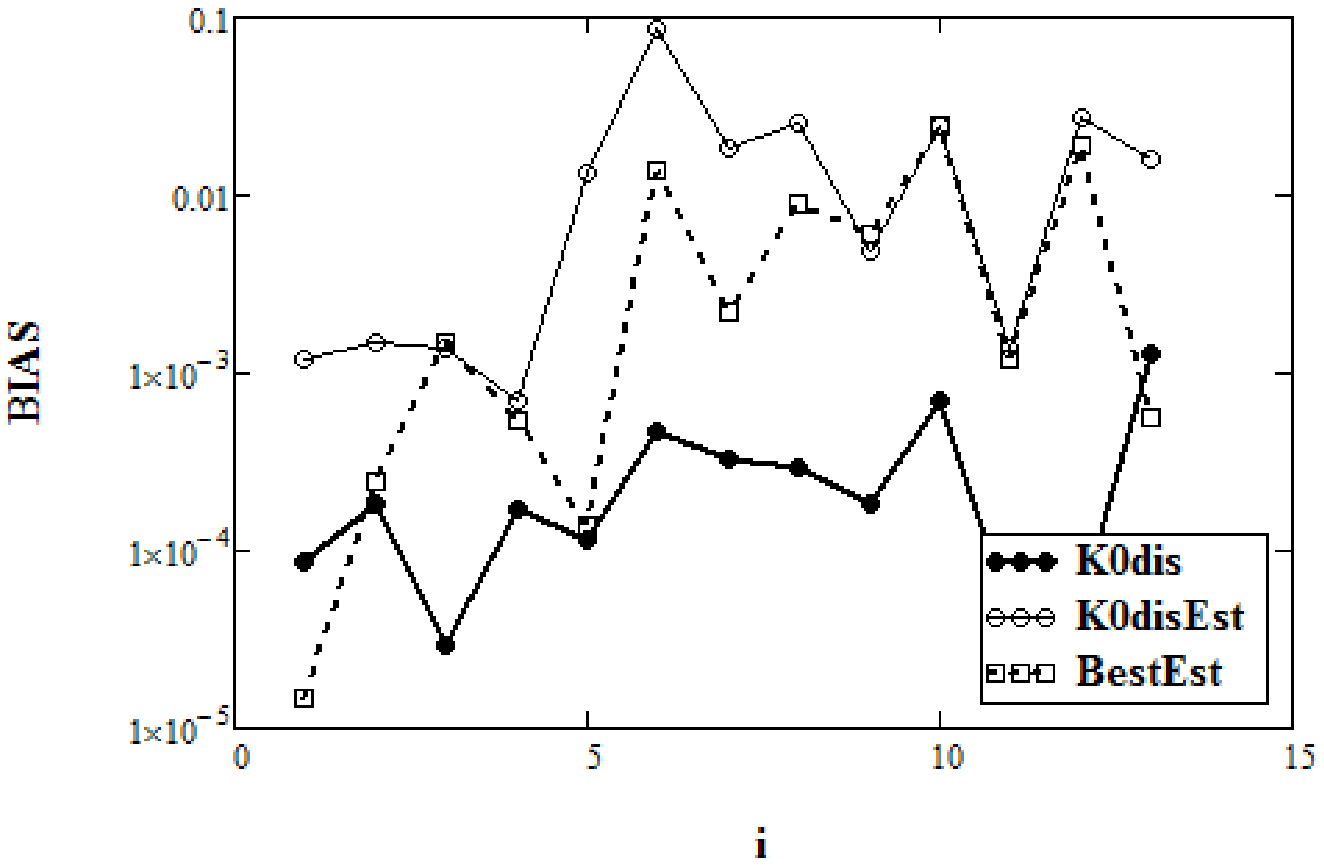}
\\
\includegraphics[width=0.5\textwidth]{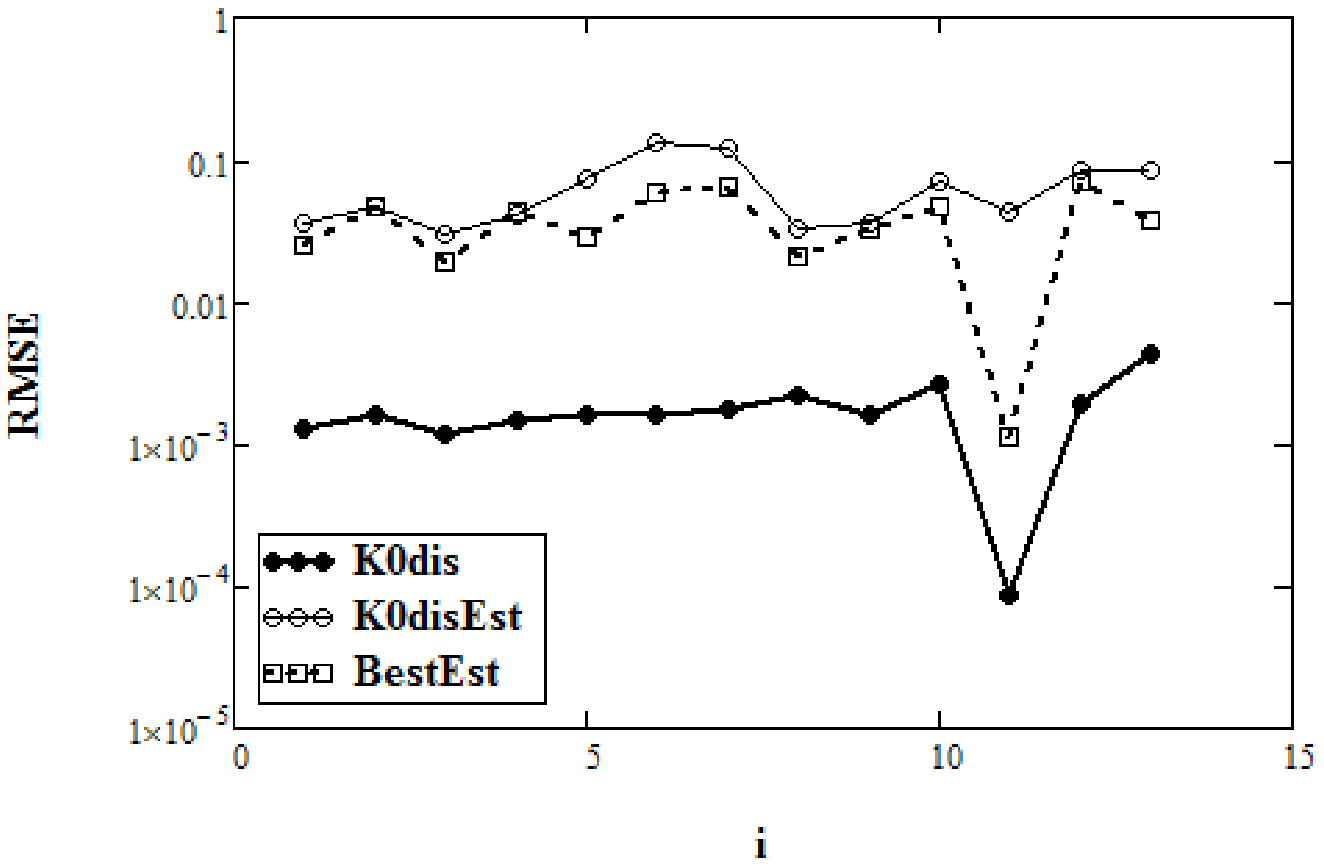}
\includegraphics[width=0.5\textwidth]{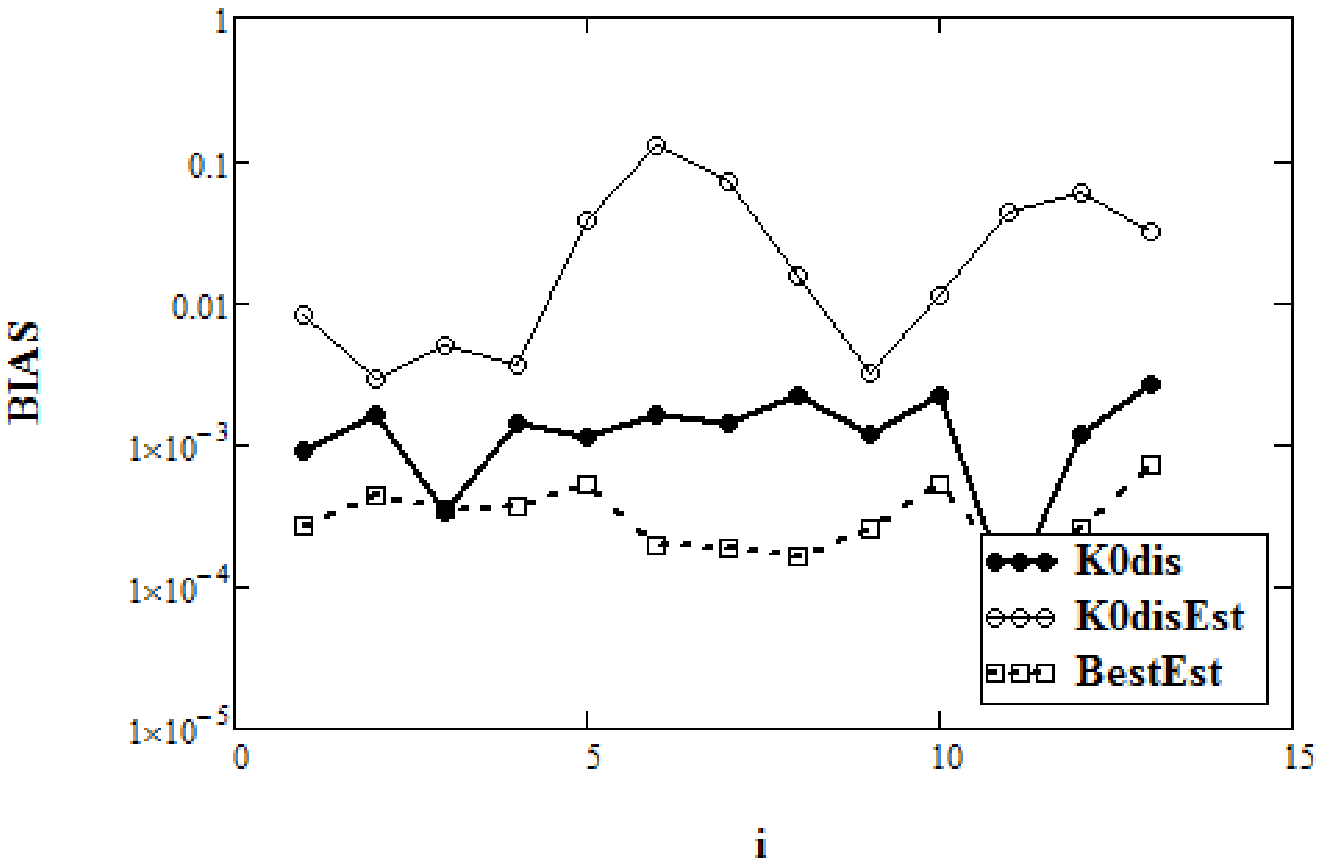}
 \end{minipage}
 \hfill
 \caption{The best RMSE and Bias for the K-gaps estimator $\widehat{\theta}^{K_0}$ 
  with the threshold $u$ selected by the discrepancy equation (\ref{4}) 
 and the corresponding inequality (\ref{5}) and with $k=\lfloor s L\rfloor$, where 'K0dis' and 'K0disEst' correspond to $s=\theta$ and $s=\widehat{\theta}_0$ in Tables \ref{Table1-KgapsK0Theta} and \ref{Table2-KgapsK0Theta}, respectively,  $\widehat{\theta}_0$ is a pilot intervals estimate; the best  RMSE and Bias among all estimates
 in Tables \ref{Table1-N}-\ref{Table2-K} notated as 'BestEst'
 against the number of processes related to the column labels in Tables \ref{Table1-KgapsK0Theta}-\ref{Table2-K}, and enumerated from left to right as in the tables for sample size $n=10^5$ (the upper row) and $n=5000$ (the lower row).}
 \label{fig:0}
\end{figure}
 \begin{figure}[tbp]
\centering
 \begin{minipage}{0.98\textwidth}
 \includegraphics[width=0.5\textwidth]{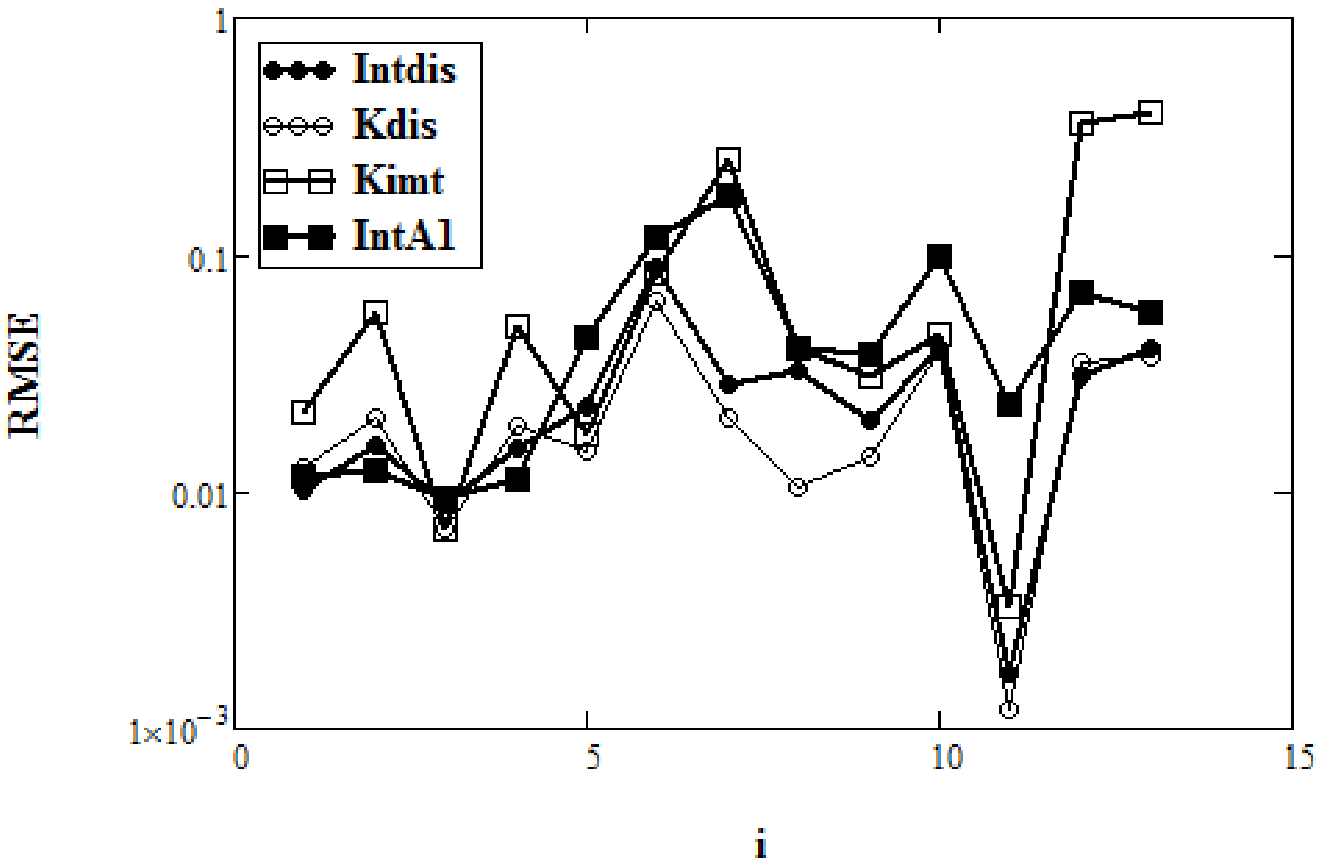}
\includegraphics[width=0.5\textwidth]{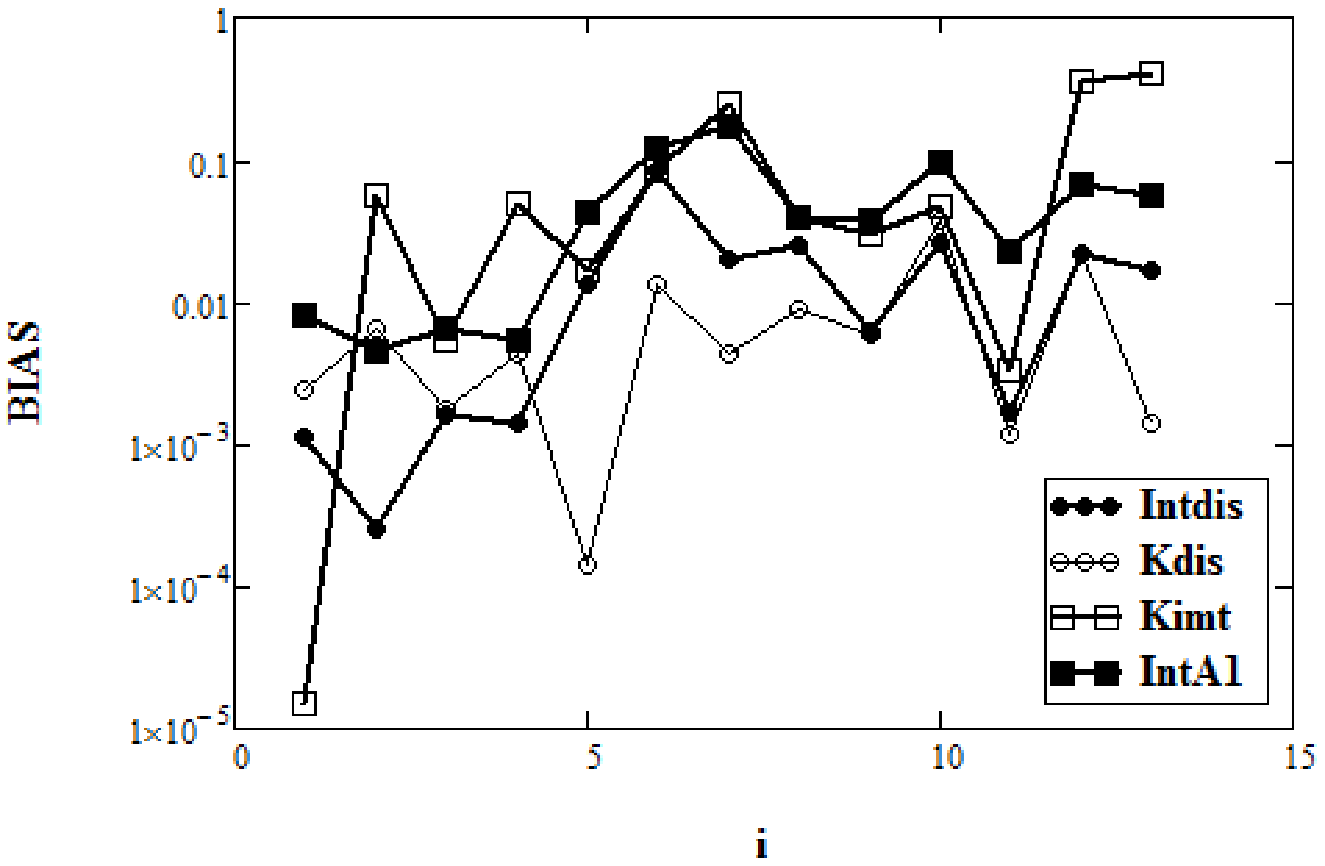}
 \end{minipage}
 \hfill
 \caption{The best RMSE and Bias for the intervals estimator ('Intdis') and K-gaps ('Kdis') estimators with threshold $u$ selected by the discrepancy equation (\ref{4}) 
 and the corresponding inequality (\ref{5}), and for the K-gaps estimator with $u$ selected by the test IMT ('Kimt') and the intervals estimator with the "plateau-finding" algorithm A1 to select $u$ ('IntA1')
 against the number of processes related to the column labels in Tables \ref{Table1-N} and \ref{Table2-N}, and enumerated from left to right as in the tables for sample size $n=10^5$.}
 \label{fig:1}
\end{figure}
\begin{figure}[tbp]
\centering
 \begin{minipage}{0.98\textwidth}
 \includegraphics[width=0.5\textwidth]{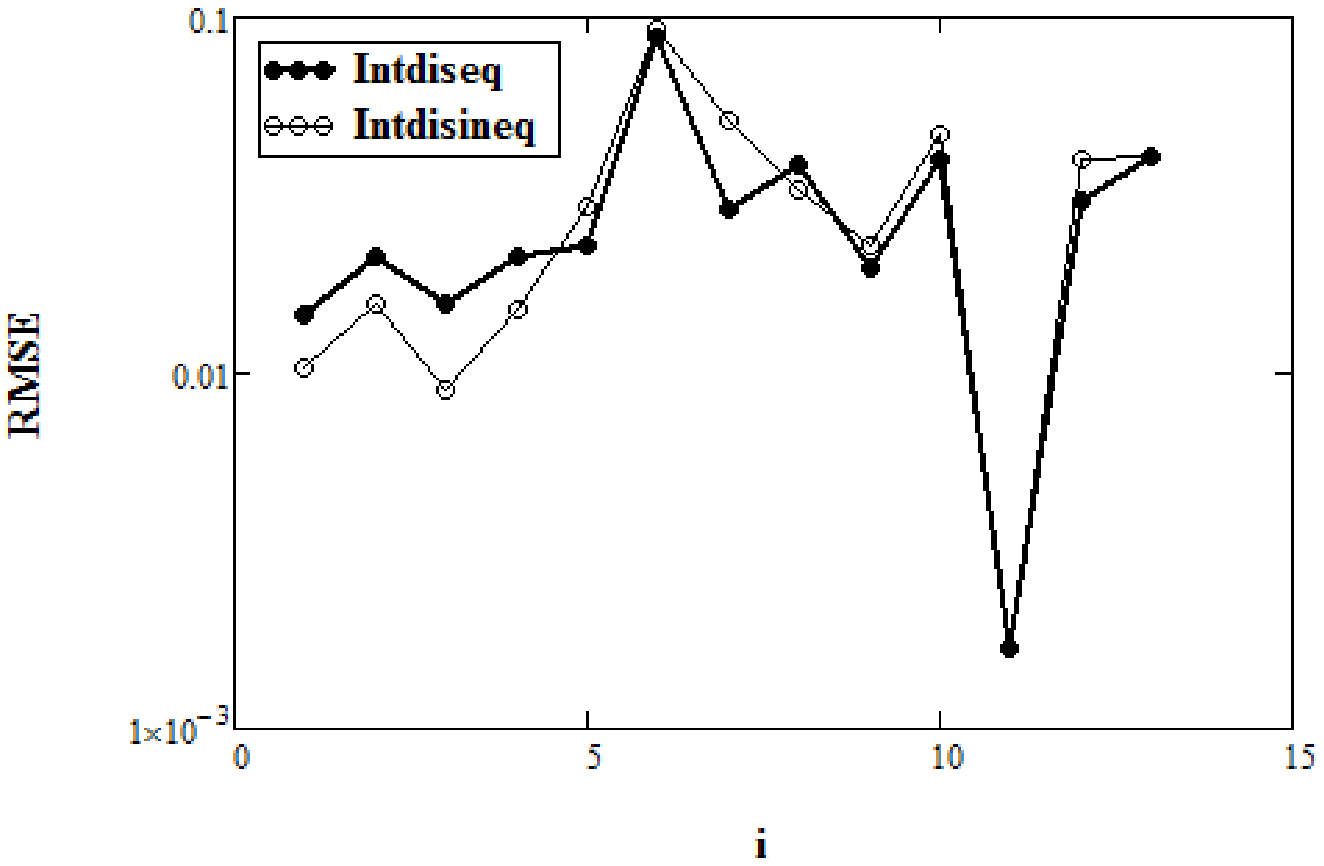}
\includegraphics[width=0.5\textwidth]{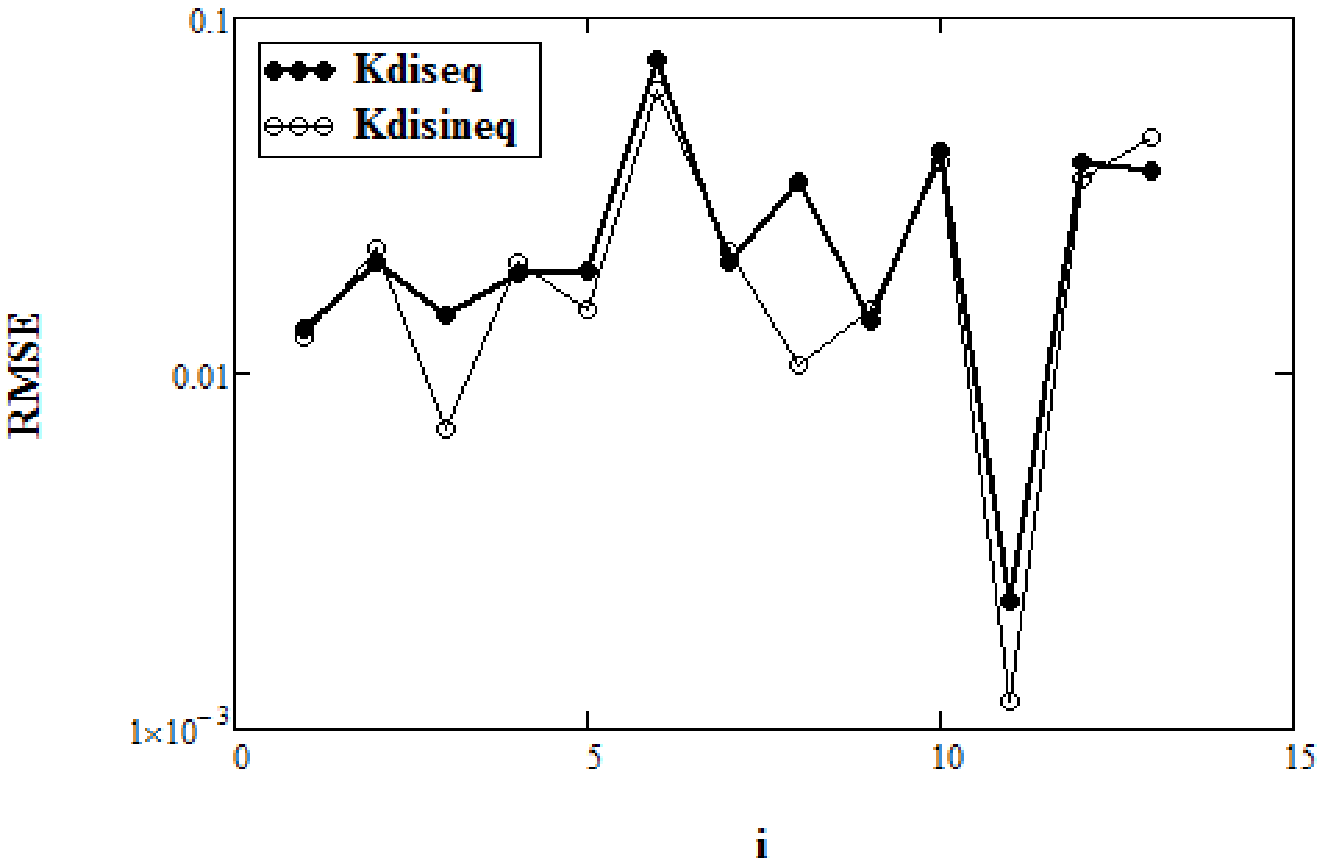}
 \end{minipage}
 \hfill
 \begin{minipage}{0.98\textwidth}
 \includegraphics[width=0.5\textwidth]{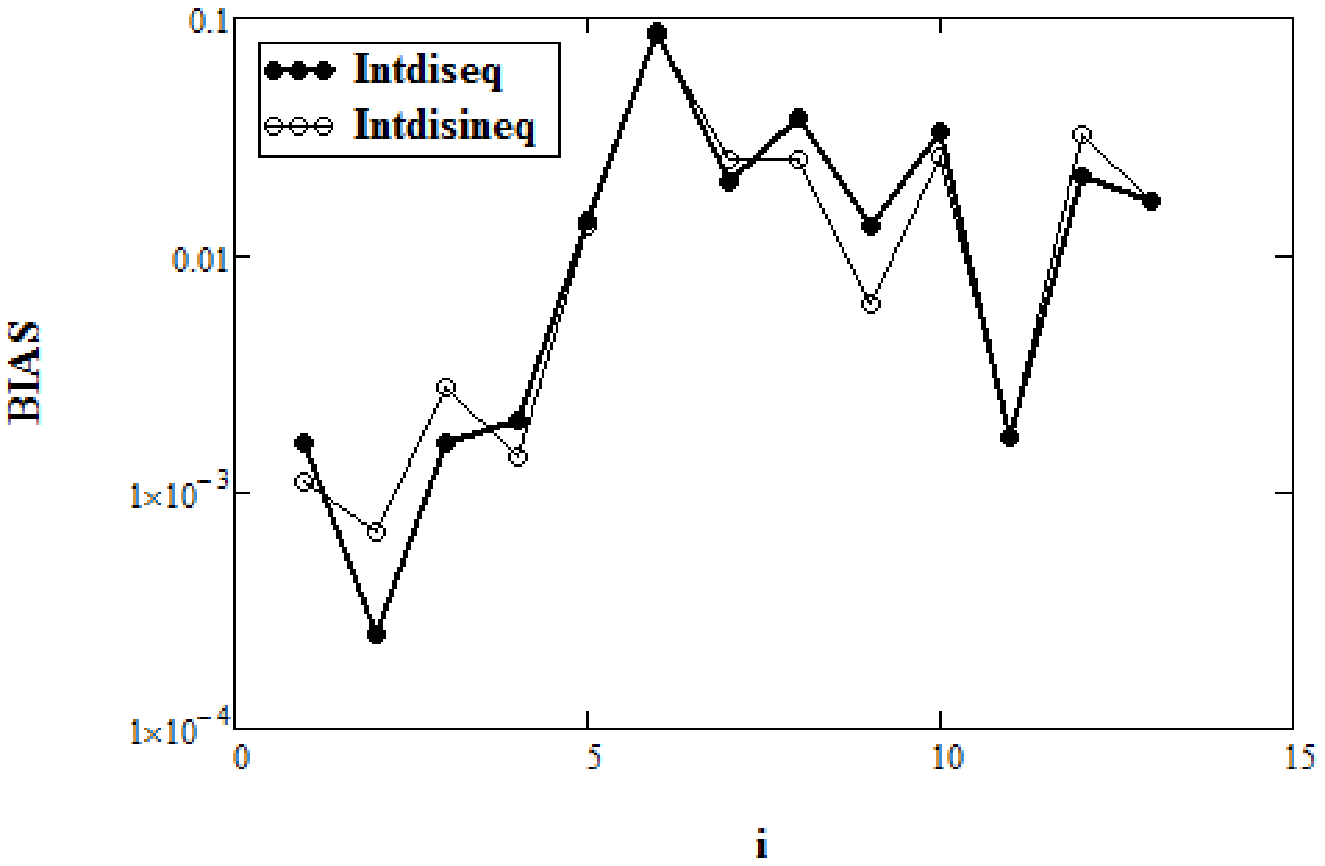}
\includegraphics[width=0.5\textwidth]{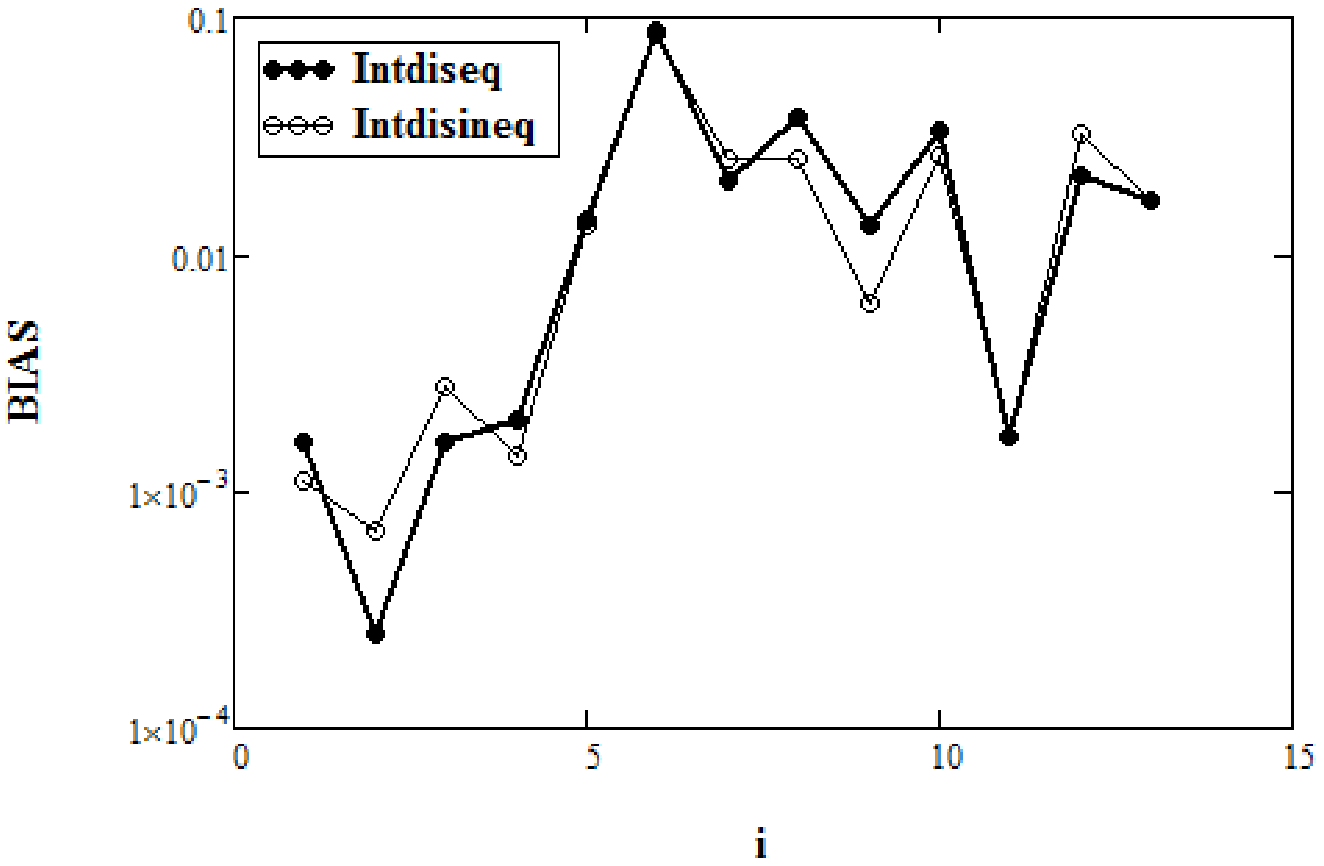}
 \end{minipage}
 \hfill
 \caption{The best RMSE and Bias for the intervals estimator (left column)
and K-gaps estimator (right column) obtained by the Algorithm with (\ref{4}) 
notated as 'Intdiseq' and  'Kdiseq', and with the inequality (\ref{5}) 
notated as 'Intdisineq' and 'Kdisineq' against the number of processes related to the column labels in Tables \ref{Table1-N} and \ref{Table2-N} for sample size $n=10^5$.}
 \label{fig:2}
\end{figure}
\begin{figure}[tbp]
\centering
 \begin{minipage}{0.98\textwidth}
 \includegraphics[width=0.5\textwidth]{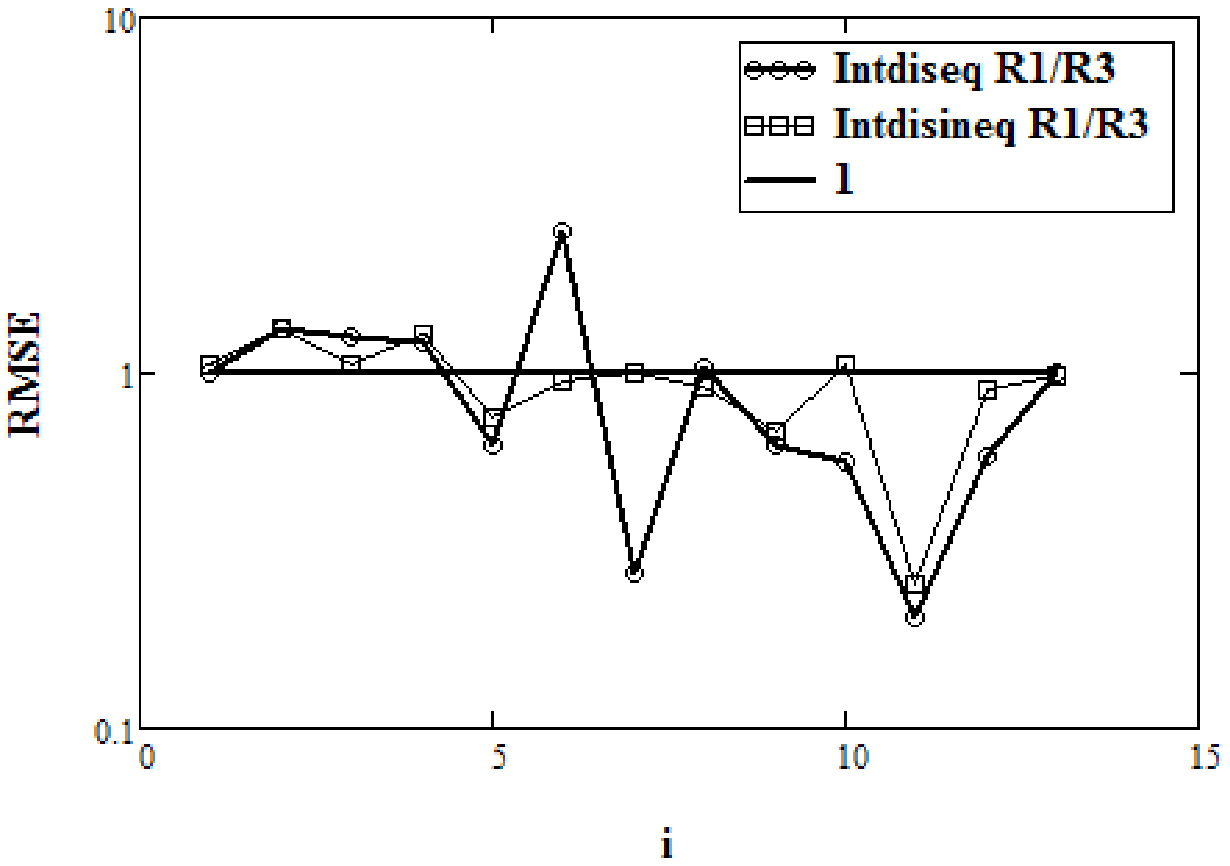}
\includegraphics[width=0.5\textwidth]{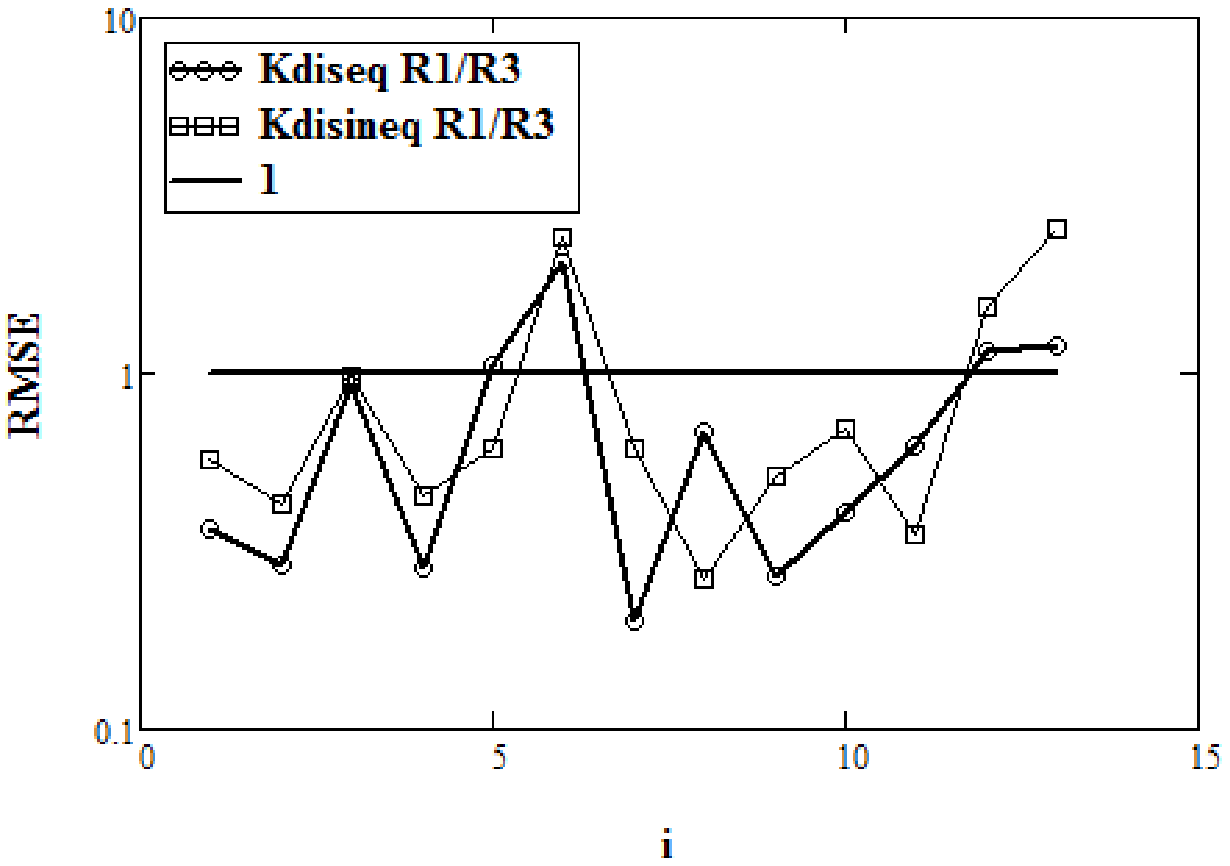}
 \end{minipage}
 \hfill
 \begin{minipage}{0.98\textwidth}
 \includegraphics[width=0.5\textwidth]{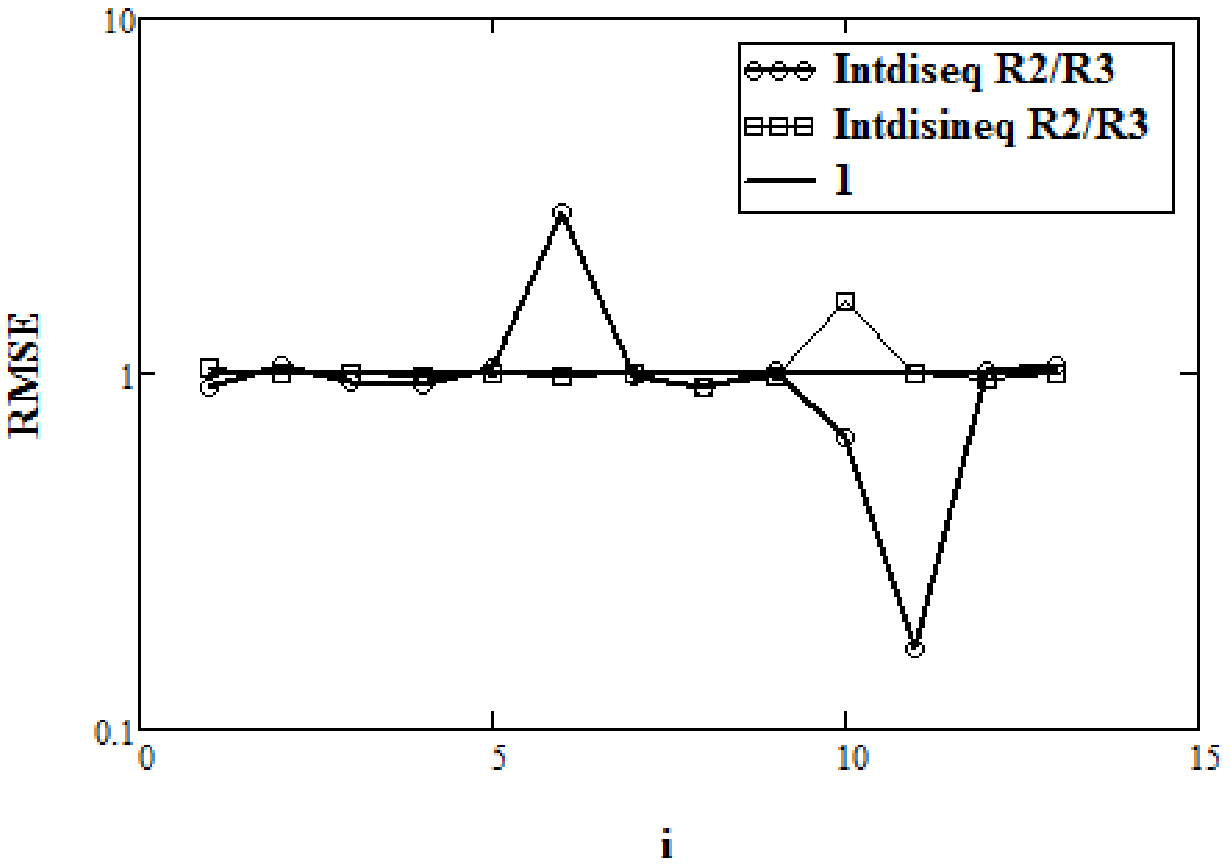}
\includegraphics[width=0.5\textwidth]{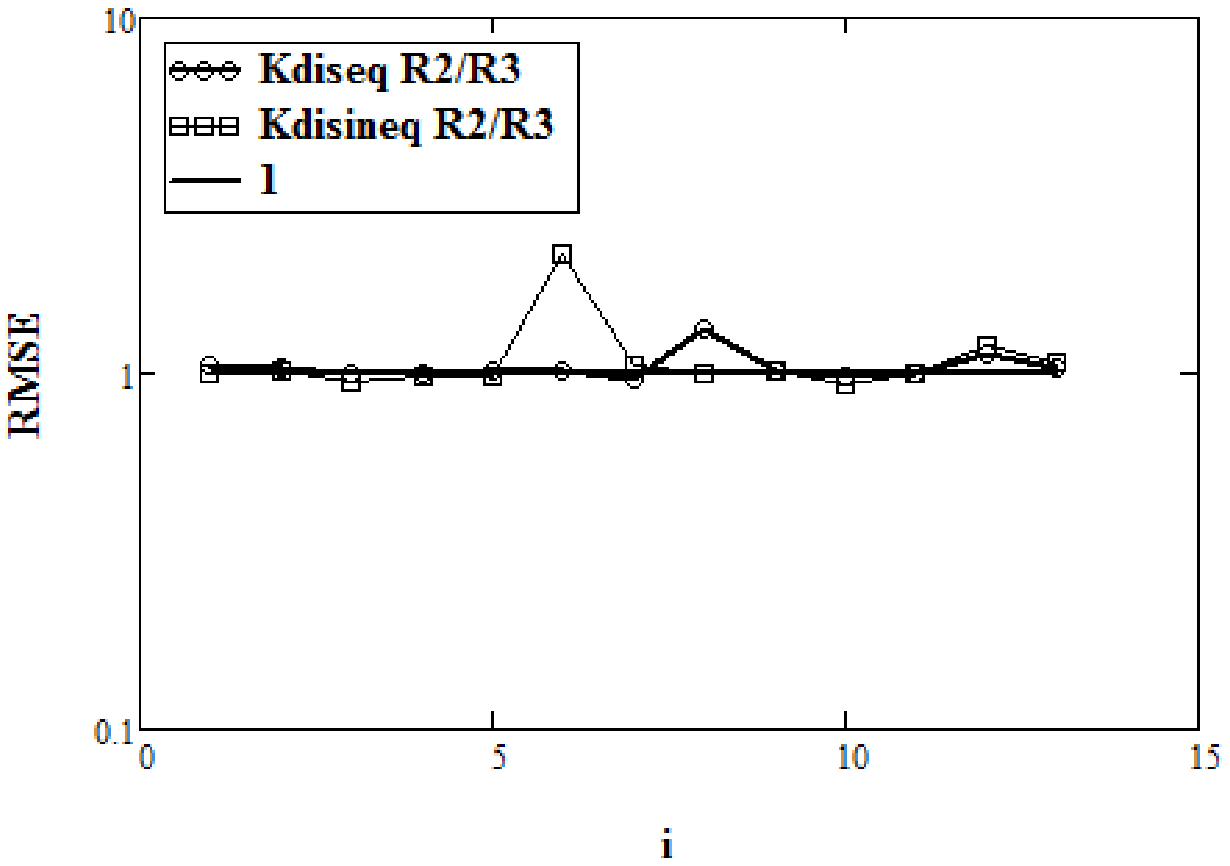}
 \end{minipage}
 \hfill
 \caption{Ratios $R1/R3$ and $R2/R3$ 
 of the best RMSE 
 for the intervals estimator (left column)
and K-gaps estimator (right column) obtained by the Algorithm with the equation (\ref{4})
notated as 'Intdiseq' and  'Kdiseq', and with the inequality (\ref{5}) notated as 'Intdisineq' and 'Kdisineq' against the number of processes related to the column labels in Tables \ref{Table1-N}, \ref{Table1-I} and \ref{Table1-K}: The $R1/R3$ corresponds to the best results in Table \ref{Table1-N}
divided to those in Table \ref{Table1-K}, and the $R2/R3$ - to those in Tables \ref{Table1-I} and \ref{Table1-K}, respectively, for sample size $n=10^5$.}
 \label{fig:3}
\end{figure}
\begin{figure}[tbp]
\centering
 \begin{minipage}{0.98\textwidth}
 \includegraphics[width=0.5\textwidth]{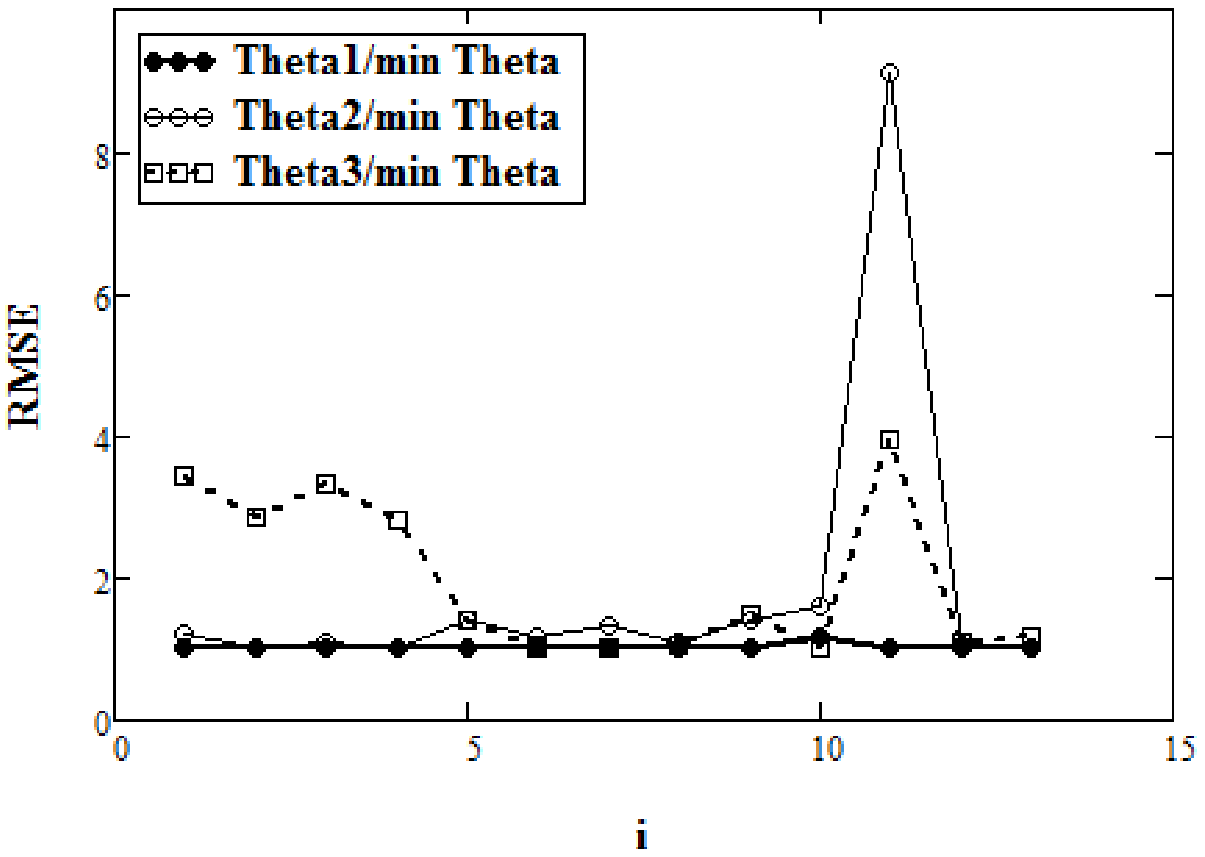}
\includegraphics[width=0.5\textwidth]{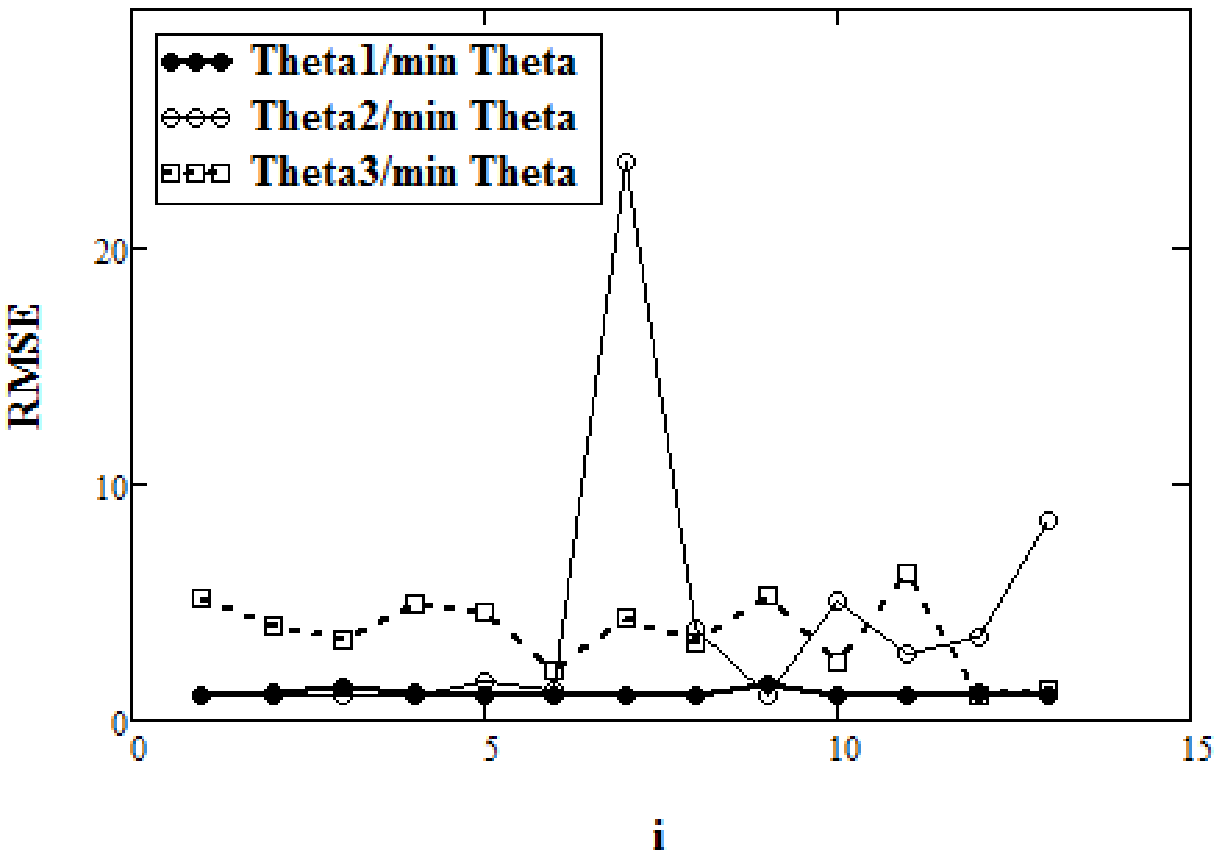}
 \end{minipage}
 \hfill
 \begin{minipage}{0.98\textwidth}
 \includegraphics[width=0.5\textwidth]{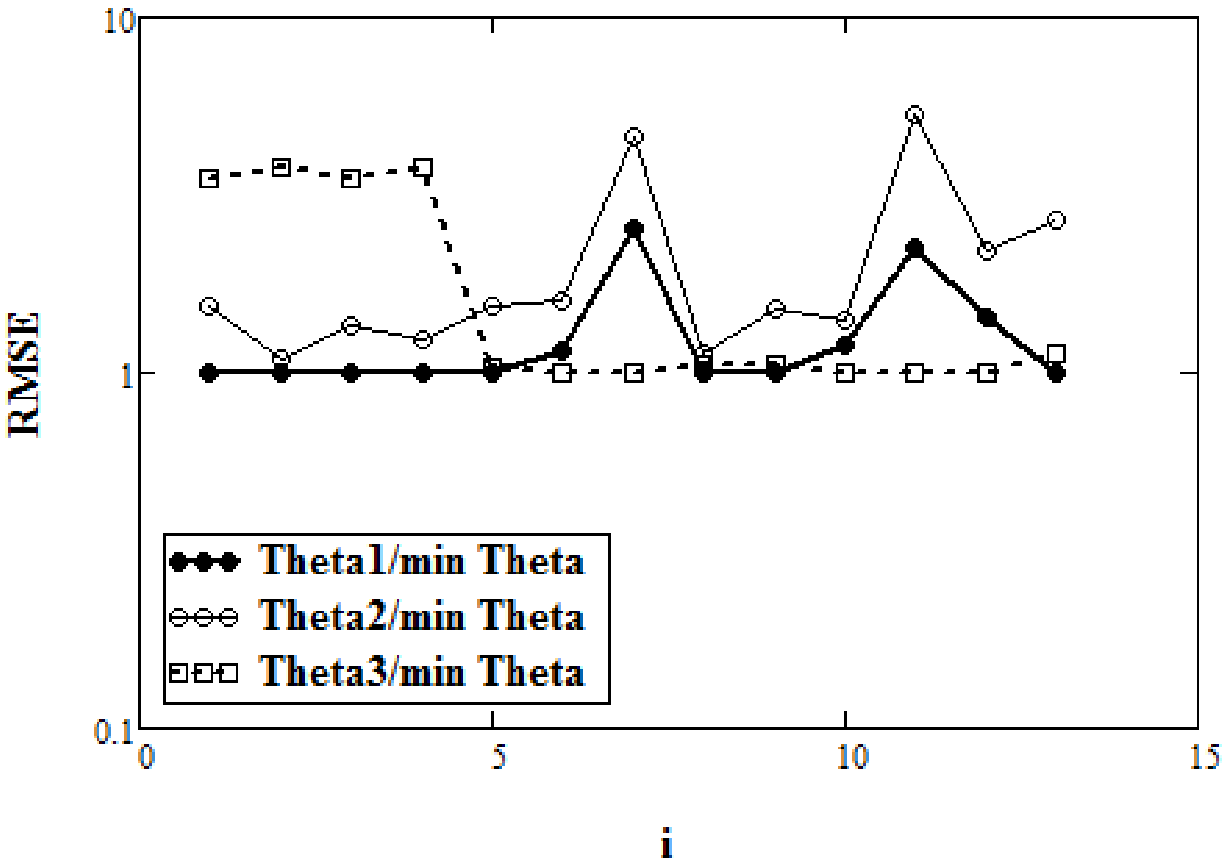}
\includegraphics[width=0.5\textwidth]{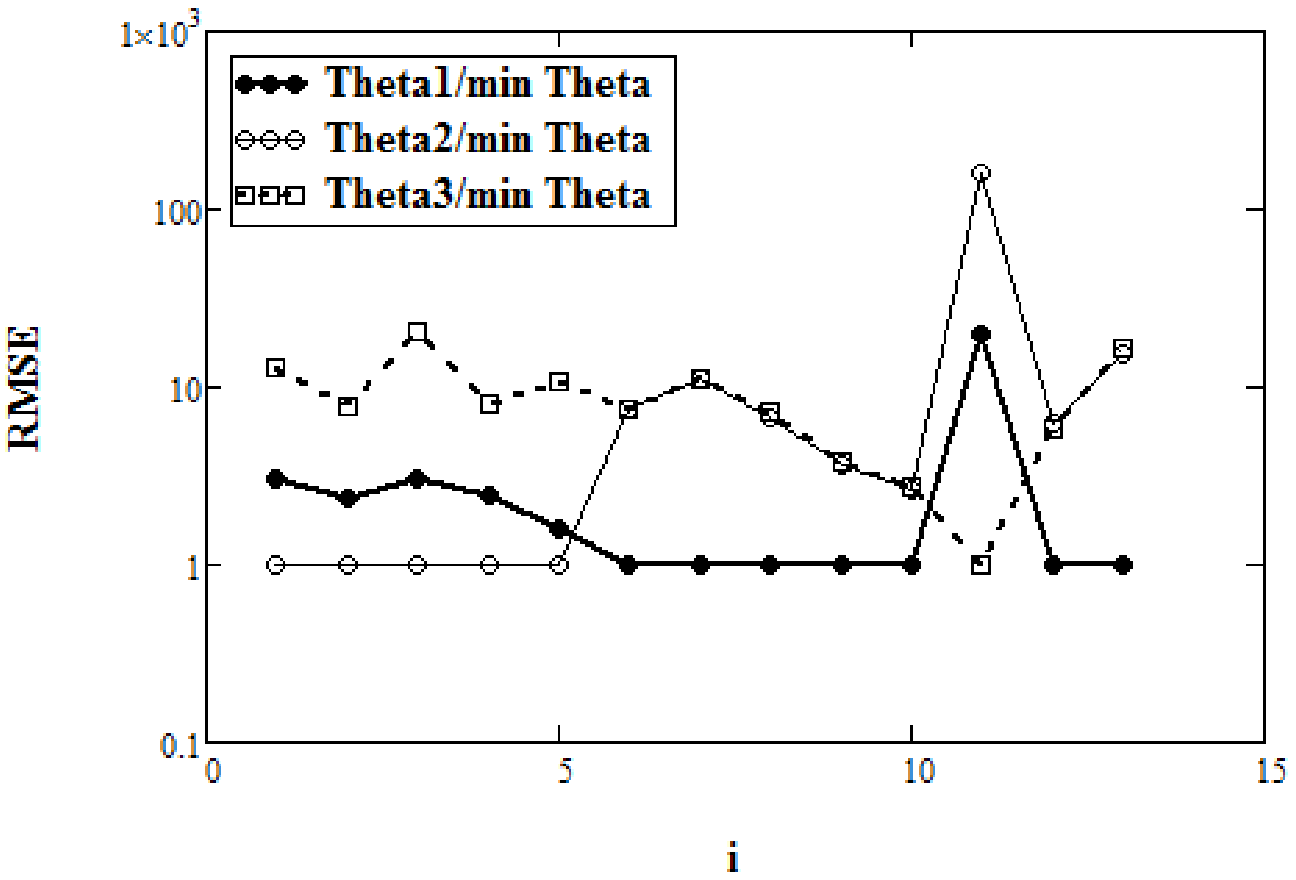}
 \end{minipage}
 \hfill
 \begin{minipage}{0.98\textwidth}
 \includegraphics[width=0.5\textwidth]{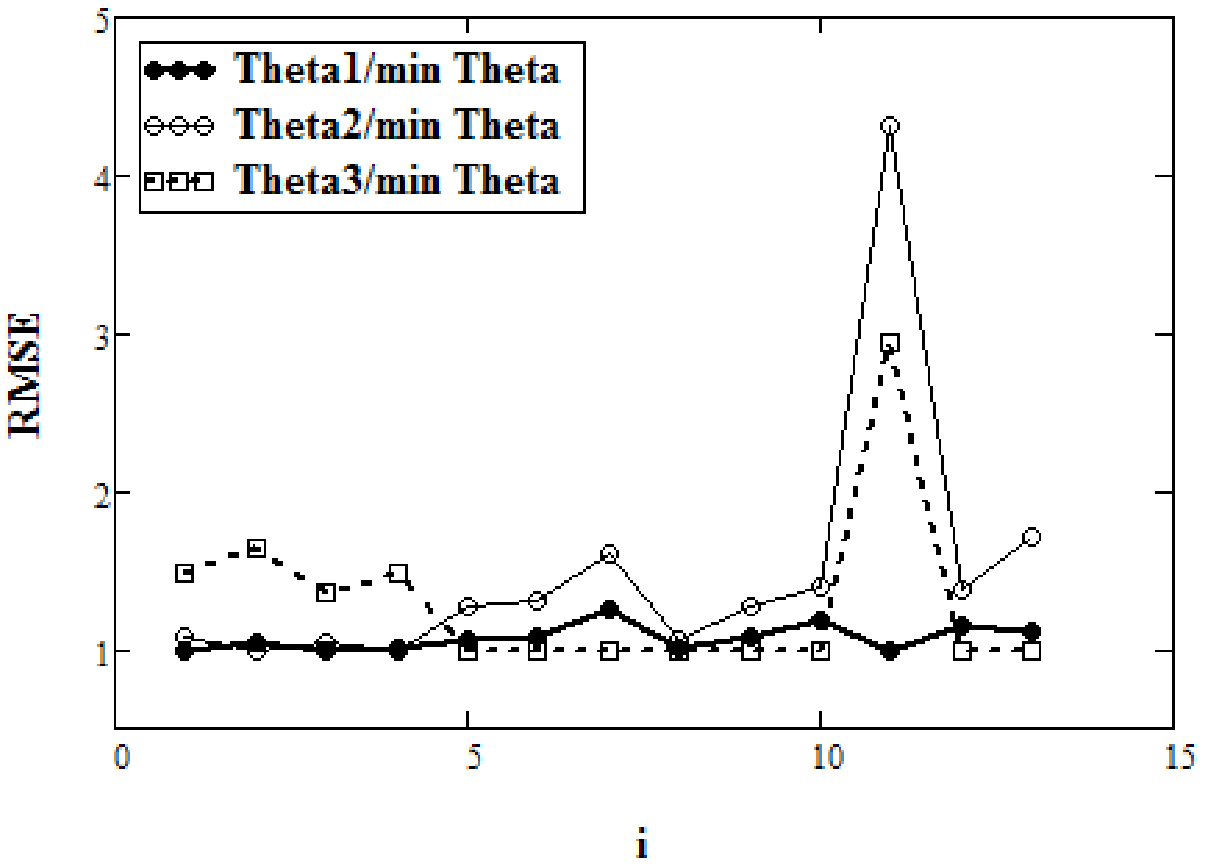}
\includegraphics[width=0.5\textwidth]{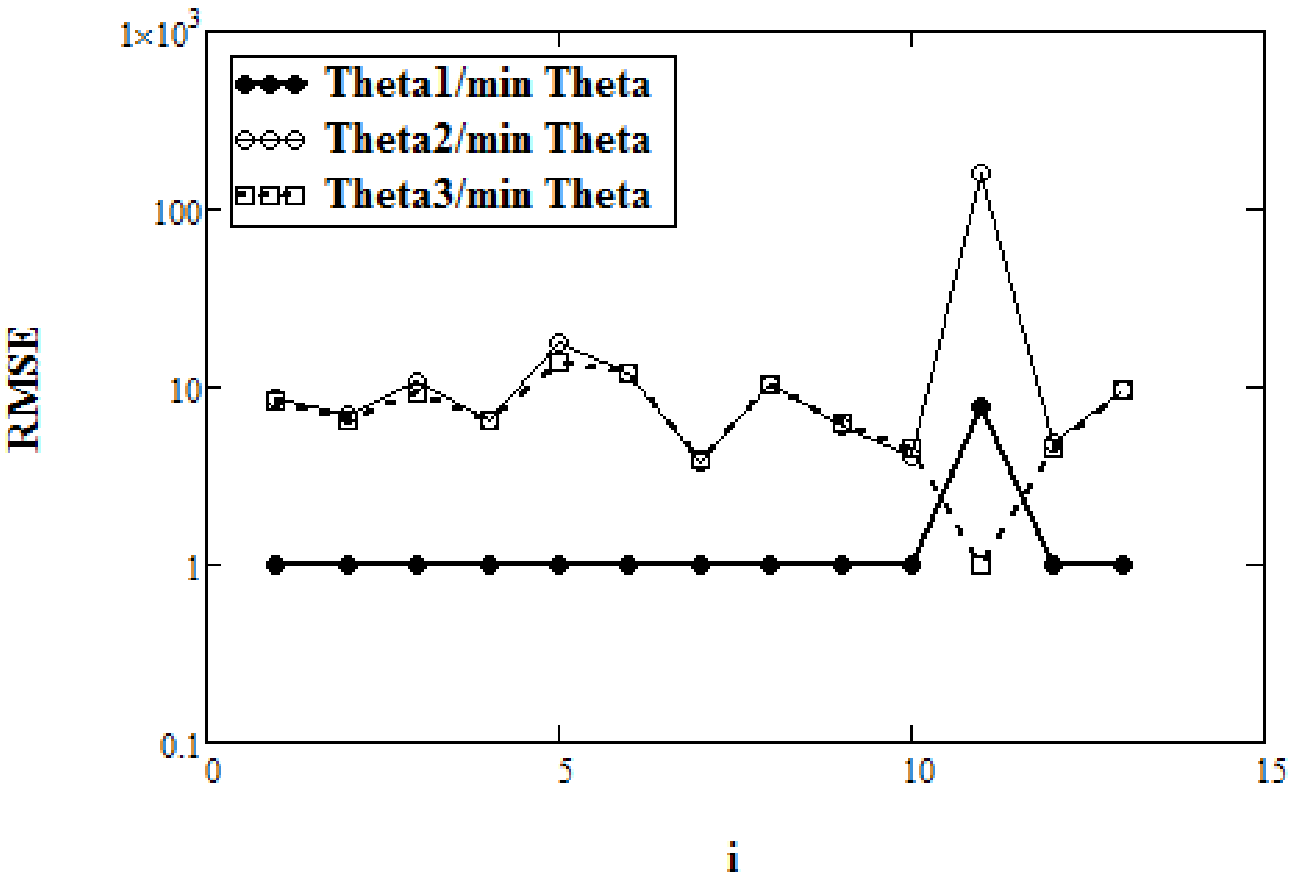}
 \end{minipage}
 \hfill
 \caption{Ratios of the RMSEs $RMSE(\widehat{\theta}_i)/\min_i(RMSE(\widehat{\theta}_i))$ corresponding to estimates $\{\widehat{\theta}_i\}$, $i\in\{1,2,3\}$ in (\ref{19})  for the intervals estimator (left column)
and the K-gaps estimator (right column)
obtained by the Algorithm with (\ref{5}) 
 against the number of processes related to the column labels in Tables \ref{Table1-N}, \ref{Table1-I} and \ref{Table1-K}: The upper figures correspond to  $k=\lfloor \widehat{\theta}_0 L\rfloor$ in Table \ref{Table1-N}, the middle figures to  $k=\lfloor\min(\widehat{\theta}_0 L, \sqrt{L})\rfloor$ in Table \ref{Table1-I} and the lower figures to  $k=\lfloor(\ln L)^2\rfloor$ in Table \ref{Table1-K}, for sample size $n=10^5$.}
 \label{fig:4}
\end{figure}
\begin{table}[t]
\centering
\caption{
The  root mean squared error 
of $\widehat{\theta}^{K_0}$  
($k=\lfloor s L\rfloor$), $\widehat{\theta}_0$ is a pilot intervals estimate.}
\tabcolsep=0.01cm
\begin{tabular}{|l|c c|c c|c c|c c|c c|c|c|c|c|}
  \hline
 {\scriptsize {$RMSE$}}         & \multicolumn{2}{c|}{\scriptsize {MM}} & \multicolumn{2}{c|}{\scriptsize {ARMAX}} & \multicolumn{2}{c|}{\scriptsize {$ARu^+$}} & \multicolumn{2}{c|}{\scriptsize {$ARu^-$}} &  \multicolumn{2}{c|}{\scriptsize {$MA(2)$}} & \scriptsize {ARc} & \scriptsize {AR(2)} & {\scriptsize {GARCH}}
 \\
 \cline{2-14}
 {\scriptsize $\cdot 10^{4}$/$\theta$}
&  $0.5$ & $0.8$            & 
$0.25$ & $0.75$               & $0.5$ & $0.8$                &   $0.75$ & $0.96$ & $0.5$ & $2/3$ &    $0.3$ & $0.25$ & $0.447$
\\
  \hline
 \multicolumn{14}{|c|}{$n=10^5$}
  \\
  \hline
$s=\theta$ &   &  &  &  &  &  &  &  &  & & &&
  \\
$\widehat{\theta}^{K0dis}_1$    & 9.018 & 11  & 8.879 & 10 &9.489  & 19 & 15 & 19  & 12 & 10 & 5.709 &9.104 &\textbf{16}
\\
$\widehat{\theta}^{K0dis}_2$  &8.701 &  11 & 8.741 & 9.985 & 9.042 & 19 & 14 & 19 & 12 & 10 & 5.709 & 9.081 &\textbf{16}
\\
$\widehat{\theta}^{K0dis}_3$   &9.866 &  12 & 9.264 & 11 & 10 & 19 & 15 & 19 & 12 & 10 & 5.709 & 9.141 &\textbf{16}
\\
$\widehat{\theta}^{K0dis*}_1$    &\textbf{\emph{3.419}} & \textbf{\emph{5.357}}  & \textbf{\emph{2.438}} & \textbf{\emph{5.016}} & \textbf{\emph{4.292}} & \textbf{\emph{9.341}} &\textbf{\emph{7.617}}  & \textbf{\emph{7.434}} & \textbf{\emph{5.507}} & \textbf{9.527} &\textbf{\emph{0.605}} &\textbf{\emph{4.696}} &19
\\
$\widehat{\theta}^{K0dis*}_2$  & \textbf{1.244} &  \textbf{1.959} & \textbf{1.073} & \textbf{1.859} & \textbf{1.681} & \textbf{5.412} & \textbf{3.712} & \textbf{2.938} & \textbf{2.643} & \textbf{\emph{9.569}} & \textbf{0.446} & \textbf{2.667} &\textbf{\emph{17}}
\\
$\widehat{\theta}^{K0dis*}_3$   &14 & 17  & 14 & 16 & 14 & 16 & 16 & 19 & 14 & 10 &5.709 & 11&22
\\
$s=\widehat{\theta}_0$  & &   &  &  &  &  &  &  &  &  & & &
\\
$\widehat{\theta}^{K0dis}_1$    & 161  & 213  & 154 & 196 & \emph{\textbf{251}} & \textbf{893} & \emph{\textbf{304}} & 386 & \textbf{196} & \textbf{376} & -&\textbf{357} &-
\\
$\widehat{\theta}^{K0dis}_2$  & 159 & 213  & 150 & 195 & \textbf{249} & \emph{\textbf{894}} & \emph{\textbf{304}} & 387 & \emph{\textbf{198}} & \emph{\textbf{377}} &- & \textbf{357}&-
\\
$\widehat{\theta}^{K0dis}_3$  & 168  &  216 & 165 & 202 & 256 & \textbf{893} & \textbf{303} & 386 & \textbf{196} & \textbf{376} &- &\textbf{357} &-
\\
$\widehat{\theta}^{K0dis*}_1$   &  \textbf{106} & \emph{\textbf{159}}  & \textbf{89} & \emph{\textbf{150}} & 299 & 931 & 525 & 329 & 229 &540 &14 &\emph{\textbf{401}} & \emph{\textbf{413}}
\\
$\widehat{\theta}^{K0dis*}_2$ & \emph{\textbf{126}}  & \textbf{157} & \emph{\textbf{103}} & \textbf{145} & 415 & 1072 & 690 & 361 & 324 & 734 &156 &424 & \textbf{404}
\\
$\widehat{\theta}^{K0dis*}_3$ & 348  & 443 & 293 & 445 & 419 & 989 & 511 & 350 & 335 & 463 &86 & 446& 474
\\
\hline
 \multicolumn{14}{|c|}{$n=5000$}
  \\
  \hline
  $s=\theta$ &   &  &  &  &  &  &  &  &  & & &&
  \\
$\widehat{\theta}^{K0dis}_1$    & 152 & 147  & 138 & 151 &162  & 189 & 199 & 211  & 173 & 254 & - &199 &106
\\
$\widehat{\theta}^{K0dis}_2$  &134 &  134 & 125 & 133 & 146 & 177 &  186 & 191 & 155 & 244 & - &186 & 106
\\
$\widehat{\theta}^{K0dis}_3$   &214 &  180 & 181 & 202 & 207 & 210 &  232 & 251 & 216 & 274 & - &231 &109
\\
$\widehat{\theta}^{K0dis*}_1$    & \textbf{\emph{60}} & \textbf{\emph{59}}  & \textbf{\emph{59}} & \textbf{\emph{51}} & \textbf{\emph{61}} & \textbf{\emph{60}} &\textbf{\emph{63}}  & \textbf{\emph{79}} & \textbf{\emph{65}} & \textbf{\emph{82}} & \textbf{\emph{65}} & \textbf{\emph{66}} & \textbf{\emph{85}}
\\
$\widehat{\theta}^{K0dis*}_2$  & \textbf{13} &  \textbf{16} & \textbf{12} & \textbf{15} & \textbf{16} & \textbf{16} & \textbf{18} & \textbf{22} & \textbf{16} & \textbf{27} & \textbf{0.87683} & \textbf{19} &\textbf{44}
\\
$\widehat{\theta}^{K0dis*}_3$   & 441 & 334  & 480 & 423 & 446 & 323 & 417 & 400 & 447 & 428 & 815 & 437 &236
\\
$s=\widehat{\theta}_0$  & &   &  &  &  &  &  &  &  &  & & &
\\
$\widehat{\theta}^{K0dis}_1$    & 1035  & 608  & 579 & 900 & 799 & 1444 & \emph{\textbf{1275}} & \emph{\textbf{391}} & 603 & \emph{\textbf{886}} & - & 1143 & 1621
\\
$\widehat{\theta}^{K0dis}_2$  & 1003 & 589  & 562 & 867 & \emph{\textbf{779}} & \emph{\textbf{1398}} & \textbf{1235} & 393 & 607 & 908 &- & 1137 & 1625
\\
$\widehat{\theta}^{K0dis}_3$  & 1117  &  707 & 690 & 1006 & 881 & 1503 & 1301 & 392 & 672 & 928 &- &1160 &1622
\\
$\widehat{\theta}^{K0dis*}_1$   &  \textbf{359} & \emph{\textbf{484}}  & \emph{\textbf{366}} & \emph{\textbf{469}} & \textbf{732} & \textbf{1310} & 1360 & \textbf{324} & \textbf{352} &\textbf{712} &\emph{\textbf{550}} &\emph{\textbf{949}} & \emph{\textbf{981}}
\\
$\widehat{\theta}^{K0dis*}_2$ & \emph{\textbf{362}}  & \textbf{466} & \textbf{299} & \textbf{411} & 803 & 1606 & 1708 & 395 & \emph{\textbf{574}} & 1174 & \textbf{422} & \textbf{848} & \textbf{851}
\\
$\widehat{\theta}^{K0dis*}_3$ & 1765  & 1573 & 1578 & 1729 & 2004 & 1722 & 1843 & 679 & 1684 & 1725 & 2242 & 1833 & 2183
\\
\hline
\end{tabular}
 \label{Table1-KgapsK0Theta}
\end{table} 
\begin{table}[t]
\centering
\caption{
The  absolute bias of $\widehat{\theta}^{K_0}$  
($k=\lfloor s L\rfloor$), $\widehat{\theta}_0$ is a pilot intervals estimate.}
\tabcolsep=0.01cm
\begin{tabular}{|l|c c|c c|c c|c c|c c|c|c|c|c|}
  \hline
    {\scriptsize $|Bias|$}         & \multicolumn{2}{c|}{MM} & \multicolumn{2}{c|}{ARMAX} & \multicolumn{2}{c|}{$ARu^+$} & \multicolumn{2}{c|}{$ARu^-$} &  \multicolumn{2}{c|}{$MA(2)$} & ARc & AR(2) & {\scriptsize GARCH}
 \\
 \cline{2-14}
 {\scriptsize $\cdot 10^{4}$/$\theta$}
&  $0.5$ & $0.8$            & 
$0.25$ & $0.75$               & $0.5$ & $0.8$                &   $0.75$ & $0.96$ & $0.5$ & $2/3$ &    $0.3$ & $0.25$ & $0.447$
\\
  \hline
   \multicolumn{14}{|c|}{$n=10^5$}
  \\
  \hline
$s=\theta$&   &  &  &  &  &  &  &  &  & & &&
\\
$\widehat{\theta}^{K0dis}_1$    &5.824 &  10 & 2.935 & 9.194 & 6.493 & 19 & 13 & 19 & 8.306 &\textbf{6.956}  & 5.709& 1.425 &\textbf{\emph{16}}
\\
$\widehat{\theta}^{K0dis}_2$  & 5.427 & 9.432  & 2.818 & 8.703 & 6.086 & 19 & 13 & 19 & 8.065 &\textbf{6.956}  & 5.709&1.470 &\textbf{\emph{16}}
\\
$\widehat{\theta}^{K0dis}_3$   &6.299 & 11  & 3.041 & 9.727 & 6.899 & 19 & 14 & 19 & 8.547 & \textbf{6.956} &5.709 & 1.381 &\textbf{\emph{16}}
\\
$\widehat{\theta}^{K0dis*}_1$    & \textbf{\emph{2.789}} & \textbf{\emph{5.198}}  & \textbf{\emph{1.167}} & \textbf{\emph{4.799}} & \textbf{\emph{3.495}} & \textbf{\emph{8.924}} & \textbf{\emph{7.232}} & \textbf{\emph{7.409}} & \textbf{\emph{4.402}} & 8.053 &\textbf{\emph{0.605}} & \textbf{\emph{1.284}} &16
\\
$\widehat{\theta}^{K0dis*}_2$  &\textbf{0.859} & \textbf{1.851}  & \textbf{0.298} & \textbf{1.708} & \textbf{1.161} & \textbf{4.739} & \textbf{3.280} & \textbf{2.914} & \textbf{1.876} & 8.286 &\textbf{0.446} & \textbf{0.519} &\textbf{13}
\\
$\widehat{\theta}^{K0dis*}_3$   & 9.880&  16 & 5.288 & 15 & 10 & 15 & 15 & 19 & 10 & \textbf{\emph{7.017}} & 5.709 & 2.693 &20
\\
$s=\widehat{\theta}_0$  & &   &  &  &  &  &  &  &  &  & & &
\\
$\widehat{\theta}^{K0dis}_1$    & \textbf{\emph{16}}  &  \textbf{\emph{19}} & \textbf{\emph{18}} & 41 & \textbf{\emph{135}} & \textbf{\emph{876}} & \textbf{182} & 382 & 125 & 313 &- & \textbf{276} &-
\\
$\widehat{\theta}^{K0dis}_2$  & \textbf{12} &  \textbf{15} & \textbf{14} & \textbf{\emph{37}} & 137 & 878 & \textbf{\emph{183}} & 383 & 127 & 315 &- & \textbf{276}&-
\\
$\widehat{\theta}^{K0dis}_3$  & 18  & 23  & 23 & 44 & \textbf{134} & \textbf{875} & \textbf{182} & 382 &\textbf{\emph{123}} &\textbf{\emph{312}} & -&\textbf{276} &-
\\
$\widehat{\theta}^{K0dis*}_1$& 28 & 15\textbf{}  & 28 & \textbf{7.082} & 268 & 919 & 481 &\textbf{\emph{315}} & 200 & 510 &\textbf{14} &357 &\textbf{\emph{177}}
\\
$\widehat{\theta}^{K0dis*}_2$ &  76 & 38  & 76 & 52 & 402 & 1065 & 668 & 348 & 302 & 707 &156 & 393&193
\\
$\widehat{\theta}^{K0dis*}_3$  & 57  &  65 & 46 & 72 & 157 & 889 & 293 & \textbf{259} & \textbf{49} & \textbf{246} & \textbf{\emph{86}} &\textbf{\emph{318}} &\textbf{161}
\\
\hline
\multicolumn{14}{|c|}{$n=5000$}
  \\
  \hline
  $s=\theta$&   &  &  &  &  &  &  &  &  & & &&
\\
$\widehat{\theta}^{K0dis}_1$   &101 &117 & 76 & 98 &112 & 162 & 138 & 182 & 119 &192& -  & 125 & \textbf{\emph{12}}
\\
$\widehat{\theta}^{K0dis}_2$  & 75 & 97 & 57 & 74 & 91 & 145 & 118 & 153 & 97 & 176  & - &105 & 13
\\
$\widehat{\theta}^{K0dis}_3$   &142 & 145  & 104&128 & 141 & 182 & 162 & 217 & 147 & 207 & - & 148 &\textbf{11}
\\
$\widehat{\theta}^{K0dis*}_1$    &  \textbf{\emph{56}} & 58 & 51  & 49 & \textbf{\emph{57}} & \textbf{\emph{59}} & \textbf{\emph{60}} & \textbf{\emph{79}} & 61 & \textbf{\emph{78}} & \textbf{\emph{65}} &\textbf{\emph{56}} & 58
\\
$\widehat{\theta}^{K0dis*}_2$  &\textbf{8.9597} & \textbf{16}  & \textbf{3.3883} & \textbf{14} & \textbf{11} & \textbf{16} & \textbf{14} & \textbf{22} & \textbf{12} & \textbf{22} &\textbf{0.87683} & \textbf{12} & 27
\\
$\widehat{\theta}^{K0dis*}_3$   & 405&  326 & 410 & 413 & 409 & 316 & 407 & 400 & 412 & 410 & 815 & 372 &102
\\
$s=\widehat{\theta}_0$  & &   &  &  &  &  &  &  &  &  & & &
\\
$\widehat{\theta}^{K0dis}_1$    & 436  &  202 & 149 & 332 & \textbf{\emph{384}} & 1312 & \textbf{\emph{750}} & 370 & \textbf{\emph{84}} & \textbf{\emph{163}} & - & \textbf{\emph{601}} & \textbf{\emph{328}}
\\
$\widehat{\theta}^{K0dis}_2$  & 365  &  \textbf{\emph{161}} & \textbf{\emph{111}} & 276 & \textbf{370} & \textbf{\emph{1255}} & \textbf{720} & 370 & \textbf{31} & 219 & - & \textbf{577} & \textbf{319}
\\
$\widehat{\theta}^{K0dis}_3$  & 505  & 243  & 203 & 400 & 404 & 1355 & 782 & 371 & 143 & \textbf{110} & - & 625 & 337
\\
$\widehat{\theta}^{K0dis*}_1$& \textbf{\emph{119}} & 168  & 124 & 153 & 576 & \textbf{1252} & 1236 & \textbf{\emph{279}} & \textbf{\emph{153}} & 571 & \textbf{\emph{550}} & 777 & 410
\\
$\widehat{\theta}^{K0dis*}_2$ &  \textbf{81} & \textbf{29}  & \textbf{49} & \textbf{37} & 684 & 1557 & 1619 & 386 & 466 & 1092 & \textbf{422} & 688 & 386
\\
$\widehat{\theta}^{K0dis*}_3$  & 1015  &  810 & 946 & 1053 & 1184 & 1421 & 1143 & \textbf{152} & 946 & 810 & 2242 & 1256 & 1086
\\
\hline
\end{tabular}
 \label{Table2-KgapsK0Theta}
\end{table} 
\begin{table}[t]
\centering
\caption{The  root mean squared error ($k=\lfloor \widehat{\theta}_0L\rfloor$), $\widehat{\theta}_0$ is a pilot intervals estimate.}
\tabcolsep=0.09cm
\begin{tabular}{|l|c c|c c|c c|c c|c c|c|c|c|c|}
  \hline
 {\scriptsize $RMSE$}         & \multicolumn{2}{c|}{MM} & \multicolumn{2}{c|}{ARMAX} & \multicolumn{2}{c|}{$ARu^+$} & \multicolumn{2}{c|}{$ARu^-$} &  \multicolumn{2}{c|}{$MA(2)$} & ARc & AR(2) & {\scriptsize GARCH}
 \\
 \cline{2-14}
 {\scriptsize $\cdot 10^{4}$/$\theta$}
&  $0.5$ & $0.8$            & 
$0.25$ & $0.75$               & $0.5$ & $0.8$                &   $0.75$ & $0.96$ & $0.5$ & $2/3$ &    $0.3$ & $0.25$ & $0.447$
\\
  \hline
  \multicolumn{14}{|c|}{$n=10^5$}
  \\
  \hline
$\widehat{\theta}_1$   & 147& 215  & 159 & 211 & 230 & 887 & 287 &383&199&\textbf{\emph{400}}  &-  & \textbf{305} &-
\\
$\widehat{\theta}_2$    & 146& 213  & 158 & 211 & 230 & 889 & 287 &384&201& 402 &-  &\textbf{305} &-
\\
$\widehat{\theta}_3$   & 154&  222 & 164 & 215 & 235 & 886 & 287&383&203 & \textbf{\emph{400}}  & - & \textbf{305} &-
\\
$\widehat{\theta}_1^*$   & \textbf{103}& 160  & 89 & \textbf{\emph{151}} &292  & 928 & 516 &\textbf{\emph{328}}& 229 & 542 &\textbf{\emph{17}}  & 400 &413
\\
$\widehat{\theta}_2^*$   &123 & \textbf{\emph{156}}  & 99 & \textbf{\emph{151}} & 410 & 1070 & 690 &354& 324 & 745 & 155 &420  &405
\\
$\widehat{\theta}_3^*$  &354 & 442  & 294 &425  & 413 & 957 & 516 &351& 336 & 467 & 67 & 443 &470
\\
$\widehat{\theta}^{Kdis}_1$    & 133& 204  & 145 & 190 & 195 & 772 & 209 &349 & \textbf{141} &423  & - & 390 &\textbf{375}
\\
$\widehat{\theta}^{Kdis}_2$    & 777&  844 & 648 & 940 & 799 & 777 & \textbf{\emph{208}} & 1396 & 799& 420 & - &519  &\textbf{\emph{382}}
\\
$\widehat{\theta}^{Kdis}_3$    &895 & 702  & 651 & 821 & 398 & \textbf{\emph{763}} & \textbf{207} & 1115 &927 & 637 & - & 477 &389
\\
$\widehat{\theta}^{Kdis*}_1$    & 136&  264 & 100 & 237 &\textbf{150}  & \textbf{631} & 220 & \textbf{105} & 229 & \textbf{390}  & \textbf{12} & 394 &464
\\
$\widehat{\theta}^{Kdis*}_2$    & 127&  226 & \textbf{\emph{70}} & 207 &249 &788  & 5189 & 404 & \textbf{\emph{152}} &1944  &34  &1231  &3929
\\
$\widehat{\theta}^{Kdis*}_3$  & 652&  897 & 238 & 1013 &690  & 1300 & 963 &345&798& 958 & 75 &\textbf{\emph{353}}  &573
\\
\hline
  & &   &  &  &  &  &  &  &  &  &&&
\\
$\widehat{\theta}^{Kimt}$  & 217 & 569  & \textbf{69} & 498 & \textbf{\emph{173}} 
& 844 
& 2501 & 401 & 309 & 466 & 33 & 3630 & 4028
\\
  & &   &  &  &  &  &  &  &  &  &&&
  \\
$\widehat{\theta}^{IA1}$  & \textbf{\emph{116}} &  \textbf{122} & 95 & \textbf{113}  & 447 & 1193 
& 1756 & 399 & 387 & 977 & 233  & 693 & 580
\\
  \hline
\multicolumn{14}{|c|}{$n=5000$}
  \\
  \hline
$\widehat{\theta}_1$   &565 & 938  &506  & 818 & 783 & 1431 & 1364 &394& 533&816& - & 900 &  1497
\\
$\widehat{\theta}_2$    &557 & 913  & 476 & 787 & 748 & 1401 & 1337 &396&537&810& - & 897 &  1497
\\
$\widehat{\theta}_3$   &633 &  1013 & 620 & 903 & 879 & 1490 & 1416 &394&593 &850& - & 910 &  1497
\\
$\widehat{\theta}_1^*$   & 359 & \textbf{\emph{496}}  & 350 & \textbf{\emph{464}} & 715 & 1294 & 1276 &\textbf{\emph{315}}&352& 760& 606 & \textbf{\emph{835}} & \textbf{\emph{955}}
\\
$\widehat{\theta}_2^*$   & 352 & \textbf{466}  & 291 & \textbf{450} & 808 & 1587 & 1666 &395&557 &1179& 422 & \textbf{794} &  \textbf{870}
\\
$\widehat{\theta}_3^*$  & 1635 & 1564  & 1377 & 1652 & 1902 & 1644 &1754  &713&1505 &1656& 1455 & 1610 &  2036
\\
$\widehat{\theta}^{Kdis}_1$    &480 &  917 & 496 & 772 &  787& 1186 & 1820 & 427&406& 807 & - & 1690 &1877
\\
$\widehat{\theta}^{Kdis}_2$    & 1525& 1880  & 982 & 1793 & 1836 & 1863 & 2672 &2981&1218& 2020 &-  &1855  &2929
\\
$\widehat{\theta}^{Kdis}_3$   &1624 & 1993 & 1286 & 1815 &1692  & 1686 & 2348 &1775&1924& 2349 &-  &1737  &2337
\\
$\widehat{\theta}^{Kdis*}_1$    & 320& 605  &299  & 507 &\textbf{\emph{453}} & \textbf{592} & \textbf{641} &\textbf{213}&404&754  &\textbf{\emph{72}}  &1528  &1491
\\
$\widehat{\theta}^{Kdis*}_2$   & \textbf{\emph{252}}& 548  & \textbf{\emph{199}} & 487 & 535 & \textbf{\emph{866}} & 2529 &423&\textbf{\emph{335}}& \textbf{\emph{488}} &\textbf{25}  &3684  &3787
\\
$\widehat{\theta}^{Kdis*}_3$  & 824& 1007  & 931 &871  & 1086 & 927 & \textbf{\emph{916}} &555&714&929  & 2321 &1106  &1324
\\
\hline
  & &   &  &  &  &  &  &  &  &  &&&
\\
$\widehat{\theta}^{Kimt}$  & \textbf{247} & 588 & \textbf{188} & 525 & \textbf{293} & 869 & 2518 & 418 & \textbf{325} & \textbf{474} & \textbf{25} & 3680 & 3900
\\
  & &   &  &  &  &  &  &  &  &  &&&
  \\
$\widehat{\theta}^{IA1}$  & 385 & 513  & 319 &478  & 694 &1388  &1985  &394  &514  & 1114 &676&980&1077
\\
\hline
\end{tabular}
 \label{Table1-N}
\end{table} 
\begin{table}[t]
\centering
\caption{The  absolute bias ($k=\lfloor \widehat{\theta}_0L\rfloor$), $\widehat{\theta}_0$ is a pilot intervals estimate.}
\tabcolsep=0.033cm
\begin{tabular}{|l|c c|c c|c c|c c|c c|c|c|c|c|}
  \hline
 {\scriptsize $|Bias|$}         & \multicolumn{2}{c|}{MM} & \multicolumn{2}{c|}{ARMAX} & \multicolumn{2}{c|}{$ARu^+$} & \multicolumn{2}{c|}{$ARu^-$} &  \multicolumn{2}{c|}{$MA(2)$} & ARc & AR(2) & {\scriptsize GARCH}
 \\
 \cline{2-14}
 {\scriptsize $\cdot 10^{4}$/$\theta$}
&  $0.5$ & $0.8$            & 
$0.25$ & $0.75$               & $0.5$ & $0.8$                &   $0.75$ & $0.96$ & $0.5$ & $2/3$ &    $0.3$ & $0.25$ & $0.447$
\\
  \hline
  \multicolumn{14}{|c|}{$n=10^5$}
  \\
  \hline
$\widehat{\theta}_1$   & 18& \textbf{\emph{2.515}}  & 20 & 21 & 142 & 872 & 203 & 380 & 132 & 333 & - & \textbf{217} &-
\\
$\widehat{\theta}_2$    &16 & \textbf{2.498}  & \textbf{16} & \textbf{\emph{20}} & 142 & 874 & 203 &381&134&335  & - & \textbf{217} &-
\\
$\widehat{\theta}_3$   & 20&  8.680 & 24 & 21 & 141 & 870 &203  &379&131&\textbf{\emph{331}}  &-  & \textbf{217} &-
\\
$\widehat{\theta}_1^*$   &29& 6.723  & 28& \textbf{14} & 259 & 915 &469  &313&200& 513 & \textbf{\emph{17}} & 354 &184
\\
$\widehat{\theta}_2^*$   &73 & 43  & 72 & 55 & 396 & 1064 & 669 &336&300& 719 & 155 & 386 &198
\\
$\widehat{\theta}_3^*$  & \textbf{\emph{11}}& 42  & 66 & 63 & 132 & 858 & 254 & 253&\textbf{\emph{62}}& \textbf{261} & 67 &320  &170
\\
$\widehat{\theta}^{Kdis}_1$    & 24&  66 & \textbf{\emph{18}} & 43 & \textbf{\emph{29}} & 757 & 48 &341&\textbf{60}&383  &-  &323  &\textbf{14}
\\
$\widehat{\theta}^{Kdis}_2$    &133 &  138 & 172 & 145 &93 & 763 & \textbf{\emph{44}} &\textbf{\emph{175}}&161&380  &-  &284  &\textbf{\emph{79}}
\\
$\widehat{\theta}^{Kdis}_3$   &188 & 133  & 181 & 135 & \textbf{1.3951} &747  & 50&185&234&423  &-  &271  &92
\\
$\widehat{\theta}^{Kdis*}_1$    & 103&  212 &40  & 194 &30  &\emph{\textbf{620}} & \textbf{43} &\textbf{90}&217&368  &\textbf{12}  &364  &346
\\
$\widehat{\theta}^{Kdis*}_2$   & 115& 204  &51 & 190 & 56 & 786 & 3497 &404&137&448  &34  &1085  &3020
\\
$\widehat{\theta}^{Kdis*}_3$  & 225 & 488  &34 & 485 & 141 & \textbf{136} & 306 &321&359&636  &75  &\textbf{\emph{222}}  &249
\\
\hline
  & &   &  &  &  &  &  &  &  &  &&&
\\
$\widehat{\theta}^{Kimt}$  & \textbf{0.148} & 567  & 54  & 496 & 165 
& 843 
& 2501 & 401 & 306 & 462 & 33 & 3627 & 4027
\\
  & &   &  &  &  &  &  &  &  &  &&&
\\
$\widehat{\theta}^{IA1}$  & 82 & 45  & 64 & 54 & 436 & 1187 
& 1752 & 399 &378 & 972 & 233 & 687 & 563
\\
  & &   &  &  &  &  &  &  &  &  &&&
  \\
\hline
\multicolumn{14}{|c|}{$n=5000$}
  \\
  \hline
$\widehat{\theta}_1$   &144 & 336  & 100 &255  & 339 & 1284 & 851 &373&46& 242 & - &345  &\textbf{\emph{11}}
\\
$\widehat{\theta}_2$    &103 & 290  & 54 &205  & 317 & 1245 & 811 &376&\textbf{\emph{8.581}}& 265 & - & 345 &\textbf{\emph{11}}
\\
$\widehat{\theta}_3$   &185 &  391 & 148 & 305 & 360 & 1322 & 891 &370&81& 221 &-  & \textbf{\emph{344}} &\textbf{\emph{11}}
\\
$\widehat{\theta}_1^*$   &99 & 142  & 109 & 112 & 560 & 1236 & 1150 &272&170&631  &606  & 665 &382
\\
$\widehat{\theta}_2^*$   &82 & \textbf{\emph{14}}  & 51 & 35 & 681 & 1534 & 1574 &385&450& 1094 &422  & 636 &378
\\
$\widehat{\theta}_3^*$  &812 & 779  & 719 & 754 & 992 & 1192 & 1021 &114&735&658  & 1455 & 1047 &842
\\
$\widehat{\theta}^{Kdis}_1$    &\textbf{\emph{16}} &  218 & 58 & 98 & 115 & 1048 & 1260 &356&41&\textbf{\emph{109}}  &-  &1186  &1254
\\
$\widehat{\theta}^{Kdis}_2$    & 437&  158 & 290 & 269 & 454 & 727 & 926 &421&305&332  &-  &853  &379
\\
$\widehat{\theta}^{Kdis}_3$   &467 & 165  & 467 & 238 & 341 & 838 & 858 &\textbf{\emph{80}}&715&564  &-  &635  &154
\\
$\widehat{\theta}^{Kdis*}_1$    & 169&  343 & 34 & 309 & \textbf{\emph{108}} & \textbf{\emph{480}} & 467 &135&343&680  &72  &1458  &1342
\\
$\widehat{\theta}^{Kdis*}_2$   & 192&  492 & 43 & 443 & 145 & 852 & 2529 &423&290&429  &\textbf{\emph{25}}  &3320  &3781
\\
$\widehat{\theta}^{Kdis*}_3$  & 129&  146 & 295 & 101& 304 & 534 & \textbf{\emph{157}} &110&25&306  &2321  &582  &241
\\
\hline
  & &   &  &  &  &  &  &  &  &  &&&
\\
$\widehat{\theta}^{Kimt}$  & \textbf{2.636}& \textbf{4.447}  & \textbf{3.448} &\textbf{3.598}  & \textbf{5.124} &\textbf{1.988}  &\textbf{1.836}  & \textbf{ 1.642}&\textbf{2.544} & \textbf{5.331} &\textbf{1.492}& \textbf{2.527}&\textbf{7.127}
\\
  & &   &  &  &  &  &  &  &  &  &&&
\\
$\widehat{\theta}^{IA1}$  &36 &  29 &  \textbf{\emph{31}}& \textbf{\emph{20}} &523  &1301  &1919  &385  &361  & 1033 &676&884&876
  \\
\hline
\end{tabular}
 \label{Table2-N}
\end{table} 
\begin{table}[t]
\centering
\caption{The  root mean squared error ($k=\lfloor\min(\widehat{\theta}_0L, \sqrt{L})\rfloor$), $\widehat{\theta}_0$ is a pilot intervals estimate.}
\tabcolsep=0.09cm
\begin{tabular}{|l|c c|c c|c c|c c|c c|c|c|c|c|}
  \hline
 {\scriptsize $RMSE$}         & \multicolumn{2}{c|}{MM} & \multicolumn{2}{c|}{ARMAX} & \multicolumn{2}{c|}{$ARu^+$} & \multicolumn{2}{c|}{$ARu^-$} &  \multicolumn{2}{c|}{$MA(2)$} & ARc & AR(2) & {\scriptsize GARCH}
 \\
 \cline{2-14}
 {\scriptsize $\cdot 10^{4}$/$\theta$}
&  $0.5$ & $0.8$            & 
$0.25$ & $0.75$               & $0.5$ & $0.8$                &   $0.75$ & $0.96$ & $0.5$ & $2/3$ &    $0.3$ & $0.25$ & $0.447$
\\
  \hline
   \multicolumn{14}{|c|}{$n=10^5$}
  \\
  \hline
$\widehat{\theta}_1$   & \textbf{\emph{135}}&  171 &120  &164  &405  &1059  &1241  &\textbf{\emph{338}}&348& 846 & \textbf{14} & 611 &479
\\
$\widehat{\theta}_2$    & 146&169   &120  &162  &471  &1150  &1524  &340&414& 1193 & 137 & 721 &700
\\
$\widehat{\theta}_3$   & 191& 245  &171  &246  &383  &997  &1029  &343&324&\textbf{467}  & 245 & 533 &427
\\
$\widehat{\theta}_1^*$   &\textbf{101} &\textbf{118}   &\textbf{\emph{85}}  &\textbf{115}  &388  &1116  &1310  &\textbf{329}&331& 834 &151  &620  &416
\\
$\widehat{\theta}_2^*$  & 155&\textbf{\emph{130}}   &116  &\textbf{\emph{142}}  &597  &1527  & 2374 &370&499& 991 & 358 & 953 &1124
\\
$\widehat{\theta}_3^*$ &359 &444   &298  &433  &399  &964  & 518 &349&348& 701 & 67 & 431 &476
\\
  $\widehat{\theta}^{Kdis}_1$   & 384& 738  & 155 & 670 & \textbf{188} & \textbf{384} &975  & 672 &537 & 1015 &\textbf{\emph{19}}&382&\textbf{\emph{322}}
\\
$\widehat{\theta}^{Kdis}_2$    & 3222& 4566  & 1672 & 4317 & 3221 & 4457 & 4252 & 4334 & 3158 & 4232 &3000&1697&3024
\\
$\widehat{\theta}^{Kdis}_3$   & 2821&4372   &1398  & 4061& 2517 &  4411& 4273 & 4596 & 3255 & 4385 &\textbf{\emph{19}}&1589&3192
\\
$\widehat{\theta}^{Kdis*}_1$    &673 &1361   & 207 & 1232 & 381 & 591 & \textbf{379} & 642& 865 &1577 &380&\textbf{273}&\textbf{195}
\\
$\widehat{\theta}^{Kdis*}_2$   & 220&572   & \textbf{67} & 501 & \textbf{\emph{237}} & 843& 2501 &404 & \textbf{\emph{315}} &\textbf{\emph{517}} &34&2531&4176
\\
$\widehat{\theta}^{Kdis*}_3$  &223 &536   & 194 & 449 &240  & \textbf{\emph{565}} & \textbf{\emph{397}} & 5588 & \textbf{303} & 596 &75&\textbf{\emph{312}}&398
\\
  \hline
  \multicolumn{14}{|c|}{$n=5000$}
  \\
  \hline
$\widehat{\theta}_1$   &535 &  636 &471  &625  &813  &1260  &1499  &424&507& \textbf{\emph{786}} & - & 943 &992
\\
$\widehat{\theta}_2$    &499 &620   &413  &593  &795  &1328  &1754  &428&560& 1203 &-  & 1000 &1122
\\
$\widehat{\theta}_3$   &918 &1008   &823  &1004  &1152  &1302  &1460  &490&828& 1655 &-  & 1093 &1242
\\
$\widehat{\theta}_1^*$   &\textbf{\emph{365}} & \textbf{\emph{499}}  &343  &\textbf{\emph{483}}  &714  & 1290 &1445  &\textbf{313}&\textbf{\emph{357}}& 936 & 606 & 897 &891
\\
$\widehat{\theta}_2^*$  & 368& \textbf{496}  &292  &\textbf{458}  &815  &1604  &2260  &\textbf{\emph{389}}&577&1104  & 422 & 1120 &1359
\\
$\widehat{\theta}_3^*$ &1668 &1555   &1359  &1697  &1887  &1662  &1716  &721&1570& 1043 & 1455 & 1649 &2030
\\
  $\widehat{\theta}^{Kdis}_1$   &598 &1147   & 354 & 1028 & \textbf{\emph{454}} & \textbf{633} & \textbf{746} & 954 & 725 &1283  &77&\textbf{\emph{844}}&\textbf{\emph{808}}
\\
$\widehat{\theta}^{Kdis}_2$    & 3437& 5398  & 1705 & 5075 & 3412 & 5268 & 5266 & 5742 & 3543 &4689  &3000&2133&4032
\\
$\widehat{\theta}^{Kdis}_3$   & 3054&5158   & 1399 & 4744 & 2832 &5335 & 5137 & 5513 & 3163 & 4776 &\textbf{11}&1859&3553
\\
$\widehat{\theta}^{Kdis*}_1$    &792 &1909   &  \textbf{\emph{276}}& 1712 & 582 & 1253 &929 &1044  &  959& 1821 &93&\textbf{836}&\textbf{377}
\\
$\widehat{\theta}^{Kdis*}_2$   & \textbf{251}&585   & \textbf{188} & 509 & \textbf{288} & \textbf{\emph{885}} & 2522 &423  &\textbf{341} & \textbf{498} &\textbf{\emph{25}}&3528&3775
\\
$\widehat{\theta}^{Kdis*}_3$  & 793&958   & 981 &876  & 1090 & 905 & \textbf{\emph{819}} &554  & 725 & 929 &2321&1095&1404
\\
  \hline
\end{tabular}
 \label{Table1-I}
\end{table} 
\begin{table}[t]
\centering
\caption{The  absolute bias ($k=\lfloor\min(\widehat{\theta}_0L, \sqrt{L})\rfloor$), $n=10^5$.}
\tabcolsep=0.035cm
\begin{tabular}{|l|c c|c c|c c|c c|c c|c|c|c|c|}
  \hline
 {\scriptsize $|Bias|$}         & \multicolumn{2}{c|}{MM} & \multicolumn{2}{c|}{ARMAX} & \multicolumn{2}{c|}{$ARu^+$} & \multicolumn{2}{c|}{$ARu^-$} &  \multicolumn{2}{c|}{$MA(2)$} & ARc & AR(2) & {\scriptsize GARCH}
 \\
 \cline{2-14}
 {\scriptsize $\cdot 10^{4}$/$\theta$}
&  $0.5$ & $0.8$            & 
$0.25$ & $0.75$               & $0.5$ & $0.8$                &   $0.75$ & $0.96$ & $0.5$ & $2/3$ &    $0.3$ & $0.25$ & $0.447$
\\
  \hline
  \multicolumn{14}{|c|}{$n=10^5$}
  \\
  \hline
$\widehat{\theta}_1$   & 61&36   &41  &\textbf{\emph{30}}  &356  &1036  &1147  &302&316& 841 & \textbf{14} &581  &315
\\
$\widehat{\theta}_2$    & 84& 48  &60  &49  &423  &1122  &1413  &304&385& 1190 &137  &687  &547
\\
$\widehat{\theta}_3$   & \textbf{35}&\textbf{15}   &\textbf{\emph{21}}  &\textbf{5.4013}  &281  &957  &862  &299&\textbf{\emph{236}}& \textbf{245} & 245 &477  &\textbf{\emph{90}}
\\
$\widehat{\theta}_1^*$   & 59&37   &42  &43  &376  &1111  &1306  &318&323&798  & 151 &613  &390
\\
$\widehat{\theta}_2^*$  & 133&76   &97  &97  &589  &1522  &2372  &364&493&956  & 358 &949  &1116
\\
$\widehat{\theta}_3^*$ & \textbf{\emph{38}}&\textbf{\emph{35}}   &45  &42  &119  &865  &\textbf{250}  &\textbf{251}&\textbf{80}&617  &67  &304  &166
\\
  $\widehat{\theta}^{Kdis}_1$   & 347&  677 & 122 & 615 & \textbf{\emph{80}} & \textbf{\emph{197}} & 429 & \textbf{\emph{276}}& 505 & 951 &\textbf{\emph{19}}&279&\textbf{5.5944}
\\
$\widehat{\theta}^{Kdis}_2$    &2177 &2874   & 1143 & 2717 & 2012 & 1932 & 562 & 1868 & 2143 & 2938 &3000&572&1503
\\
$\widehat{\theta}^{Kdis}_3$   &1840 &2909   & 853 &2669  & 1405 & 2460 & 2262 & 2414 & 2375 & 3330 &\textbf{\emph{19}}&821&2413
\\
$\widehat{\theta}^{Kdis*}_1$    & 669&1355   & 200 &1227  & 372 & 577 & 365 & 476& 862 &1574  &380&\textbf{\emph{243}}&91
\\
$\widehat{\theta}^{Kdis*}_2$   & 217&570   & 52 & 499 &  163& 842 & 2501 & 404 & 313 &\textbf{\emph{476}}  &34&2271&1100
\\
$\widehat{\theta}^{Kdis*}_3$  &151 &494   & \textbf{18} &405  & \textbf{74} & \textbf{161} & \textbf{\emph{ 261}}&  3538 &260  &568&75&\textbf{223}&235
    \\
 \hline
 \multicolumn{14}{|c|}{$n=5000$}
  \\
  \hline
$\widehat{\theta}_1$   &101& 75  &105  &105  &541  &1132  &1331  &201&170&687  &-  & 787&584
\\
$\widehat{\theta}_2$    &\textbf{15} &\textbf{6.5514}   & \textbf{5.8965} &\textbf{22}  &538  &1203  &1601  &206&317&1139  &-  & 856 &798
\\
$\widehat{\theta}_3$   & 268& 203  &263  &210  &605  &1114  & 1068 &199&\textbf{39}&595  &-  & 776 &423
\\
$\widehat{\theta}_1^*$   & 111&125   &99  &150  &545  &1222  & 1384 &260&165&685  &606  &805  &625
\\
$\widehat{\theta}_2^*$  &\textbf{\emph{54}} &\textbf{\emph{26}}   &50  &\textbf{\emph{29}}  & 682 &1549  &2249  &381&457&914  &422  &1042  &1221
\\
$\widehat{\theta}_3^*$ & 925&697   &706  &897  &977  &1310  &975  &\textbf{102}&745& \textbf{\emph{374}} & 1455 &1085  &857
\\
  $\widehat{\theta}^{Kdis}_1$   & 476&967   &101  &881 &\textbf{148} & \textbf{135} & \textbf{61} & 418& 630 & 1151 &77&\textbf{484}&\textbf{\emph{269}}
\\
$\widehat{\theta}^{Kdis}_2$    & 2435&3873   &1157  & 3633 & 2266 & 2983 & 2087 & 3242 & 2620 &3460  &3000&935&2777
\\
$\widehat{\theta}^{Kdis}_3$   & 2044&  3748 &689  & 3427 & 1649 & 3572 & 3578 & 3254 &2202 &3771  &\textbf{11}&1015&2735
\\
$\widehat{\theta}^{Kdis*}_1$    &775 & 1891  & 143 & 1696 & 522 & 1201 & 868&779  &  950& 1814 &93&797&\textbf{94}
\\
$\widehat{\theta}^{Kdis*}_2$   &197 & 549  & \textbf{\emph{36}} & 473 & \textbf{\emph{181}}& 874 & 2522 & 423 & 299 & 447& \textbf{\emph{25}}&2874&3446
\\
$\widehat{\theta}^{Kdis*}_3$  & 153& 103  &305  & 55& 289 & \textbf{\emph{555}} & \textbf{\emph{104}} &\textbf{\emph{131}}  & \textbf{\emph{43}} & \textbf{306} &2321&\textbf{\emph{596}}&271
    \\
    \hline
 \end{tabular}
 \label{Table2-I}
\end{table} 
\begin{table}[t]
\centering
\caption{The  root mean squared error 
($k=\lfloor(\ln L)^2\rfloor)$).}
\tabcolsep=0.09cm
\begin{tabular}{|l|c c|c c|c c|c c|c c|c|c|c|c|}
  \hline
 {\scriptsize $RMSE$}         & \multicolumn{2}{c|}{MM} & \multicolumn{2}{c|}{ARMAX} & \multicolumn{2}{c|}{$ARu^+$} & \multicolumn{2}{c|}{$ARu^-$} &  \multicolumn{2}{c|}{$MA(2)$} & ARc & AR(2) & {\scriptsize GARCH}
 \\
 \cline{2-14}
 {\scriptsize $\cdot 10^{4}$/$\theta$}
&  $0.5$ & $0.8$            & 
$0.25$ & $0.75$               & $0.5$ & $0.8$                &   $0.75$ & $0.96$ & $0.5$ & $2/3$ &    $0.3$ & $0.25$ & $0.447$
\\
  \hline
  \multicolumn{14}{|c|}{$n=10^5$}
  \\
  \hline
 $\widehat{\theta}_1$     &\textbf{\emph{147}} &  167 & 126 & 176 & 396 & 383 & 1310 & 373 & 345 & 834 & 83 & 601& 453
\\
$\widehat{\theta}_2$   &158 &  159 & 131 & 173 & 468 & 463 & 1669 & 387 & 409 & 982 &358 & 720&693
\\
$\widehat{\theta}_3$    &218 &  261 & 172 & 257 & 370 & \textbf{\emph{353}} & 1039 & \textbf{\emph{367}} & 319 & 706 &245 & 519 & 406
\\
$\widehat{\theta}_1^*$   &\textbf{98} &  \textbf{118} & \textbf{\emph{85}} & \textbf{117} & 386 & 1122 & 1444 & 378 & 334 & 846 &151 & 619 &416
\\
$\widehat{\theta}_2^*$  & 151 & \textbf{\emph{134}}  & 113 & \textbf{\emph{141}} & 599 & 1534 & 2500 & 400 & 504 & 1197 &358 & 948 &1131
\\
$\widehat{\theta}_3^*$  & 343 & 478  & 293 & 421 & 412 & 980 & 515 & \textbf{359} & 336  & \textbf{439}  &67 & 444 &455
\\
$\widehat{\theta}^{Kdis}_1$   & 366&710   & 155 & 675 & \textbf{186} & 374 & 1018 & 508 & 527 & 1032 &150 &338 &\textbf{\emph{311}}
\\
$\widehat{\theta}^{Kdis}_2$   & 3148&  4925 & 1659 & 4403 & 3336 & 4512 & 3895 & 5259 & 3208 & 4254 & 3000& 1636&2962
\\
$\widehat{\theta}^{Kdis}_3$   & 3003& 4608  & 1439 & 4441 & 2598 & 4428 & 4004 & 5216 & 3255 & 4597 &\textbf{19} & 1559 &3034
\\
$\widehat{\theta}^{Kdis*}_1$   & 696& 1368 & 210 & 1255 &398  & 598 &\textbf{356}  & 1248 & 880 & 1602 &380 &\textbf{230} &\textbf{181}
\\
$\widehat{\theta}^{Kdis*}_2$  & 220&  573 &  \textbf{71}& 500 & 362 & 844 &2501  &401  & \textbf{\emph{315}} & \textbf{\emph{558}} & \textbf{\emph{34}}& 2512& 4101
\\
$\widehat{\theta}^{Kdis*}_3$  & 225 & 526  & 183 & 458 & \textbf{\emph{243}} & \textbf{263} &\textbf{\emph{462}}  & 498 & \textbf{297} & 584 &75 &\textbf{\emph{304}} & 392
\\
  \hline
\multicolumn{14}{|c|}{$n=5000$}
  \\
  \hline
  $\widehat{\theta}_1$      &493 & 683  &399  &652  &808  &1301  &1619  &\textbf{\emph{382}}&535&898  &-  &905  &988
\\
$\widehat{\theta}_2$    &484 & 662  &348  &641  &784  &1384  &1910  &407&593& 1087 &-  &953  &1132
\\
$\widehat{\theta}_3$     &796 &  1026 &624  &997  &1091  &1365  &1544  &486&814& 1037 &-  &1011  &1146
\\
$\widehat{\theta}_1^*$   &357& \textbf{480}  & 340 &\textbf{\emph{472}}  &689  &1274  &1631  &\textbf{325}&\textbf{\emph{367}}&\textbf{\emph{764}}  &559  &909  &908
\\
$\widehat{\theta}_2^*$   &\textbf{\emph{353}} & \textbf{\emph{482}}  &295  &\textbf{432}  &790  &1611  &2441  &399&554&1211  &422  &1138  &1358
\\
$\widehat{\theta}_3^*$   & 1572& 1568  &1295  &1735  &1871  & 1658 & 1754 &716&1547& 1617 & 1455 & 1618 &1979
\\
$\widehat{\theta}^{Kdis}_1$   & 486&915   & 351 & 822 &\textbf{\emph{504}} & 557 &914  &673&596&1123  &-  &\textbf{710}  &\textbf{\emph{825}}
\\
$\widehat{\theta}^{Kdis}_2$    &3214&4704   & 1543&  4618& 3238 &4581  &4483  &5477&3290& 4262 &-  &1746  &3571
\\
$\widehat{\theta}^{Kdis}_3$    &2979 &4482   &1718  & 4283 & 2765 & 4564 &4559  &5671&3218& 4425 &-  &1894  &3264
\\
$\widehat{\theta}^{Kdis*}_1$    &767 &1706   &\textbf{\emph{263}}  & 1535 & 549 & 976 & \textbf{702} &1715&943&1757  &\textbf{\emph{93}}  &\textbf{\emph{726}}  &\textbf{406}
\\
$\widehat{\theta}^{Kdis*}_2$   &\textbf{257} &  589&\textbf{194}  &516 & \textbf{397} & \textbf{874} &2518  &418&\textbf{345}&\textbf{502}  &\textbf{25}  &3299  &3817
\\
$\widehat{\theta}^{Kdis*}_3$   &812 &1002   &943  & 913 &1058  &  \textbf{\emph{921}}&\textbf{\emph{897}}  &573&720&860  &2321  &1253  &1424
\\
  \hline
\end{tabular}
 \label{Table1-K}
\end{table} 
\begin{table}[t]
\centering
\caption{The  absolute bias  
($k=\lfloor(\ln L)^2\rfloor)$.}
\tabcolsep=0.035cm
\begin{tabular}{|l|c c|c c|c c|c c|c c|c|c|c|c|}
  \hline
 {\scriptsize $|Bias|$}         & \multicolumn{2}{c|}{MM} & \multicolumn{2}{c|}{ARMAX} & \multicolumn{2}{c|}{$ARu^+$} & \multicolumn{2}{c|}{$ARu^-$} &  \multicolumn{2}{c|}{$MA(2)$} & ARc & AR(2) & {\scriptsize GARCH}
 \\
 \cline{2-14}
 {\scriptsize $\cdot 10^{4}$/$\theta$}
&  $0.5$ & $0.8$            & 
$0.25$ & $0.75$               & $0.5$ & $0.8$                &   $0.75$ & $0.96$ & $0.5$ & $2/3$ &    $0.3$ & $0.25$ & $0.447$
\\
  \hline
  \multicolumn{14}{|c|}{$n=10^5$}
  \\
  \hline
$\widehat{\theta}_1$   & 67& 37  & 45 & \textbf{\emph{39}} & 354 & 343 & 1212 & 361 & 309 & 798 & 83 & 569& 285
\\
$\widehat{\theta}_2$   & 95& 56  & 69 & 60 & 428 & 425 & 1543 & 378 & 377 & 944 &358 & 683& 522
\\
$\widehat{\theta}_3$   &\textbf{42} & \textbf{8.2991}  & \textbf{15} &\textbf{17}  & 277 & 253 & 861 & 338 & \textbf{\emph{232}} & 622 &245 & 454& 57
\\
$\widehat{\theta}_1^*$  &57 & \textbf{\emph{30}}  & \textbf{\emph{42}}  & 42 & 371 & 1117 & 1440 & 377 & 326 & 841 & 151&613 &391
\\
$\widehat{\theta}_2^*$ &126 & 78  & 95 & 96 & 589 & 1530 & 2500 & 400 & 498 & 1193 & 358& 944 &1124
\\
$\widehat{\theta}_3^*$ & \textbf{\emph{56}} & 48  & \textbf{\emph{42}} & 44 & 120 & 872 & \textbf{\emph{262}} & \textbf{235} & \textbf{82} & \textbf{242} & 67 & 316 &175
\\
$\widehat{\theta}^{Kdis}_1$   &331 &  653 & 122 & 615 & \textbf{69} & \textbf{208} & 486 & 321 & 492 & 966 &150 &249 &\textbf{13}
\\
$\widehat{\theta}^{Kdis}_2$   & 2084& 3275  & 1128 & 2821 &2169  & 2006 & \textbf{22} & 2675 & 2200 & 2978 & 3000& 521&1366
\\
$\widehat{\theta}^{Kdis}_3$   & 2022& 3119  & 891 & 3047 & 1478 & 2458 & 1895 & 3226 &2359  & 3575 &\textbf{19} & 802&2215
\\
$\widehat{\theta}^{Kdis*}_1$   &693 &  1362 & 203 &1249  & 388 & 584 & 337 & 1239 & 878 & 1598 &380 &\textbf{191} &\textbf{\emph{49}}
\\
$\widehat{\theta}^{Kdis*}_2$  & 217& 571  & 56 & 499 & 149 & 843 & 2501 & 401 & 314 & \textbf{\emph{482}} & \textbf{\emph{34}}&1676 &309
\\
$\widehat{\theta}^{Kdis*}_3$ & 159&  483 & \textbf{\emph{17}} & 417 & \textbf{\emph{87}} &\textbf{\emph{211}}  & 265 & 484 & 258 & 556 &75 & \textbf{\emph{234}}&225
\\
  \hline
\multicolumn{14}{|c|}{$n=5000$}
  \\
  \hline
  $\widehat{\theta}_1$      &91 & 147  & 112 & 128 &553  & 1173 & 1470 &268&144&629  &-  &781  &632
\\
$\widehat{\theta}_2$    & \textbf{17}& \textbf{\emph{41}}  &\textbf{28}  & \textbf{12} &543  & 1255 &1767  &331&299&890  &-  & 827 &812
\\
$\widehat{\theta}_3$     & 255&  290 & 235 & 326 & 619 & 1106 & 1212 &210&\textbf{\emph{79}}&\textbf{\emph{260}}  &-  &789  &536
\\
$\widehat{\theta}_1^*$   &105 &  121 & 103 & 146 & 529 &  1274&1580  &294&175&674  &559  &816  &630
\\
$\widehat{\theta}_2^*$   & \textbf{\emph{67}}&  \textbf{20} &55  &\textbf{\emph{26}}  &662  &1560  &2435  &398&445&1148  &422  &1065  &1219
\\
$\widehat{\theta}_3^*$   & 814&   722&761  &840  &959  &1277  &1088  &\textbf{89}&762&536  &1455  &1111  &830
\\
$\widehat{\theta}^{Kdis}_1$   & 322&731   &47  &610  & \textbf{26} & \textbf{133} &\textbf{\emph{421}} &425&491&947  & -& 479 &253
\\
$\widehat{\theta}^{Kdis}_2$    &2140&  3018 &944 &3037  &2001  &2049  & 543 &2890&2277&2926  &-  &\textbf{359}  &2115
\\
$\widehat{\theta}^{Kdis}_3$    &1871& 2941  &1072  &2757  & 1477 & 2572& 2638 &3702&2215&3232  & - &952  &2264
\\
$\widehat{\theta}^{Kdis*}_1$    &747 & 1680  &139  &1511  &477  &943  &645 &1699&933&1748 & \textbf{\emph{93}} & 648 &\textbf{118}
\\
$\widehat{\theta}^{Kdis*}_2$   &203 &55   &\textbf{\emph{39}}  &478  &\textbf{\emph{158}}  &863  &2518  &418&305& 454 & \textbf{25} & 1489 &2966
\\
$\widehat{\theta}^{Kdis*}_3$   &195& 173  &352  &91  &239  &\textbf{\emph{559}}  &\textbf{106}    &\textbf{\emph{136}}& \textbf{77} & \textbf{203}&2321 &\textbf{\emph{378}} &\textbf{\emph{237}}
\\
  \hline
\end{tabular}
 \label{Table2-K}
\end{table} 
 We take the intervals estimator as $\widehat{\theta}_0$ since it requires only $u$ as parameter. 
\begin{remark}
For the $K$-gaps estimator the algorithm is the same, but instead of $\{Y_i\}$ one has to use the normalized $K$-gaps $\{\overline{F}(u)S(u)_{i}^{(K)}\}$. For each value of $u$ one can examine different values of $K$, for instance, $K\in\{1,2,...,20\}$ may be taken. For $K=0$ $\{Y_i\}$ are still used.
\end{remark}
\begin{remark}
\label{Rem5} As the solutions of (\ref{4}) 
 may not exist among considered quantiles for  given $k$ and $K$, we propose to use the inequality
\begin{eqnarray}\label{5}
\widetilde{\omega}^2_L(\widehat{\theta})\le \delta_2 
\end{eqnarray}
as an alternative, where $\delta_2=1.49$ is
 the $99.98\%$ quantile  of the C-M-S statistic.
\end{remark}
\begin{remark}
The discrepancy method is somewhat similar to the multilevel approach by Sun and  Samorodnitsky (2018), where a fixed number of levels $u_n^1<...<u_n^m$ is selected such that $\overline{F}(u_n^s)/\overline{F}(u_n^m)\to\tau_s/\tau_m$ for some $\tau_1>...>\tau_m>0$. In our case, the number of thresholds, which are the solutions of the discrepancy equation, is random and
thus, it cannot be considered as 
an additional parameter.
\end{remark}
The discrepancy method is universal and any estimator depending on $u$ can substitute $\widehat{\theta}$  in (\ref{4}). In case of the free-threshold estimators one can express a cluster identification parameter such as the block size as depending on $u$ and find the latter by the discrepancy method. For example, the block size can be selected as $b(u)=\lfloor n/L(u)\rfloor$. The simulation study of this case is out of scope of our paper.
Comparison of the threshold-based intervals and $K$-gaps estimators with other estimators based on other tuning parameters (like the block size or the length for runs of non-exceedances)
is very complicated since the numbers $L(u)$ and $N_C=N_C(u)$ used for calculations are random.
That is the reason our comparison concerns only intervals and $K$-gaps estimators coupled with different threshold choice methods.
\subsection{Models} In our simulation study we consider the processes MM, ARMAX, 
 AR(1), AR(2), MA(2) and GARCH(1,1) 
 all with known values $\theta$.
 The simulation is repeated $1000$ 
 times with the sample size $n= 10^5$ 
 of initial measurements $\{X_1,...,X_n\}$. 
 Big  sample sizes may lead, however, 
 to moderate size samples $L(u)$ of normalized inter-exceedance 
 times $\{Y_1,...,Y_{L(u)}\}$. We recall the definitions of the processes.
The $m$th order MM process is 
$X_t=\max_{i=0,...,m}\{\alpha_i Z_{t-i}\}$, $t\in \textmd{Z}$,
where $\{\alpha_i\}$ are constants with $\alpha_i \ge 0$, $\sum_{i=0}^m \alpha_i=1$, and
 $Z_t$ are i.i.d. standard Fr\'{e}chet distributed r.v.s with the cdf
$F(x)=\exp\left(-1/x\right)$, for $x>0$.
 The extremal index of the process is equal to $\theta=\max_i\{\alpha_i\}$  (Ancona-Navarrete and Tawn 2000). The distribution of $\{X_t\}_{t\ge 1}$ is standard
Fr\'{e}chet. Values $m=3$ and $\theta\in\{0.5,0.8\}$ corresponding to $\alpha\in\{0.5,0.3,0.15,0.05\}$ and $\alpha\in\{0.8,0.1,0.008,0.02\}$, respectively, are taken for our study.
\\
The ARMAX process is determined as $X_t=\max\{\alpha X_{t-1},(1-\alpha)Z_t\}$, $t\in \textmd{Z},$
where $0\le \alpha < 1$,
 $\{Z_t\}$ are i.i.d standard Fr\'{e}chet distributed r.v.s and $P\{X_t\le x\}=\exp\left(-1/x\right)$ holds assuming $X_0=Z_0$. The extremal index of the process is given by $\theta=1-\alpha$,  Beirlant et al. (2004). $P\{X_t\le x\}=\exp(-1/x)$ holds assuming $X_0=Z_0$. We consider $\theta\in\{0.25,0.75\}$.
\\
The positively correlated AR(1) process with uniform noise ($ARu^+$) is defined by $X_j=(1/r) X_{j-1}+\epsilon_j$, $j \geq 1$ and $X_0\sim  U(0,1)$
with $X_0$ independent of  $\epsilon_j$. $X_1\sim  U(0,1)$ holds. For a fixed integer $r\geq 2$ let $\epsilon_n$, $n\geq 1$ be i.i.d. r.v.s with
$P\{\epsilon_1=k/r\}=1/r$, $k\in\{0,1, \ldots , r-1\}$. The extremal index of $ARu^+$ is $\theta=1-1/r$ (Chernick et al. 1991). $\theta\in\{0.5,0.8\}$ corresponding to $r\in\{2,5\}$ are taken. The negatively correlated AR(1) process with uniform noise ($ARu^-$) is defined by $X_j=-(1/r) X_{j-1}+\epsilon_j$ with the similar distributed $\epsilon_n$ but with $k\in\{1, \ldots , r\}$. Its extremal index is $\theta=1-1/r^2$ (Chernick et al. 1991). The same $r$ were taken corresponding to $\theta\in\{0.75,0.96\}$.
 \\
We simulate  the MA(2) process (Sun and  Samorodnitsky 2018) $X_i= p Z_{i-2} + q Z_{i-1} + Z_i$, $i\ge 1,$
 with $p>0$, $q<1$, and i.i.d. Pareto random variables $Z_{-1}, Z_0, Z_1,...$ with $P\{Z_0>x\}=1$ if $x<1$, and
 $P\{Z_0>x\}=x^{-\alpha}$ if $x\ge 1$.
for some $\alpha>0$. The extremal index of the process is $\theta=\left(1+p^{\alpha}+q^{\alpha}\right)^{-1}$. The cases $\alpha=2$, $(p,q)=(1/\sqrt{2}, 1/\sqrt{2}), (1/\sqrt{3}, 1/\sqrt{6})$ with corresponding $\theta\in\{1/2, 2/3\}$ are considered. The distribution of the sum of weighted i.i.d. Pareto r.v.s behaves like a Pareto distribution in the tail and its exact form may be obtained by Ramsay (2008). 
\\
We consider also 
processes studied in (Ferreira 2018b; Northrop 2015; S\"{u}veges and Davison 2010). These comprise the
AR(1) process $X_j=0.7 X_{j-1}+\epsilon_j$, where $\epsilon_j$ is standard Cauchy distributed and $\theta=0.3$ (ARc);
the $AR(2)$ process $X_j=0.95 X_{j-1}-0.89X_{j-2} +\epsilon_j$, where $\epsilon_j$  is Pareto distributed with tail index $2$ and $\theta=0.25$;  $GARCH(1,1)$, $X_j=\sigma_j\epsilon_j$, with $\sigma_j^2=\alpha+\lambda X_{j-1}^2+\beta \sigma_{j-1}^2$, $\alpha=10^{-6}$, $\beta=0.7$, $\lambda=0.25$, with $\{\epsilon_j\}_{j\ge 1}$ an i.i.d. sequence of standard Gaussian r.v.s and $\theta=0.447$ (see Laurini and Tawn 2012).

\subsection{Notations}
Tables \ref{Table1-KgapsK0Theta} and \ref{Table2-KgapsK0Theta} contain the statistics (\ref{19}) for the $K-$gaps estimates with $K=0$ coupling with the discrepancy method (\ref{4}) denoted as $\widehat{\theta}^{K0dis}_i$, $i\in\{1,2,3\}$. $\lfloor \theta L\rfloor$ and $\lfloor \widehat{\theta}_0 L\rfloor$ are considered as options for $k$, where $\widehat{\theta}_0$ is a pilot intervals estimate. The sign $'-'$ means that there are no solutions of the discrepancy equation.
\\
The rest of the tables is a partition regarding $k$ for the intervals and $K$-gaps estimators coupled with the discrepancy method. We study
$k=\lfloor \widehat{\theta}_0L\rfloor$ 
in Tables \ref{Table1-N} and \ref{Table2-N}, $k=\lfloor\min(\widehat{\theta}_0L, \sqrt{L})\rfloor$  
in Tables \ref{Table1-I} and \ref{Table2-I}, and $k=\lfloor(\ln L)^2\rfloor)$  in Tables \ref{Table1-K} and \ref{Table2-K}. In Tables \ref{Table1-N} and \ref{Table2-N}
the statistics (\ref{19}) corresponding to the intervals estimates coupled with the discrepancy method (\ref{4}) are denoted as $\{\hat{\theta}_i\}$, $i\in\{1,2,3\}$. 
The 
$K-$gaps estimates with pairs $(u,K)$ selected by  (\ref{4}) are denoted as $\widehat{\theta}^{Kdis}_i$, $i\in\{1,2,3\}$, and 
with IMT-selected pairs $(u,K)$ 
as $\widehat{\theta}^{Kimt}$.
Statistics (\ref{19}) relating to the intervals and $K-$gaps estimators and corresponding to solutions of the discrepancy inequality 
(\ref{5}) 
are denotated by asterisks in all tables. The intervals estimate with the threshold $u$ selected by the "plateau-finding" Algorithm 1  by Ferreira (2018a) is denoted as $\widehat{\theta}^{IA1}$. This algorithm seems to be the best one for the intervals estimator among other algorithms proposed in Ferreira (2018a) according to the provided simulation study. For this algorithm we consider the bandwidth $d=[w n]$ with $w=0.25$ and compute the moving average of $2d+1$ successive  points of $\widehat{\theta}$. The value $w=0.005$ used in Ferreira (2018a) demonstrates slightly worse accuracy uniformly for all processes and we do not show it in  Tables \ref{Table1-N} and \ref{Table2-N}.
\\
The values in bold and italic bold correspond to the
first and second best performances.
\subsection{Conclusions and practical recommendations}
We propose to select a threshold of the threshold-based intervals and $K-$gaps estimators and a tuning parameter  of free-threshold procedures as solutions of the $\omega^2$ discrepancy equation, where the discrepancy value is equal to the mode of the $\omega^2$- statistic, i.e. to its most likelihood value.
\\
On the first view, the intervals threshold-based estimator does not require another parameter to be specified apart of the threshold. The intervals estimator coupled with the discrepancy method works the same way as the $K-$gaps estimator. An additional regularization parameter such as the moving window size for the "plateau-finding" algorithm A1 (Ferreira 2018a) or
the number of the largest order statistics $k$ 
is  required to choose the threshold anyway. It is shown in our paper that there is a potential benefit in choosing
$k$ 
jointly with the threshold.
\\
It is proposed in Ferro and Segers (2003) to select the largest $C-1=\lfloor\theta L(u)\rfloor$ interexceedance times which are approximately independent intercluster times as associated with the threshold $u$.
We follow a similar way, i.e. $k=\lfloor \widehat{\theta}_0L\rfloor$ is used as one of the choices of $k$.
\\
The best 'ideal estimator' $\widehat{\theta}^{K_0}$ coupled with discrepancy method (\ref{4}), where $k$ is taken equal to $\lfloor\theta L\rfloor$ and $\lfloor\hat{\theta}_0 L\rfloor$, is presented in Tables \ref{Table1-KgapsK0Theta} and \ref{Table2-KgapsK0Theta}.
 Fig. \ref{fig:0} shows that $\widehat{\theta}^{K_0}$ with $k=\lfloor\theta L\rfloor$ outperforms other estimators in Tables \ref{Table1-N}-\ref{Table2-K} with the best RMSE. The results degrade if one selects the intervals estimate as $\hat{\theta}_0$. A small deviation from $\theta$ does not worsen the best estimate much.
 \\
 The discrepancy method is competitive with threshold  choices such as the IMT and "plateau-finding" algorithms and it improves substantially the existing intervals and $K-$gaps estimates coupled with the mentioned 
adjustment methods. 
Fig. \ref{fig:1} corresponding to Tables \ref{Table1-N} and \ref{Table2-N} shows that the K-gaps estimator works better, if $u$ is selected by the discrepancy method than by the IMT method. According to our simulation study the K-gaps estimator coupled with the IMT method demonstrates a slow convergence as the sample size increases. The IMT method requires more computation time due to a  full search among pairs $(u,K)$. Generally,  the K-gaps estimator  works better than the intervals estimator
  both coupled with the discrepancy method. The intervals estimator coupled with algorithm A1 provides the RMSE similar to the discrepancy method coupled with both intervals and $K$-gaps estimates only for MM and ARMAX processes, see Fig. \ref{fig:1}.
\\
The  discrepancy inequalities can be applied when the solutions of the discrepancy equation do not exist among the considered quantiles for given $k$ and $K$. This may slightly improve the RMSE and the absolute bias of both intervals and K-gaps estimates in comparison with the usage of the discrepancy equations, see Fig. \ref{fig:2}. This property is due to a larger number of solutions.
\\
Fig. \ref{fig:3} aims to compare the impact of the choice of $k$.
It shows that $k=\lfloor\min(\widehat{\theta}_0 L, \sqrt{L})\rfloor$
and $k=\lfloor(\ln L)^2\rfloor$ (both satisfy Theorem \ref{TR2}) provide similar
values of the best
RMSE. $k=\lfloor \widehat{\theta}_0 L\rfloor$ provides the best accuracy. 
\\
Fig. \ref{fig:4} aims to find the best measure from (\ref{19}). Ratios $\{RMSE(\theta_j)/\min_{i\in\{1,2,3\}}RMSE(\theta_i)\}$, $j\in\{1,2,3\}$ are compared. One may conclude that $\widehat{\theta}_1$ provides consistently better accuracy than $\widehat{\theta}_2$ and $\widehat{\theta}_3$.
 \\
 By the simulation study we recommend the $K$-gaps estimator 
 coupled with the discrepancy method (\ref{5}) with $k=\lfloor\hat{\theta}_0L\rfloor$ and an accurate pilot estimate $\hat{\theta}_0$, and the measure $\widehat{\theta}_1$. 
 \\
The impact of the heaviness of tail on the accuracy of the discrepancy method remains an open problem. Intuitively, the heaviness of the distribution tail may impact on the rate of  convergence of the  exceedance point process to a compound Poisson process and hence,  on the convergence of the  distribution of the discrepancy statistic to the limit distribution of the C-M-S statistic.

\section{Application to real data}\label{Sec5}
\subsection{First example}
Following Ferreira (2018a) we consider two data sets of the daily maximum temperatures (in $0.1$ degrees Celsius) of July at Uccle (Belgium), from $1833$ to $1999$ and from $1900$ to $1999$ with sample sizes $n\in\{5177, 3100\}$, respectively. The data are available at $"http://lstat.kuleuven.be/Wiley/Data/ecad00045TX.txt"$.  The extremal index of the  smaller sample was shown to be  ranged between $0.49$ and $0.56$ in Beirlant et al. (2004); a reduced-bias version of Nandagopalan's runs estimator applied in Ferreira (2018a) has shown  $0.41$ and $0.57$; and the wide range of estimators has shown $0.10$ and $0.57$ in Ferreira (2018a).  We have analyzed the intervals
and $K-$gaps estimators coupled with the discrepancy method based on the algorithm in Section \ref{Sec4}. The $K-$gaps estimator  with the IMT method and the intervals with "plateau-finding" Algorithm 1 with  $\omega=0.3$ were also applied here and in the next example.
\begin{table}[t]
\centering
\caption{Extremal index estimates for Uccle data.}
\tabcolsep=0.035cm
\begin{tabular}{|l|c|c|c|c| c |c| c| c |c|}
  \hline
  {\scriptsize $n$}   &'Kimt' & 'IntA1' &   &  \multicolumn{2}{c|}{'Intdis'} &  \multicolumn{2}{c|}{'Kdis'} &  \multicolumn{2}{c|}{'K0dis'}
  \\
  \cline{1-3}\cline{5-10} 
  & & &  &$s=\widehat{\theta}_0$ & $s=0.51$
     &$s=\widehat{\theta}_0$ & $s=0.51$
     &$s=\widehat{\theta}_0$ & $s=0.51$
  \\ \cline{5-10}
  $3100$ & 0.5133 & 0.4625 &  $\widehat{\theta}_1^*$ & 0.5329 & 0.5741& 0.5670 &0.5879 & 0.5383 & 0.5148 
  \\
  & && $\widehat{\theta}_2^*$ &0.4199&0.4637& 0.5232 &0.5232&0.4186&0.5104
  \\
  & && $\widehat{\theta}_3^*$ &0.9575&0.9575& 0.7244 & 0.7244 &1& 0.5520
    \\
     \cline{5-10}
 & & & &$s=\widehat{\theta}_0$ & $s=0.56$  &$\widehat{\theta}_0$ & $s=0.56$
  &$s=\widehat{\theta}_0$ & $s=0.56$  
    \\
  \cline{5-10}
   $5177$ & 0.5695 & 0.4392 & $\widehat{\theta}_1^*$ & 0.4655 & 0.4837& 0.5632 &0.6251 & 0.4741 & 0.5691
   \\
    & && $\widehat{\theta}_2^*$ &0.4184 &0.4919 & 0.5285 & 0.5662 &0.4201& 0.5604
  \\
   & && $\widehat{\theta}_3^*$ &0.5618 & 0.5618 & 0.7024 & 0.7024& 0.6524 &0.6524
    \\
    \hline
\end{tabular}
\label{Table9}
\end{table}

 'Kdis', 'K0dis'
 and 'Intdis' are calculated with $k=\lfloor sL\rfloor$, where $s$ was taken equal to the pilot intervals estimate $\widehat{\theta}_0$ for each threshold value $u$ or to values $\{0.51,0.56\}$ for $n\in\{3100, 5177\}$, respectively,  based on previous estimation of $\theta$ and the 'Kimt' estimates. 
The discrepancy inequality method (\ref{5}) was used.
One may trust more $\widehat{\theta}^*_1$ and $\widehat{\theta}^*_2$ as well as 'Kdis', 'K0dis' estimates since they provide better results on the simulation. The results are shown
in Table \ref{Table9}.
\subsection{Second example}
We use the data corresponding to Figure S18 in Raymond et al. (2020) and kindly provided by the authors, which represent daily-maximum dewpoint temperatures at  station Dhahran, Saudi Arabia. This station is among several selected stations where a wet-bulb temperature (TW) has exceeded $TW=33^oC$ at least $5$ times.
The dates span from 1 Jan 1979 to 31 Dec 2017. The sample size is equal to $n=13866$ due to missing observations. The estimated values of $\theta$ are shown in Table \ref{Table10}.
\begin{table}[t]
\centering
\caption{Extremal index estimates for dewpoint temperatures data.}
\tabcolsep=0.035cm
\begin{tabular}{|l|c|c|c|c|c| c |}
  \hline
  {\scriptsize $n$}   &'Kimt' & 'IntA1' &   &  'Intdis' &  'Kdis' &  'K0dis' 
  \\
  \cline{1-3}\cline{5-7}
  & & &  &\multicolumn{3}{c|}{$k=\lfloor\widehat{\theta}_0L\rfloor$}
  \\\cline{5-7}
  $13866$ & 0.4753 & 0.1541 &  $\widehat{\theta}_1^*$ & 0.2489 & 0.3178& 0.2749 
  \\
  & &
  & $\widehat{\theta}_2^*$ &0.2003&0.2016& 0.2021 
  \\
  & && $\widehat{\theta}_3^*$ &0.4092&0.4765&0.5181 
      \\
    \hline
\end{tabular}
\label{Table10}
\end{table}
 
\section{Proofs}\label{Proofs}
\subsection{Proof of Theorem \ref{TR1}}\label{S-5.1}

Consider the conditional distribution of 
$\omega^2_k(\theta)$
given $U_{L-k, L} = t_k.$ 
According to Lemma \ref{Lem1} and the condition $\limsup_{n\to\infty}k/L<\theta$
the conditional joint distribution of the set of the order statistics
$\{U_{i, L}\}_{i=L-k+1}^L$ asymptotically equals to the joint distribution of the set
of order statistics $\{U^\ast_{i, k}\}_{i=1}^k$ of a sample
$\{U_i^\ast\}_{i=1}^k$ from the uniform distribution on $[t_k,1].$
Therefore, it holds 
\begin{eqnarray*}&&\omega^2_k(\theta) \stackrel{d}{=} \frac{1}{(1 - t_k)^2}
\left(\sum\limits_{i=1}^k\left(U^\ast_{i, k} -
t_k-\frac{i-0.5}{k}(1-t_k)\right)^2\right) + \frac{1}{12k}.\end{eqnarray*}
Moreover,
$V^\ast_{i,k} = U^\ast_{i,k}-t_k$ are the order statistics of a sample
$\{V^\ast_i\}$ from the uniform distribution on $[0, 1-t_k]$. Hence, it follows
 \[\omega^2_k(\theta) \stackrel{d}{=}
\frac{1}{(1 - t_k)^2} \left(\sum\limits_{i=1}^k\left(V^\ast_{i, k} -
\frac{i-0.5}{k}(1-t_k)\right)^2\right) + \frac{1}{12k}.\] Finally,
$W^\ast_{i,k} = V^\ast_{i,k}/(1-t_k)$ are the order statistics of a
sample
 $\{W^\ast_i\}$ from the uniform distribution on $[0, 1]$. Therefore, we get
 \begin{eqnarray*}\omega^2_k(\theta) &\stackrel{d}{=}&
\sum\limits_{i=1}^k\left(W^\ast_{i,
k} - \frac{i-0.5}{k}\right)^2 + \frac{1}{12k}.\end{eqnarray*}  It is easy to see,
that the last expression is the C-M-S 
statistic and it
converges in distribution to the r.v. $\xi$ with the
cdf 
$A_1$ independently of the value of $t_k$.

\subsection{Proof of Theorem \ref{TR2}}\label{S-5.2}
Let $\{E_i^{(L)}\}_{i=1}^{\lfloor\theta L\rfloor}$ be the sequence of r.v.s satisfying condition (\ref{newcondition}).
Let us denote $t_i = 1 - \theta e^{-\theta E^{(L)}_{l-i, l}},
~~\widehat{t}_i = 1 - \widehat{\theta}_{n}
e^{-\widehat{\theta}_{n} E^{(L)}_{l-i, l}},$ 
\begin{eqnarray*}
&& a_{i} = ((L-i) -(L-k) - 0.5)/k, 
~~0\leq i< k.\end{eqnarray*}
Turning back to (\ref{omegasquared2}) and (\ref{omegasquared}), we consider the following 
difference
\begin{eqnarray}\nonumber&&\widetilde{\omega}^2_L(\widehat{\theta}_{n}) -
\omega^2_k(\theta) =
\left(\widetilde{\omega}^2_L(\widehat{\theta}_{n}) - \omega^2_k(\widehat{\theta}_{n})\right) + \left(\omega^2_k(\widehat{\theta}_{n}) - \omega^2_k(\theta)\right) \\
\label{9} &=& \left(\widetilde{\omega}^2_L(\widehat{\theta}_{n}) - \omega^2_k(\widehat{\theta}_{n})\right) + \left(\omega^2_k(\widehat{\theta}_{n}) -
\frac{1}{12k}\right)\left(1 - \frac{(1 - \widehat{t}_k)^2}{(1 -
t_k)^2}\right)\nonumber
\\
&+&
\left(\left(\omega^2_k(\widehat{\theta}_{n}) -
\frac{1}{12k}\right)\frac{(1 - \widehat{t}_k)^2}{(1 - t_k)^2} -
\left(\omega^2_k(\theta) -
\frac{1}{12k}\right)\right).\end{eqnarray}
The third term on the
right-hand side of (\ref{9}) is equal to 
\begin{eqnarray*}&&
\frac{1}{(1-t_k)^2}\left(\sum_{i=0}^{k-1}\left(\widehat{t}_i -
 \widehat{t}_k - a_i(1-\widehat{t}_k)\right)^2 - \sum_{i=0}^{k-1}\left(t_i -
 t_k - a_i(1-t_k)\right)^2\right).\end{eqnarray*}
  Using the 
  relation $x^2 -
y^2 = (y-x)^2 - 2y(y-x),$ we obtain
\begin{eqnarray*}&&\frac{1}{(1-t_k)^2}\left(\sum_{i=0}^{k-1}\left(\widehat{t}_i -
 \widehat{t}_k - a_i(1-\widehat{t}_k)\right)^2 - \sum_{i=0}^{k-1}\left(t_i -
 t_k - a_i(1-t_k)\right)^2\right)
 \\
 &=&
\frac{1}{(1-t_k)^2}\sum_{i=0}^{k-1}\left(t_i - \widehat{t}_i - (t_k - \widehat{t}_k)\left(1 - a_{i}\right)\right)^2
\\
&-&\frac{2}{(1-t_k)^2}\sum_{i=0}^{k-1} \left(t_i - \widehat{t}_i
- (t_k - \widehat{t}_k)\left(1 - a_{i}\right)\right)\left(t_i -
 t_k - a_i(1-t_k)\right).\end{eqnarray*}
Thereby, we can rewrite
\begin{eqnarray}&&\omega^2_k(\widehat{\theta}_{n}) -
\omega^2_k(\theta) =
\left(\omega^2_k(\widehat{\theta}_{n}) -
\frac{1}{12k}\right)\left(1 - \frac{(1 - \widehat{t}_k)^2}{(1 -
t_k)^2}\right)\nonumber
\\
& +&
\frac{1}{(1-t_k)^2}\sum_{i=0}^{k-1}\left(t_i - \widehat{t}_i - (t_k - \widehat{t}_k)\left(1 - a_{i}\right)\right)^2\nonumber
\\
&-&
\frac{2}{(1-t_k)^2}\sum_{i=0}^{k-1} \left(t_i - \widehat{t}_i
- (t_k - \widehat{t}_k)\left(1 - a_{i}\right)\right)\left(t_i -
 t_k - a_i(1-t_k)\right).\label{main}\end{eqnarray}
Let us find the asymptotics of the difference $t_i - \widehat{t}_i = \theta
e^{-\theta E^{(L)}_{l-i, l}} - \widehat{\theta}_{n}
e^{-\widehat{\theta}_{n} E^{(L)}_{l-i, l}},$ $0\leq i\leq k.$ We have
\begin{eqnarray}\label{hatminus}\widehat{\theta}_{n} e^{-\widehat{\theta}_{n} E^{(L)}_{l-i, l}} -
\theta e^{-\theta E^{(L)}_{l-i, l}} &=& (\widehat{\theta}_{n} - \theta)
e^{-\widehat{\theta}_{n} E^{(L)}_{l-i, l}} \nonumber
\\
&+& \theta e^{-\theta E^{(L)}_{l-i, l}}(e^{-E^{(L)}_{l-i, l}(\widehat{\theta}_{n} - \theta)} -
1).\end{eqnarray}
By (\ref{2}) 
it holds
\begin{eqnarray}\label{10a}&&\widehat{\theta}_{n} - \theta =
O_P\left(\frac{1}{\sqrt{m_n}}\right).\end{eqnarray}
We have $U_{i+1, l}\stackrel{d}{=}
e^{-E^{(L)}_{l-i, l}\theta}\leq e^{-E^{(L)}_{l-k, l}\theta}
\stackrel{d}{=} U_{k+1, l},$ $0\leq i< k<l$, where $\{U_{i,l}\}$ are the
order statistics arising from a standard uniform distribution. 
By Lemma \ref{Lem2} we get
\[\frac{U_{k, l} -\frac{k-1}{l-1}}{\sqrt{\frac{(k-1)(l-k)}{(l-1)^3}}}
\stackrel{d}{\to} N(0,1)\]
as $k\to\infty,$ $L\to\infty,$
$L-k\to\infty.$ Since both
\begin{eqnarray*}\frac{k-1}{l-1} - \frac{k}{l} = -\frac{l-k}{l(l-1)}&&~ \mbox{ and }~
\frac{l-k}{l(l-1)} \sqrt{\frac{(l-1)^3}{(l-k)(k-1)}} =
\sqrt{\frac{(l-1)(l-k)}{(k-1)l^2}} \end{eqnarray*}
are equivalent to $o(1)$,
then under the same
conditions as in Lemma \ref{Lem2} and using Slutsky's theorem, we obtain
\begin{eqnarray*}\label{Smirnov}&&\frac{U_{k, l} -
\frac{k}{l}}{\sqrt{\frac{k(l-k)}{l^3}}} \stackrel{d}{\to}
N(0,1).\end{eqnarray*}
This result implies
\begin{equation}1-t_k = \theta
e^{-E^{(L)}_{l-k, l}\theta} =
\frac{k}{L}(1+o_P(1)),\label{t}\end{equation}
and $E^{(L)}_{l-i, l} = \ln (l/i)/\theta(1+o_P(1))$  as $i\leq k,$
$i\to\infty$ holds.
Using the condition $(\ln L)^2 = o(m_n),$ we have $E^{(L)}_{l-i, l}(\widehat{\theta}_{n} - \theta) = o_P(1)$ and
\begin{eqnarray}\label{3}&&e^{-E^{(L)}_{l-i, l}(\widehat{\theta}_{n} - \theta)} - 1 = -E^{(L)}_{l-i, l}(\widehat{\theta}_{n} - \theta)(1+o_P(1)).\end{eqnarray}
Then the
first term on the right-hand side in (\ref{hatminus}) is
asymptotically smaller than the second one. 
Hence, from (\ref{hatminus}), (\ref{t}) and (\ref{3}) we obtain
\begin{eqnarray*}&&\hat{t}_i - t_i  =
-\frac{i}{L}E^{(L)}_{l-i, l}(\widehat{\theta}_{n} - \theta)(1+o(1))=O_P\left(\frac{i\ln(i/L)}{L\sqrt{m_n}}\right).\end{eqnarray*}
Therefore, the asymptotics of the expression $t_i - \widehat{t}_i -
(t_k - \widehat{t}_k)\left(1 - a_i\right)$ is the following
\begin{eqnarray}
&& \frac{i}{L}E^{(L)}_{l-i, l}(\widehat{\theta}_{n} - \theta)(1+o_P(1)) -
\frac{i+0.5}{k}\frac{k}{L}E^{(L)}_{l-k, l}(\widehat{\theta}_{n} -
\theta)(1+o_P(1)) \nonumber
\\
&=&
 \label{substr}
-\frac{i}{L}(E^{(L)}_{l-k, l} - E^{(L)}_{l-i, l})(\widehat{\theta}_{n} -
\theta)(1+o_P(1)) = O_P\left(\frac{i\ln(k/i)}{L \sqrt{m_n}}\right),
\end{eqnarray}
since (\ref{10a}) and
  \begin{equation}E^{(L)}_{l-i, l} - E^{(L)}_{l-k, l} \stackrel{d}{=}
\frac{\ln(k/i)}{\theta}(1+o_P(1))\label{renyi}\end{equation}
hold. 
Note that the maximum of the function $f(x) =
\ln (a/x) x,$ where $a$ is a positive constant, is achieved in the
point $x_0 = a/e,$ so $i \ln (k/i) \leq k/e$, 
$1\leq i\leq
k.$ Hence, from (\ref{t}) and (\ref{substr}) the asymptotics of the second
 term on the right-hand side of (\ref{main}) is given by
\begin{eqnarray*}&&\frac{1}{(1-t_k)^2}\sum_{i=0}^{k-1} \left(t_i - \widehat{t}_i -
(t_k - \widehat{t}_k)\left(1 - a_i\right)\right)^2
= O_P\left(\frac{k}{m_n}\right) = o_P(1)\end{eqnarray*}
due to
(\ref{asymp}). Now we estimate the asymptotics of the third summand on
the right-hand side of (\ref{main}). An appeal to (\ref{substr}) and the
Cauchy-Schwarz inequality gives us the following
\begin{eqnarray*}&&\frac{2}{(1-t_k)^2}\left|\sum_{i=0}^{k-1} \left(t_i - \widehat{t}_i -
(t_k - \widehat{t}_k)\left(1 - a_i\right)\right)\left(t_i -
 t_k - a_i(1-t_k)\right)\right|\nonumber
 \\
  &\leq &
\max_{i} \left(\frac{i}{L}(E^{(L)}_{l-i, l} - E^{(L)}_{l-k, l})|\widehat{\theta}_{n} -
\theta|\right)\cdot
\frac{2}{(1-t_k)^2} \nonumber
\\
&\cdot &\sum_{i=0}^{k-1}
\left|1 - \theta\exp(-E^{(L)}_{l-i, l}\theta) -
 t_k - \frac{k - i - 0.5}{k}(1-t_k)\right|(1+o_P(1))\nonumber
 \\
 &\leq &
O_P\left(\frac{k}{L\sqrt{m_n}}\right) \frac{2\sqrt{k}}{(1-t_k)^2}\cdot \nonumber
\\
&\cdot &\sqrt{\sum_{i=0}^{k-1} \left(1 - \theta\exp(-E^{(L)}_{l-i, l}\theta) -
 t_k - \frac{k-i -0.5}{k}(1-t_k)\right)^2}\nonumber
 \\
 & = &
O_P\left(\frac{k^{3/2}}{L\sqrt{m_n}}\right)\frac{2}{(1-t_k)}
\sqrt{\omega^2_k(\theta)-\frac{1}{12k}} =
O_P\left(\frac{k^{3/2}}{L\sqrt{m_n}}\cdot\frac{L}{k}\right) \nonumber
\\
&=&
O_P\left(\frac{\sqrt{k}}{\sqrt{m_n}}\right) =
o_P(1),\label{long}\end{eqnarray*} where
the last two strings follow from (\ref{asymp}) and
(\ref{t}) and since $\omega^2_k(\theta) - \frac{1}{12k} =
O_P(1)$ holds from Theorem \ref{TR1}. 
\\
Thus, 
by (\ref{9}) 
the sum of the second and the third
terms in (\ref{main}) is equal to
\begin{equation}\left(\omega^2_k(\widehat{\theta}_{n}) -
\frac{1}{12k}\right)\frac{(1 - \widehat{t}_k)^2}{(1 - t_k)^2} -
\left(\omega^2_k(\theta) -
\frac{1}{12k}\right)=o_P(1).\label{23}\end{equation}
Now we derive that the asymptotic of the first term on the
right-hand side of (\ref{main}) is $O_P(\frac{\ln L}{\sqrt{m_n}}).$
Let us prove the following
\begin{eqnarray*}\frac{(1 - \widehat{t}_k)^2}{(1 - t_k)^2} - 1 &=& o_P(1).\end{eqnarray*}
Using (\ref{hatminus}),
(\ref{t}) and (\ref{substr}), we obtain
\begin{eqnarray}\label{tt}&&\frac{(1 -
\widehat{t}_k)^2}{(1 - t_k)^2} - 1 = \frac{(\widehat{\theta}_{n}
e^{-\widehat{\theta}_{n} E^{(L)}_{l-k, l}})^2}{(\theta e^{-\theta
E^{(L)}_{l-k, l}})^2} - 1= (\widehat{\theta}_{n}
e^{-\widehat{\theta}_{n} E^{(L)}_{l-k, l}} - \theta e^{-\theta E^{(L)}_{l-k, l}})\cdot\nonumber
\\
&\cdot &(\widehat{\theta}_{n} e^{-\widehat{\theta}_{n} E^{(L)}_{l-k, l}}
- \theta e^{-\theta E^{(L)}_{l-k, l}} + 2\theta e^{-\theta E^{(L)}_{l-k, l}}
)/(\theta e^{-\theta E^{(L)}_{l-k, l}})^2 \nonumber
\\
&=&\frac{O_P\left(\frac{k \ln
L}{L\sqrt{m_n}}\right)\left(O_P\left(\frac{k \ln
L}{L\sqrt{m_n}}\right) + O_P\left(\frac{k}{L}\right)\right)}{
O_P\left(\frac{k^2}{L^2}\right)} = O_P\left(\frac{\ln
L}{\sqrt{m_n}}\right) = o_P(1).\end{eqnarray}
It remains to show that
$\omega^2_k(\widehat{\theta}_{n}) - \frac{1}{12k} =
O_P(1).$ It follows from (\ref{tt}) that $\frac{(1 -
\widehat{t}_k)^2}{(1 - t_k)^2} = O_P(1).$ Dividing the expression
(\ref{23}) by $\frac{(1 - \widehat{t}_k)^2}{(1 - t_k)^2},$ we obtain
again the expression, that is equal to $o_P(1)$. Namely, we get 
\begin{eqnarray}
&&
\left(\omega^2_k(\widehat{\theta}_{n}) -
\frac{1}{12k}\right) - \left(\omega^2_k(\theta) -
\frac{1}{12k}\right)\frac{(1 - t_k)^2}{(1 - \widehat{t}_k)^2}=\left(\omega^2_k(\widehat{\theta}_{n}) -
\frac{1}{12k}\right) \nonumber
\\
&-&
 \left(\omega^2_k(\theta) -
\frac{1}{12k}\right)\left(\frac{(1 - t_k)^2}{(1 -
\widehat{t}_k)^2}-1\right) - \left(\omega^2_k(\theta) -
\frac{1}{12k}\right) = o_P(1).\label{last}\end{eqnarray}
Using (\ref{tt}), we obtain that the second term on the
left-hand side of (\ref{last}) is $o_P(1),$ hence
\[\omega^2_k(\widehat{\theta}_{n}) - \frac{1}{12k} = O_P\left(\omega^2_k(\theta) -
\frac{1}{12k}\right) = O_P(1)\] holds. Therefore, the first term in
(\ref{main}) is $o_P(1)$.
\\
It remains to prove, that the first term in (\ref{9}) is $o_P(1).$ We have
\begin{eqnarray*}&&\widetilde{\omega}^2_L(\widehat{\theta}_{n}) - \omega^2_k(\widehat{\theta}_{n}) = \sum_{i=0}^{k-1}\left(\frac{i+0.5}{k} - y_i
\right)^2
-\sum_{i=0}^{k-1}\left(\frac{i+0.5}{k} - e_i
\right)^2
\\
&=& \sum_{i=0}^{k-1} \left(2\frac{i+0.5}{k} - e_i - y_i\right)(e_i - y_i),\end{eqnarray*}
where $e_i=e_i(\widehat{\theta}_{n})  = \exp\left(-\widehat{\theta}_{n}(E^{(L)}_{l-i,l}-E^{(L)}_{l-k, l})\right)$ and $y_i=y_i(\widehat{\theta}_{n})=$\\$= \exp\left(-\widehat{\theta}_{n}(Y_{L-i, L}-Y_{L-k, L})\right),$ $i \in\{0, \ldots, k-1\}.$ It easily follows from (\ref{newcondition}), (\ref{asymp}) and (\ref{10a}), that for all $i,$ $0\leq i\leq k-1,$
\begin{eqnarray*}\!\!e_i - y_i & =& \!\exp\!\left(\!-\widehat{\theta}_{n}\!(E^{(L)}_{l-i,l}-E^{(L)}_{l-k, l})\right)\! - \! \exp\!\left(\!-\widehat{\theta}_{n}\!(E^{(L)}_{l-i,l}-E^{(L)}_{l-k, l}+ o_P(1/\sqrt{k}))\right)
\\
&=& o_P(1/\sqrt{k}).\end{eqnarray*}
Further, from the latter, (\ref{asymp}), (\ref{10a}) and (\ref{renyi}) we obtain for all $i,$ $0\leq i\leq k-1,$
\[2\frac{i+0.5}{k} - e_i - y_i = o_P(1/\sqrt{k}).\] Thus, we derive
\begin{eqnarray}\label{6}&&\widetilde{\omega}^2_L(\widehat{\theta}_{n}) - \omega^2_k(\widehat{\theta}_{n}) = \sum_{i=0}^{k-1} o_P(1/k) = o_P(1),\end{eqnarray} the required result.

\subsection{Proof of Theorem \ref{TR4}}\label{S-5.3}
Note that formula (\ref{6}) in the proof of Theorem \ref{TR2}
\begin{equation}\widetilde{\omega}^2_L(\theta) - \omega^2_k(\theta) = o_P(1)\label{omega_small_difference}\end{equation}
is valid after the replacement $\hat{\theta}_{n}$ by $\theta$.
Simplifying (\ref{omegasquared2}), we obtain
\[\widetilde{\omega}^2_L(\theta) = \sum_{i=0}^{k-1} \left(\exp\left(-\theta(Y_{L-i, L} - Y_{L-k,L})\right)-\frac{i+0.5}{k}\right)^2 + \frac{1}{12k}.\]
Let us consider the difference $\widetilde{\omega}^2_L(\hat{\theta}_{n}) - \widetilde{\omega}^2_L(\theta)$. 
We have
\[
\widetilde{\omega}^2_L(\hat{\theta}_{n})- \widetilde{\omega}^2_L(\theta) = \sum_{i=0}^{k-1} y_i (d_i - 1) \left(y_i (d_i + 1) - 2\frac{i+0.5}{k}\right),\]
where $y_i =y_i(\theta)$ is taken as in the proof of Theorem \ref{TR2} 
and $d_i = \exp(-(\hat{\theta}_{n} - \theta) (Y_{L-i, L} - Y_{L-k,L}))$. 
It follows from (\ref{newcondition}) and (\ref{renyi}), that
\[y_i - \frac{i+0.5}{k} = y_i + o_p\left(\frac{1}{\sqrt{k}}\right) = O_P\left(\frac{1}{\sqrt{k}}\right)\] and
\[y_id_i - \frac{i+0.5}{k} = y_i(d_i - 1) + O_P\left(\frac{1}{\sqrt{k}}\right) = \left(\frac{i}{k} + O_P\left(\frac{1}{\sqrt{k}}\right) \right)(d_i-1) + O_P\left(\frac{1}{\sqrt{k}}\right).\]
Thus, the latter two equations imply
\begin{equation}
\widetilde{\omega}^2_{L_s}(\hat{\theta}_{n_s}) - \widetilde{\omega}^2_{L_s}(\theta) = \sum_{i=0}^{k_{n_s}-1}\frac{i^2}{k_{n_s}^2} (d_i-1)^2 \left(1 + O_P\left(\frac{1}{\sqrt{k_{n_s}}}\right)\right)
\label{omega_sum}\end{equation} under the condition that $\sqrt{k_{n_s}}|d_i-1| \to \infty$ as $n\to\infty$ holds.
Indeed, in terms of the subsequence $\{k_{n_s}\}_{s\ge 1}$, from (\ref{newcondition}) and (\ref{renyi}) it follows
\begin{eqnarray*} d_i - 1 &=& \exp\{-(\hat{\theta}_{n_s} - \theta) (Y_{L_s-i, L_s} - Y_{L_s-k_{n_s},L_s})\} - 1 = \frac{\ln (k_{n_s}/i)}{\theta}O_P(\hat{\theta}_{n_s} - \theta) \\
&=& \ln (k_{n_s}/i)\Omega_P(k_{n_s}^{-\alpha}) \end{eqnarray*}
for $\alpha>0$ and $d_i - 1 = \Omega_P(1)$ for $\alpha=0$. 
Here, $\xi_n := \Omega_P(\eta_n)$ means that $|\xi_n/\eta_n|\stackrel{P}{\rightarrow}\infty$ holds as $n\to\infty$.
 From the latter, (\ref{omega_small_difference}), (\ref{omega_sum}) and Theorem 2 we finally obtain
\begin{eqnarray*}\label{11}\widetilde{\omega}^2_{L_s}(\hat{\theta}_{n_s}) - \widetilde{\omega}^2_{L_s}(\theta)&=& k_{n_s}  O_P(\hat{\theta}_{n_s} - \theta)^2
\end{eqnarray*}
and
\begin{eqnarray}\label{12}\widetilde{\omega}^2_{L_s}(\hat{\theta}_{n_s}) &=& \left(\widetilde{\omega}^2_{L_s}(\hat{\theta}_{n_s}) - \widetilde{\omega}^2_{L_s}(\theta) \right) + (\widetilde{\omega}^2_{L_s}(\theta) - \omega^2_{k_{n_s}}(\theta)) + \omega^2_{k_{n_s}}(\theta)\nonumber
\\
&=& \Omega_P\left(k_{n_s}^{1-2\alpha}\right) + o_P(1) + O_P(1) = \Omega_P\left(k_{n_s}^{1-2\alpha}\right),\nonumber\end{eqnarray} the required result.

\section{Acknowledgements}
The authors were partly supported by the Russian Foundation for Basic
Research (grant \mbox{No.\,19-01-00090}).

\end{document}